\documentclass[preprint,3p,12pt]{elsarticle}

\usepackage{bm}
\usepackage{amsmath,amssymb}
\usepackage{IEEEtrantools}

\usepackage{xcolor}
\usepackage{verbatim}
\usepackage{epsfig}
\usepackage{graphicx}
\usepackage{makeidx}
\usepackage{multicol,color}
\usepackage{subfigure}
\usepackage{algorithm}
\usepackage{algorithmic}

\definecolor{ForestGreen}{RGB}{34,139,34}

\journal{Arxiv}

\begin{document}

\begin{frontmatter}
\begin{center}

\title{Adaptive model reduction of high-order solutions of compressible flows via optimal transport}
\end{center}

            
\author[inst1]{R. Loek Van Heyningen}
\author[inst1]{Ngoc Cuong Nguyen}
\author[inst2]{Patrick Blonigan}
\author[inst1]{Jaime Peraire}

\affiliation[inst1]{organization={Center for Computational Engineering, Department of Aeronautics and Astronautics, Massachusetts Institute of Technology},
            addressline={77 Massachusetts
Avenue}, 
            city={Cambridge},
            state={MA},
            postcode={02139}, 
            country={USA}}

\affiliation[inst2]{organization={Sandia National Laboratories},
            city={Livermore},
            postcode={94450}, 
            state={CA},
            country={USA}}
\begin{abstract}
The solution of conservation laws with parametrized shock waves presents challenges for both high-order numerical methods and model reduction techniques. We introduce an $r$-adaptivity scheme based on optimal transport and apply it to develop reduced order models for compressible flows. The optimal transport theory allows us to compute high-order $r$-adaptive meshes from a starting reference mesh by solving the Monge–Ampère equation. A high-order discretization of the conservation laws enables high-order solutions to be computed on the resulting $r$-adaptive meshes. Furthermore, the Monge-Ampère solutions contain mappings that are used to reduce the spatial locality of the resulting solutions and make them more amenable to model reduction. We use a non-intrusive model reduction method to construct reduced order models of both the mesh and the solution. The procedure is demonstrated on three supersonic and hypersonic test cases, with the hybridizable discontinuous Galerkin method being used as the full order model.
\end{abstract}

\end{frontmatter}
\section{Introduction}
Numerical simulations of high-speed flow simulations have a key role to play in the development of supersonic and hypersonic technologies, for which the setup and monitoring of experiments can be costly and challenging. 
Many-query workflows such as design optimization or uncertainty quantification may require the repeated solution of partial differential equations (PDEs) with sharp features and shocks at different sets of operating parameters.
When parametric uncertainty is considered for external hypersonic flow with strong shocks, it is often found that uncertainty in the free-stream Mach number has an outsized impact on the uncertainty in the quantities of interest, especially wall heat flux \cite{cortesi2020forward, tryoen2014bayesian, ray2023assessment}. 
Therefore, it is essential that the solution of a parametric system of equations can accurately capture dynamics with moving shocks, and that any surrogate models used to make these workflows tractable can also handle parametric shock movement. 
Unfortunately, the parametric variation of shocks poses challenges for  classical numerical methods and model reduction techniques. 

Low-order finite volume (FV) methods are still the most common choice for hypersonic simulations. 
However, it is well-known that for external flow problems with strong bow shocks, the quality of the computational mesh has a large impact on the solution when using most FV methods, with poorly designed but seemingly reasonable meshes resulting in non-physical results such as the carbuncle phenomenon \cite{candler2007unstructured, maccormack2013carbuncle, vila2023benchmarking}.
Mesh adaptation can be used to avoid such issues, with automatic mesh refinement (AMR), anisotropically refined grids, and $r$-adaptive shock tracking methods having shown success for high-speed flow. 
Of particular note is the grid tailoring procedure of \cite{saunders2007approach}, which iteratively aligns grid lines with a strong bow shock and is implemented in state-of-the-art hypersonics codes such as US3D \cite{candler2015development} and the Sandia Parallel Aerodynamics Reentry Code (SPARC) \cite{howard2017towards}, among others. This procedure retains a fixed mesh topology and is thus an $r$-adaptive method. 

An alternative to standard FV methods are high-order numerical methods, which have the potential to reach high levels of accuracy needed for high-fidelity simulations of complex phenomena. 
Discontinuous Galerkin (DG) methods are an attractive approach for CFD, due to their geometric flexibility, suitability for modern high-performance computer architectures, and amenability to adaptivity. 
Limiting the use of DG for high-speed flow is the fact that the presence of shocks can be particularly challenging for high-order methods, where the Gibbs phenomena can severely impact accuracy and robustness. The literature for shock capturing for DG methods is vast, with approaches including flux and slope limiting \cite{burbeau2001problem, cockburn1989tvb, krivodonova2004shock, krivodonova2007limiters, lv2015entropy}, locally low-order fluxes \cite{huerta2012simple, sonntag2017efficient, persson2022discontinuous}, and artificial viscosity (AV) methods \cite{persson2006sub, persson2013shock, fernandez2018physics, barter2010shock, ching2019shock, bai2022continuous}. 

Manipulation of the computational mesh can also greatly improve resolution of sharp features when using high-order methods. 
Numerous works have coupled goal-oriented mesh adaptation with high-order methods for high-speed flow problems \cite{barter2010shock, yano2011importance, may2021hybridized, sabo2022investigation, bai2022continuous}. 
The relative benefits of $r$-adaptive shock-aligned meshes are arguably greater for high-order discontinuous methods; curved shocks can be tracked and high-order convergence is able to be regained when the shock alignment is performed exactly. 
By formulating the solution as an optimization problem over solution and mesh degrees of freedom, the High-Order Implicit Shock Tracking (HOIST) \cite{zahr2020implicit} and Moving DG with Interface Condition Enforcement (MDG-ICE) \cite{corrigan2019moving} methods can capture highly accurate solutions on coarse meshes. 
Both approaches have been applied to challenging hypersonic problems and returned accurate solutions on remarkably coarse meshes, including for viscous and reacting problems \cite{ching2023moving, huang2022robust, zahr2021high}. %

While accurate and robust evaluations of the full-order model (FOM) at different Mach numbers are important, so too is the development of accurate surrogate models. 
In parametric CFD, projection-based reduced order models (ROMs) have shown promise in delivering surrogate models that are efficient, can be coupled with rigorous error estimates, and exhibit low error even in high-dimensional parameter spaces and in regions of extrapolation. 
This context is where $r$-adaptivity can be especially impactful. If each snapshot of the FOM has a different grid topology, the definition of inner products between snapshots becomes non-trivial. Although some works \cite{sleeman2022goal, little2023nonlinear} have addressed this for AMR mesh adaptation, generalizing these approaches to anisotropically adapted meshes remains a challenge. 
A more fundamental issue is that variation in the free-stream Mach number will usually cause shocks to move location; this means each FOM snapshot will have different locally coherent structures and the parametric solution manifold will suffer from a large Kolmogorov n-width. 
If the parametric variation is such that the shock does not move very much, classical linear model reduction techniques such as the proper orthogonal decomposition (POD) method can give acceptable answers \cite{blonigan2021model}.
As the variation increases, these methods will not be able to efficiently characterize the solution variation without a large amount of training data. 
Some degree of nonlinearity is usually required, such as with adaptively constructed or refined linear ROMs \cite{carlberg2015adaptive}, a combination of FOM and ROM information in different parts of the domain \cite{legresley2003dynamic, huang2023predictive, peherstorfer2020model, zucatti2023adaptive}, or a choice of a nonlinear basis. Variations include linear bases augmented with nonlinear closure terms \cite{barnett2022quadratic, barnett2023neural}, fully nonlinear bases \cite{lee2020model, kim2022fast, romor2023non}, or bases that consist of combining a parameter-dependent mapping and a linear basis. 

Examples of the latter for steady high-speed flow include the TSMOR method of \cite{nair2019transported}, wherein mappings are derived based on wave speed arguments. The ROM-IFT \cite{mirhoseini2021model} method minimizes over a reduced solution space and admissible domain mappings to align features, using similar arguments from the $r$-adaptive HOIST method. End-to-end efficiency gains are improved with greedy sampling of the solution and mapping space along with empirical quadrature-based hyperreduction \cite{mirhoseini2023accelerated}. The registration methods of \cite{taddei2020registration} also compose a mapping and a linear basis. FOM snapshots are evaluated, a reference state is chosen, and a regularized optimal transport problem is solved to align the snapshots to some reference state. Recent improvements of the method include the incorporation of goal-oriented metric-based adaptation, empirical quadrature hyperreduction, and greedy multifidelity sampling methods \cite{ferrero2022registration, barral2023registration}. The grid-tailored ROMs of \cite{ching2022gridtailor} use the outputs of a grid-tailored FV method simulation to construct a non-intrusive ROM for grid deformations and an intrusive ROM for the solution field.

Many of the aforementioned works emphasize the connection between $r$-adaptivity methods for the FOM and the learning of mappings that can be composed with a linear basis to form accurate ROMs for problems with moving features. In a similar vein, we make use of a recently developed high-order method for $r$-adaptivity based on optimal transport to aid in snapshot collection and data compression for model reduction of problems with parametrically varying shock structures. 
We propose that each FOM snapshot will include a solution field and mesh deformation. 
The mesh deformations are learned via solutions of the Monge-Ampère equation, an $r$-adaptation method derived from the application of optimal transport formulations to an equidistribution principal for mesh adaptation. 
Solutions of the Monge-Ampère equation have been used extensively for $r$-adaptivity \cite{budd2009adaptivity, budd2015geometry, budd2018scaling, mcrae2018optimal, dipietro2019adaptive, ramani2023fast}. The equation can be challenging to solve with standard numerical methods, but the resulting meshes tend to avoid entanglements and smoothly transition between regions of refinement.

This is an elliptic equation that naturally handles curved boundaries and is itself amenable to model reduction \cite{hou2022RBMonge}.

A number of authors have used optimal transport ideas to remove or reduce the convective nature of parametrized problems \cite{iollo2014advection, cagniart2019model, torregrosa2022surrogate, blickhan2023registration}. 
The most akin to ours are the registration methods of \cite{taddei2020registration}. 
Registration methods form a regularized optimization statement of an optimal transport problem which explicitly maps snapshots onto a single reference configuration. We solve the Monge-Ampère equation to map uniform mesh densities to relevant features; while we expect our method will do less to mitigate the Kolmogorov barrier than most optimal transport formulations since alignment is not explicitly targeted, it is a practical formulation that reuses our mesh adaptation procedure. 

These mesh mappings can be applied to the FOM snapshots after training or they can be created as part of the FOM solution. We demonstrate the latter approach here. The mappings can be used as a mesh refinement strategy and an enabler for model reduction. This is similar to the examples of ROM-IFT coupled with HOIST, but is most similar to the grid-tailored ROMs of \cite{ching2022gridtailor}. 
Grid tailoring is a shock alignment method while ours is not. Grid tailoring more explicitly aims for shocks to be fixed at a specific mesh degree of freedom.
We sacrifice exact alignment for flexibility over grid tailoring, which is restricted to bow shocks and requires the ability to compute lines from the wall to the inflow boundary. 
The Monge-Ampère $r$-adaptivity can be used with structured and unstructured grids, tetrahedral or hexahedral elements, and can refine along various curved or interacting features, as long as they can be tracked by a scalar value monitor function. 

The use of a high-order method for the flow solution allows for the PDE to be iteratively solved on anisotropically adapted meshes generated by the solution of the Monge-Ampère equation. 

Some simplifications are made for this first study of the effectiveness of mappings constructed with Monge-Ampère mesh adaptation procedures in aiding model reduction. First, we only consider parametric variation in one parameter, the Mach number. 
This can still challenge linear basis ROMs, but the effects of the Kolmogorov barrier are in general less pronounced in this parameter regime and are more easily able to be overcome with sampling. 
Also, we only consider nonintrusive interpolation-based ROMs and test cases without extrapolation. 
With the small dimensional parameter spaces considered here,we expect nonintrusive interpolation ROMs to perform sufficiently well, while still providing insight into the effectiveness of this approach. 
It is worth mentioning that for steady high-speed flow problems, interpolation-based ROMs can be surprisingly effective for cases without extrapolation \cite{blonigan2021model, ching2022gridtailor}.

The paper is organized as follows: Section \ref{sec:methods} details the full order model discretization and proposed $r$-adaptivity scheme. 
The nonintrusive ROM formulations are presented in Section \ref{sec:ROM}, while Section \ref{sec:results} shows results on three supersonic and hypersonic flow problems with varying degrees of shock strength and complexity. Conclusions and future work are described in Section \ref{sec:conclusions}.

\section{Methods}
\label{sec:methods}
\subsection{Steady parametrized conservation laws} 
\label{sec:fom}
We consider a system of steady conservation laws on a domain $\Omega \in \mathbb{R}^d$
\begin{equation}   
    \nabla \cdot \bm F(\bm u(\bm x; \bm \mu); \bm \mu) = 0, \; \; \forall \bm x \in \Omega
    \label{eq:Fu}
\end{equation}
where $\bm u \in \mathbb{R}^m$ is the solution with $m$ components, $\bm \mu \in \mathcal{D} \subset \mathbb{R}^{q}$ is a vector of parameters, and $\bm F \in \mathbb{R}^{m \times d}$ is a physical flux function. The domain $\Omega$ is partitioned into a collection of disjoint elements $\mathcal{T}_h$ and parametrically varying boundary conditions are specified on the boundary $\Gamma \in \partial \Omega$. 

For this work, we exclusively solve the compressible Euler equations in conservative form,  
\begin{equation}
    \bm u = \begin{pmatrix}
        \rho \\
        \rho v_i \\
        \rho E \\
    \end{pmatrix}, \; \; 
    \bm F_i = \begin{pmatrix}
        \rho v_i \\ 
        \rho v_i v_j + \delta_{ij} p \\ 
        \rho v_i H
    \end{pmatrix}
    \label{eq:euler}
\end{equation}
with density $\rho$, velocity $\bm v$, total energy $E$, total specific enthalpy $H = E + p/\rho$ and pressure $p$ given by the ideal gas law.  Let $\Gamma_{\rm wall} \subset \partial \Omega$ be the wall boundary. The boundary condition at the wall boundary $\Gamma_{\rm wall}$ is $\bm v \cdot \bm n = 0$, where $\bm v$ is the velocity field and $\bm n$ is the unit normal vector outward the boundary. Because supersonic and hypersonic flows are considered in this paper, supersonic inflow and outflow conditions are imposed on the inflow and outflow boundaries, respectively.

\subsection{Adaptive viscosity regularization}
\label{sec:av}
For discretizing and solving \eqref{eq:Fu}, we use the adaptive viscosity regularization approach of \cite{nguyen2023adaptive}. For clarity, we drop the parametric dependence for this section. While not the focus of this work, we give a brief overview for context, as this method allows for the rapid and reliable solution of the shock-dominated examples below. 
Instead of solving \eqref{eq:Fu} directly, we augment it with an artificial viscosity regularization term $\bm G$ and an equation to smooth out the resulting artificial viscosity, in line with PDE-based artificial viscosity approaches \cite{barter2010shock}
\begin{subequations}
\label{eq:viscsystem}
\begin{alignat}{2}
\nabla \cdot \bm F(\bm u) - \lambda_1 \nabla \cdot \bm G(\bm u,   \nabla \bm u,  \eta) = 0  \quad \mbox{in }\Omega,\label{eq:1}\\
\eta - \lambda_2^2 \nabla \cdot \left( \ell^2 \nabla \eta \right) - s(\bm u, \nabla \bm u) = 0 \quad \mbox{in }\Omega \label{eq2}
\end{alignat}
\end{subequations}
where the Helmholtz equation in \eqref{eq2} is assigned homogeneous boundary conditions  $\eta=0$ at the wall boundary and $\nabla \eta \cdot \bm n = 0$ at the remaining  boundary. The term $s(\bm u, \nabla \bm u)$ is defined as 
\begin{equation}
    s(\bm u, \nabla \bm u) = g_{\text{clip}}\left(\mathcal{S}(\bm u, \nabla \bm u) ; s_{\text{min}}, s_{\text{max}} \right)
    \label{eq:clip}
\end{equation}
where $g_{\text{clip}}(\cdot; s_{\text{min}}, s_{\text{max}})$ is a smooth function that limits the first argument between $s_{\text{min}}$ and $s_{\text{max}}$. 
The term $\mathcal{S}(\bm u, \nabla \bm u)$
is a shock sensor, here chosen to be the negative divergence of the velocity \cite{moro2016dilation}
\begin{equation}
    \mathcal{S}(\bm u, \nabla \bm u) = -\nabla \cdot \bm v. 
\end{equation}
The limiting parameters are chosen as in \cite{nguyen2023adaptive}, with $s_{\text{min}} = 0$ to avoid negative viscosity and $s_{\text{max}}$ is set iteratively to $0.5\|\mathcal{S}\|_{\infty}$.  
Artificial dissipation is added with $\bm G(\bm u, \nabla \bm u)$, defined as
\begin{equation}
\bm G(\bm u,   \nabla \bm u,  \eta) = \mu(\eta) \nabla \bm u \ . 
\end{equation} 
The parameters $\lambda_1, \lambda_2$ which control the magnitude and width of the artificial viscosity respectively, are iteratively ramped down until regularity or physical constraints are violated. In other words, the method aims to minimize the amount of artificial viscosity while maintaining smooth and accurate solutions. The method acts as a nonlinear solver and procedure for controlling artificial viscosity amounts \cite{nguyen2023adaptive}.

Equation \eqref{eq:1} is discretized with a hybridized discontinuous Galerkin (HDG) method while \eqref{eq2} is discretized with continuous Galerkin (CG).
HDG is chosen for the flow system due to its stability for convective problems and reduced globally coupled degrees of freedom compared to standard DG methods.  See \cite{nguyen2011adaptive,Nguyen2012,Fernandez2018a,Moro2011a,Peraire2010,Vila-Perez2021,Fidkowski2016,Fernandez2017a} for more details on the HDG discretization of the Euler equations, including the precise weak form statements, static condensation solution procedures, and boundary conditions.  CG is used for the Helmholtz equation \eqref{eq2} so that the artificial viscosity field is continuous, which tends to result in more robust artificial viscosity methods \cite{bai2022continuous, ching2019shock, barter2010shock}.

As emphasized in \cite{nguyen2023adaptive}, the amount of artificial viscosity depends on the mesh. In the extreme case of perfect alignment of shocks with element faces for discontinuous schemes with a suitable numerical Riemann solver, no viscosity is needed \cite{corrigan2019moving, zahr2020implicit}.
In \cite{nguyen2023adaptive}, a solution on an initial mesh was used to locate shock structures and a new shock-aligned mesh was generated. 
While capable of producing high-quality solutions, when solving parametric systems of PDEs, this procedure would generate meshes with different topologies for each snapshot.

\subsection{Mesh adaptation} 
\label{sec:mesh}

\subsubsection{Monge-Ampère equation for mesh adaptation}
\label{sec:monge}

Mesh adaptation is accomplished via high-order solutions of the Monge-Ampère equation of optimal transport, based on the method introduced in \cite{nguyen2023monge}.  
Mesh adaptation using the Monge-Ampère is usually oriented around an equidistribution principle. Suppose we are given a target $\bm \mu$-dependent density function $\rho'(\bm x; \bm \mu)$ and a source $\bm \mu$-dependent density function $\rho(\bm x; \bm \mu)$ defined on $\Omega$.  The optimal transport problem seeks to find a mapping $\phi(\bm x)$ such that
\begin{equation}
    \inf_{\bm \phi \in \mathcal{M}_{\bm \mu}} \int_{\Omega} \| \bm x - \bm \phi(\bm x; \bm \mu) \|^2 \rho (\bm x; \bm \mu) d \bm x
    \label{eq:MAstatement}
\end{equation}
with
\begin{equation}
    \mathcal{M}_{\bm \mu} = \left\{ \bm \phi \, : \, \Omega \rightarrow \Omega, \; \rho'(\bm \phi(\bm x; \bm \mu); \bm \mu) \det\left( \nabla \bm \phi(\bm x; \bm \mu))\right) = \rho(\bm x; \bm \mu), \; \forall \bm x \in \Omega \right\}
    \label{eq:MAset}
\end{equation}
being defined as the set of feasible mappings.
Under certain regularity conditions described by Brenier  \cite{brenier1991polar} and Caffarelli \cite{caffarelli1990interior}, the optimal map of \eqref{eq:MAstatement} is equal to the gradient of a scalar convex potential $w$. Inserting $\bm \phi(\bm x; \bm \mu) = \nabla w(\bm x; \bm \mu)$ into the feasible set of $\mathcal{M}_{\bm \mu}$ gives the parametrized Monge-Ampère equation
\begin{equation}
    \rho' \left( \nabla w(\bm x; \bm \mu); \bm \mu\right) \det \left( D^2 w(\bm x; \bm \mu)\right) = \rho(\bm x; \bm \mu).
\end{equation}

We require that that $\bm \phi$ conforms to the physical boundary $\partial \Omega$. If the physical boundary $\partial \Omega$ can be expressed as the equation $g(\bm x) = 0$ for  $\bm x \in \partial \Omega$, we impose the following nonlinear Neumann boundary condition: 
\begin{equation}
    g(\nabla w( \bm x; \bm \mu)) = 0, \; \forall \bm x \in \partial \Omega .
\end{equation}
Mesh adaptation by the equidistribution principle means that the source density is a uniform redistribution of the target density on the domain. Therefore, we have $\rho(\bm x; \bm \mu) = \theta(\bm \mu)$ for $\theta(\bm \mu) = \int_{\Omega} \rho'(\bm x; \bm \mu) d\bm x / \int_{\Omega} d \bm x$.  The  Monge-Ampère equation can then be written as 
\begin{subequations}
 \label{eq:MA}
\begin{align}
\det(D^2w(\bm x; \bm \mu)) &= f(\nabla w(\bm x; \bm \mu); \bm \mu), \; \forall \bm x \in \Omega, \\
g(\nabla w(\bm x; \bm \mu)) &= 0, \, \qquad 
\; \, \qquad \qquad\forall \bm x \in \partial \Omega, \\
\int_{\Omega} w(\bm x; \bm \mu) d \bm x &= 0,
\end{align}
\end{subequations} 
where $f(\nabla w(\bm x; \bm \mu); \bm \mu) = \theta(\bm \mu) / \rho'(\nabla w(\bm x; \bm \mu); \bm \mu)$. The Monge-Ampère equation (\ref{eq:MA}) is a nonlinear elliptic PDE with a nonlinear Neumann boundary condition. 

For two dimensions, the Monge-Ampère equation can also be written as
\begin{subequations}
\label{eq2DMA}
    \begin{align}
    \bm H - \nabla \bm q &= 0, \qquad\text{in } \Omega \\
    \bm q- \nabla w &= 0, \qquad \text{in } \Omega \\
    S(\bm H, \bm q; \bm \mu) - \nabla \cdot \bm q &= 0, \qquad \text{in } \Omega \\
    g(\bm q) &= 0, \qquad \text{on } \partial \Omega \label{eq:MAbc} \\
  \int_{\Omega} w(\bm x) d \bm x &= 0
\end{align}
\end{subequations}
with $S(\bm H, \bm q; \bm \mu) = \sqrt{H_{11}^2 + H_{22}^2 + H_{12}^2 + H_{21}^2 + 2f(\bm q; \bm \mu)}$. 
Equation \ref{eq2DMA} 
is discretized in space with an HDG method described in \cite{nguyen2023monge}.
The static condensation procedure of HDG means that the size of the global unknowns scales with the trace of the scalar variable, not of the gradient or Hessian terms. For smooth analytic choices of target density, it was found that the solution, gradient, and Hessian all converge with the same order of accuracy. 
Finally, the use of HDG allows for solutions on unstructured grids and curved boundaries, meaning that the starting reference mesh can be high-order and unstructured.

In the remainder of the paper, the subscript $h$ will indicate the numerical solution computed by using the HDG method. Hence, the mapping $\bm q_h$ is an approximation to the exact mapping $\bm \phi(\bm x)$, which is a discontinuous field. It  means that a single vertex will have two separate mappings that are not required to be equal.  In practice, we use a continuous approximation of $\bm q_h$ by averaging the duplicate degrees of freedom of $\bm q_h$, which shall be denoted as $\bm \phi_h$. A more subtle issue is the presence of corners inside the domain.
Other $r$-adaptivity procedures have found it useful to constrain corner nodes \cite{fortunato2016high, huang2022robust}. 
Since the mesh mapping is given by the \textit{gradient} of a scalar potential, it is not obvious how to constrain corner nodes. 
It was found that geometries with internal corners posed challenges for the HDG Monge-Ampère solver. 
This can be avoided by making sure that $g(\cdot) = 0$ is a \textit{global} description of the geometry. 
For example, suppose we have a boundary with intersecting lines $g_1(\bm x) = 0$ for $\bm x\in \Gamma_1$ and $g_2(\bm x) = 0$ for $\bm x \in \Gamma_2$.
When evaluating \eqref{eq:MAbc} on a face that lies on $\Gamma_1$, rather than always letting $g_1(\nabla u_h) = 0$, we switch to enforcing $g_2(\nabla u_h) = 0$ if $\nabla u_h \in \Gamma_2$. This approach, though allowing the solver to converge, lets boundary nodes transition from one boundary to another. This means corner nodes can get detached from corners, which can be manually corrected in the mapping. It also means that nodes that were once a wall boundary can move to an outflow boundary, for example. For our solver, this just means that the mapped mesh needs to be reinitialized and reassigned boundary conditions. 

\subsubsection{Adaptation Procedure}
\label{sec:adapt}

We still need to define the mesh density function $\rho'(\bm x; \bm \mu)$ in order to obtain a good $r$-adaptive mesh $\bm \phi_h(\bm x; \bm \mu)$. The optimal transport will drive mesh nodes to concentrate around high values of the mesh density function. Therefore, $\rho'(\bm x; \bm \mu)$ should be large in the region where refinement is needed and small elsewhere. It is possible to choose the artificial viscosity field as a target density so that the $r$-adaptive mesh $\bm \phi_h(\bm x; \bm \mu)$ is fine in the shock region. However, the dilatation sensor used can fail to detect other sharp features including contact discontinuities. Instead, we use the gradient of a scalar function $\xi_h(\bm x; \bm \mu)$ to define a resolution sensor: 
\begin{equation}
    s_h(\bm x; \bm \mu) = \sqrt{1 +  g_{\text{clip}}(\|\nabla \xi_h (\bm x; \bm \mu)\|_\Omega^2; s_{\text{min}}, s_{\text{max}})}
    \label{eq:sensor}
\end{equation}
where $s_{\text{min}} = 0 $ and $s_{\text{max}} = 0.5 \| \nabla \xi_h(\bm x; \bm \mu) \|_{\infty}$. 
The mesh density function is computed by solving the following Helmholtz equation
\begin{equation}
    \rho'(\bm x; \bm \mu) - \nabla \cdot \left( \ell^2 \nabla \rho'(\bm x; \bm \mu) \right) = s_{h}(\bm x; \bm \mu) \quad \mbox{in }\Omega,  \label{eq:sensorsmooth}
\end{equation}
with homogeneous Neumann boundary conditions on the whole domain. This is a feature-based adaptation indicator. Other indicators are possible, such as  those based on adjoint information or some combination of physics-based artificial viscosity sensors that can distinguish between shocks, large temperature gradients, and other sharp features.
The only requirement is that an indicator function be expressed as a scalar-valued positive function.

It remains to describe how to choose the scalar function $\xi_h(\bm x; \bm \mu)$. In this paper, we choose $\xi_h(\bm x; \bm \mu)$ as the approximate fluid density obtained by using the adaptive viscosity regularization method described in Subsection \ref{sec:av}. More specifically, let $\mathcal{T}_h^0$ denote the  high-order reference mesh on the physical domain $\Omega$. Note that the reference mesh $\mathcal{T}_h^0$ is independent of $\bm \mu$. For any given $\bm \mu \in \mathcal{D}$, we repeatedly solve the system (\ref{eq:viscsystem}) and reduce the values of $\lambda_1, \lambda_2$ by using the HDG/CG scheme with the reference mesh $\mathcal{T}_h^0$ until any of the smoothness and positivity constraints is violated. This yields $\bm u_h^0(\bm x; \bm \mu)$ as the numerical approximation to the exact solution $\bm u$ of the conservation laws (\ref{eq:Fu}). Then $\xi_h(\bm x; \bm \mu)$ is nothing but the first component of $\bm u_h^0(\bm x; \bm \mu)$.

The fact that the target mesh density $\rho'$ is the solution of the Helmholtz equation (\ref{eq:sensorsmooth}) makes the  Monge-Ampère equation more expensive to be solved numerically. Specifically, the fixed point scheme is used to solve the Monge-Ampère equation (\ref{eq2DMA}). At the $l$-th iteration of the fixed point scheme, we  evaluate $S(\bm H_h^{l-1}, \bm q_h^{l-1}; \bm \mu)$ which in turn requires us to evaluate $\rho'(\bm q_h^{l-1}; \bm \mu)$.  Because the mesh density function $\rho'$ is constructed from local polynomials of high degree on the reference mesh $\mathcal{T}_h^0$, we must interpolate it onto  $\bm q_h^{l-1}(\bm \mu)$. The interpolation of a scalar field from one high-order mesh to another high-order mesh can be computationally expensive. The fixed-point scheme reaches a convergence when $\|\bm q_h^{l}(\bm \mu) - \bm q_h^{l-1}(\bm \mu)\|_{\Omega}$ is less than a specified tolerance. At convergence, $\bm q_h^{l}(\bm \mu)$ is averaged at the duplicate degrees of freedom to yield the mapping $\bm \phi_h(\bm \mu)$, which defines the resulting $r$-adaptive mesh $\mathcal{T}_h^{\bm \mu}$. The mesh adaptation procedure is described in Algorithm \ref{alg:meshadapt}.

\setcounter{algorithm}{-1}
\begin{algorithm}
\caption{Monge-Ampère mesh adaptation for steady parametrized problem}
\begin{algorithmic}[1]
\REQUIRE{ The reference mesh $\mathcal{T}_h^0$ and parameters $\bm \mu \in \mathcal{D}$.}
\ENSURE{The mapping $\bm \phi_h(\bm \mu)$  that maps the reference mesh $\mathcal{T}_h^0$ to the adaptive mesh $\mathcal{T}_h^{\bm \mu}$, and the numerical solution $\bm u_h(\bm \mu)$ on $\mathcal{T}_h^{\bm \mu}$.}
\STATE{Solve for $\bm u_h^0(\bm \mu)$ on $\mathcal{T}_h^0$ using the adaptive viscosity regularization method.}
\STATE{Compute the target mesh density $\rho'(\bm x; \bm \mu)$ based on $\bm u_h^0$ by using \eqref{eq:sensor} and \eqref{eq:sensorsmooth}.}
\STATE{Solve the Monge-Ampére equation \eqref{eq2DMA} on $\mathcal{T}_h^0$ using the fixed-point HDG method.}
\STATE{Average $\bm q_h$ to obtain the mapping $\bm \phi_h(\bm \mu)$ and the associated $r$-adaptive mesh $\mathcal{T}_h^{\bm \mu}$.}
\STATE{Interpolate $\bm u_h^0(\bm \mu)$ onto $\mathcal{T}_h^{\bm \mu}$ and use it as the initial guess.}
\STATE{Solve for $\bm u_h(\bm \mu)$ on $\mathcal{T}_h^{\bm \mu}$ using the adaptive viscosity regularization method.}

\end{algorithmic}
\caption{Monge-Ampère mesh adaptation for steady parametrized problems.}
\label{alg:meshadapt}
\end{algorithm}

\begin{figure}[htb!]
\begin{center}
\subfigure[Input: reference mesh $\mathcal{T}_h^0.$]
{
\includegraphics[height=3.9cm,width=0.17\textwidth]{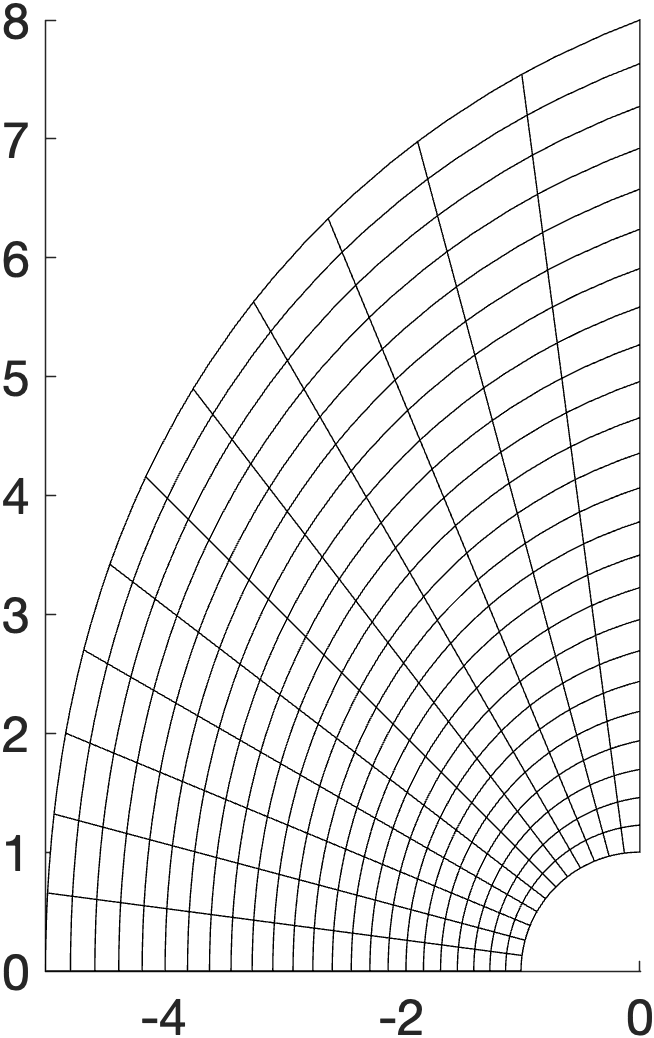}\label{fig:algmesh}} \qquad
\subfigure[Step 1: Compute solution $\bm u_h^0(\bm \mu)$ on $\mathcal{T}_h^0$]
{
\includegraphics[width=0.18\textwidth]{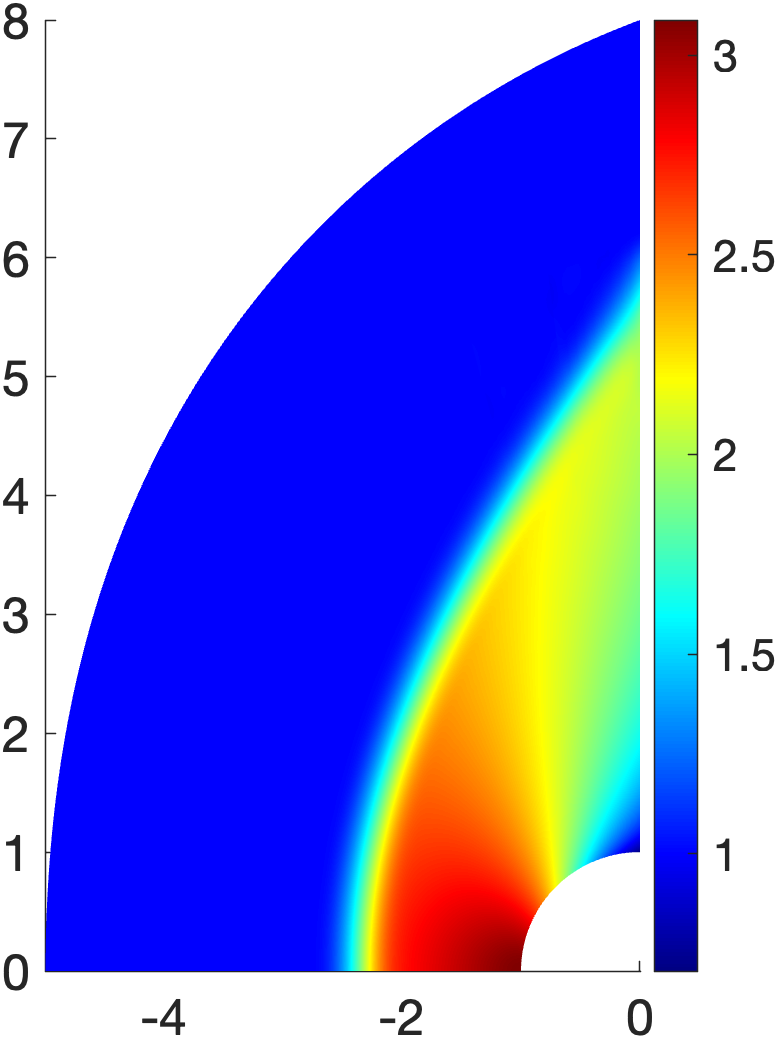}\label{fig:algsol0}} 
\qquad
\subfigure[Step 2: Compute target mesh density $\rho'(\bm x; \bm \mu).$]
{\includegraphics[width=0.18\textwidth]{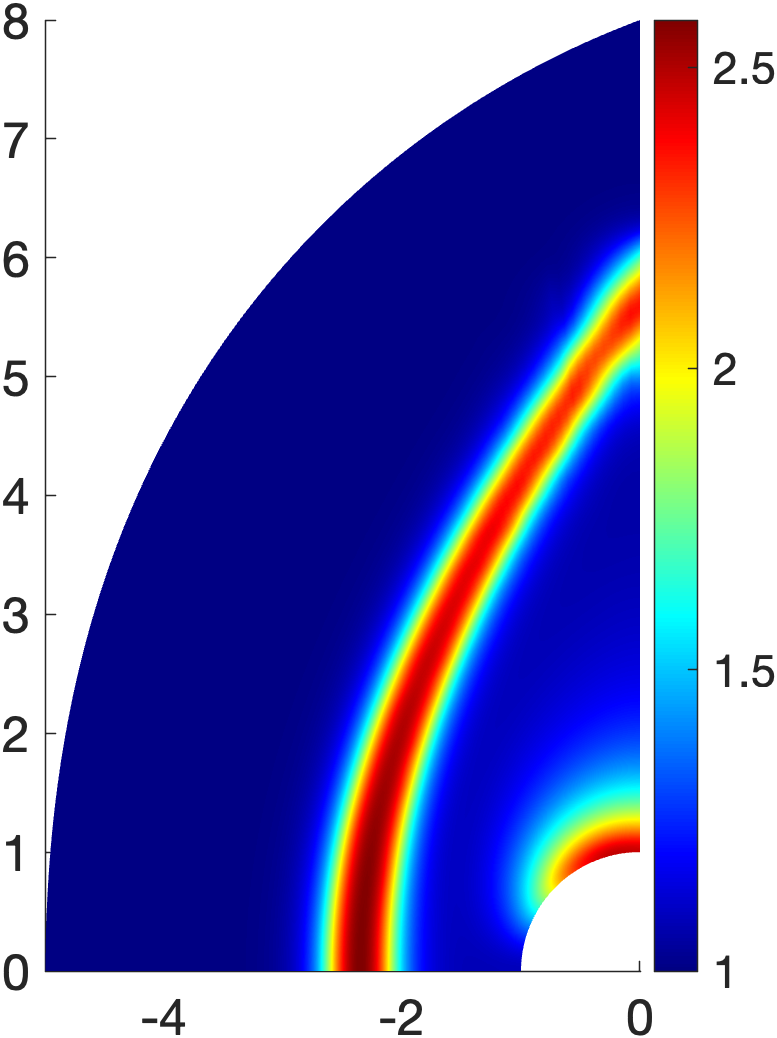}\label{fig:algsensor}} \\
\subfigure[Step 3: Solve Monge-Ampère equation to get $q_{h,x}$ (left) and $q_{h,y}$ (right).]
{\includegraphics[width=0.18\textwidth]{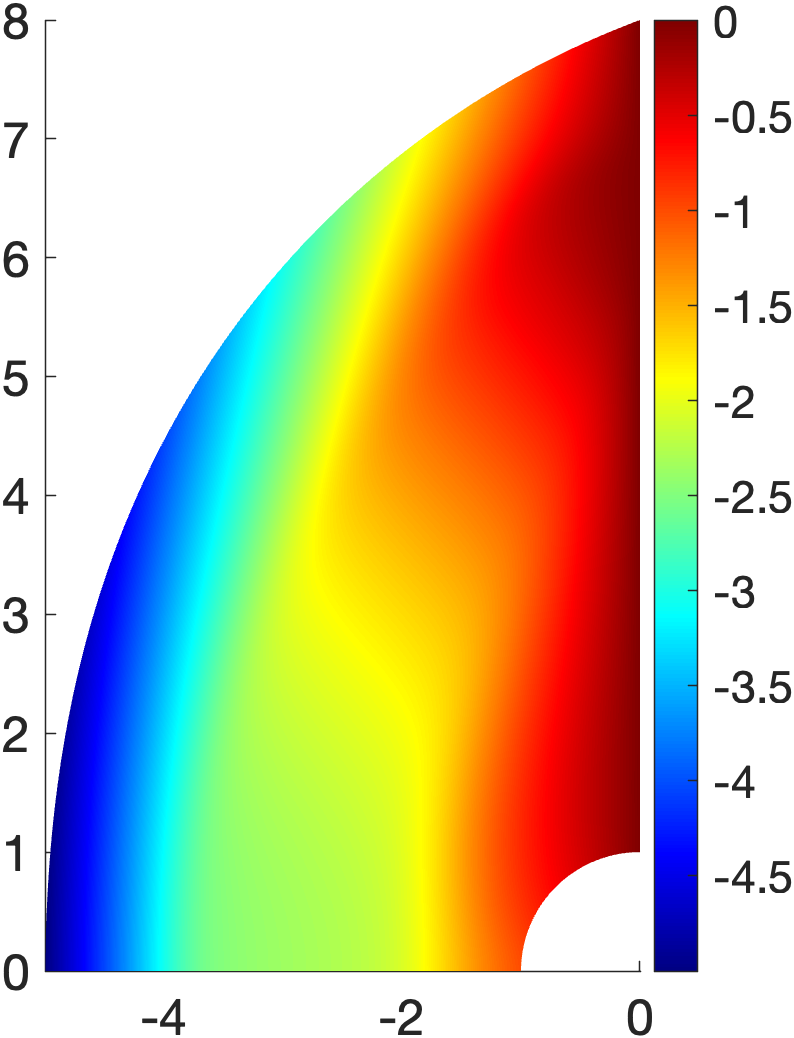}\label{fig:algq1}
\includegraphics[width=0.17\textwidth]{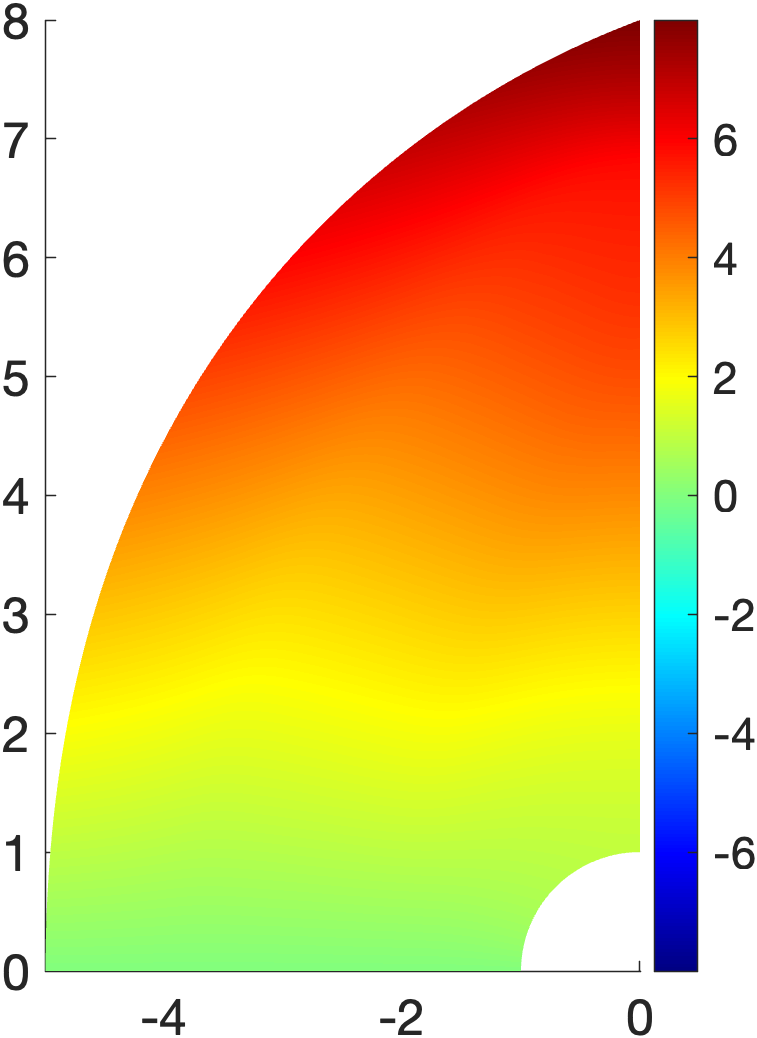}\label{fig:algq2}} \qquad
\subfigure[Step 4: Average $\bm q_h$ to obtain adapted mesh $\mathcal{T}_h^{\bm \mu}.$]
{\includegraphics[height=3.9cm,width=0.17\textwidth]{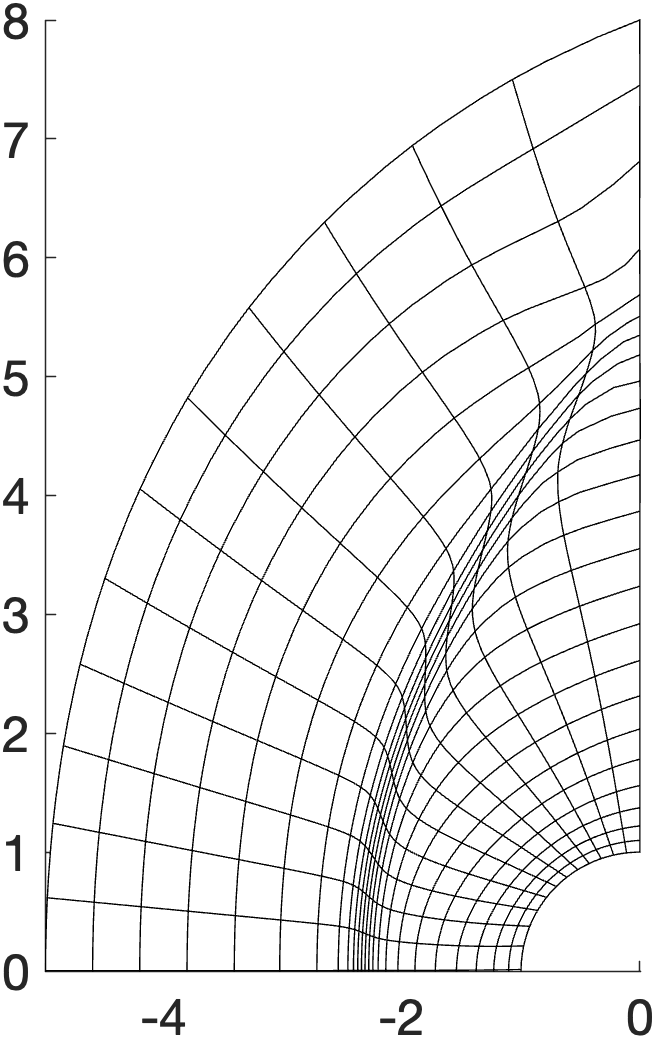}\label{fig:algmeshadapt}} 
\qquad
\subfigure[Steps 5-6: Solve for $\bm u_h(\bm \mu)$ on $\mathcal{T}_h^{\bm \mu}. $]
{\includegraphics[width=0.18\textwidth]{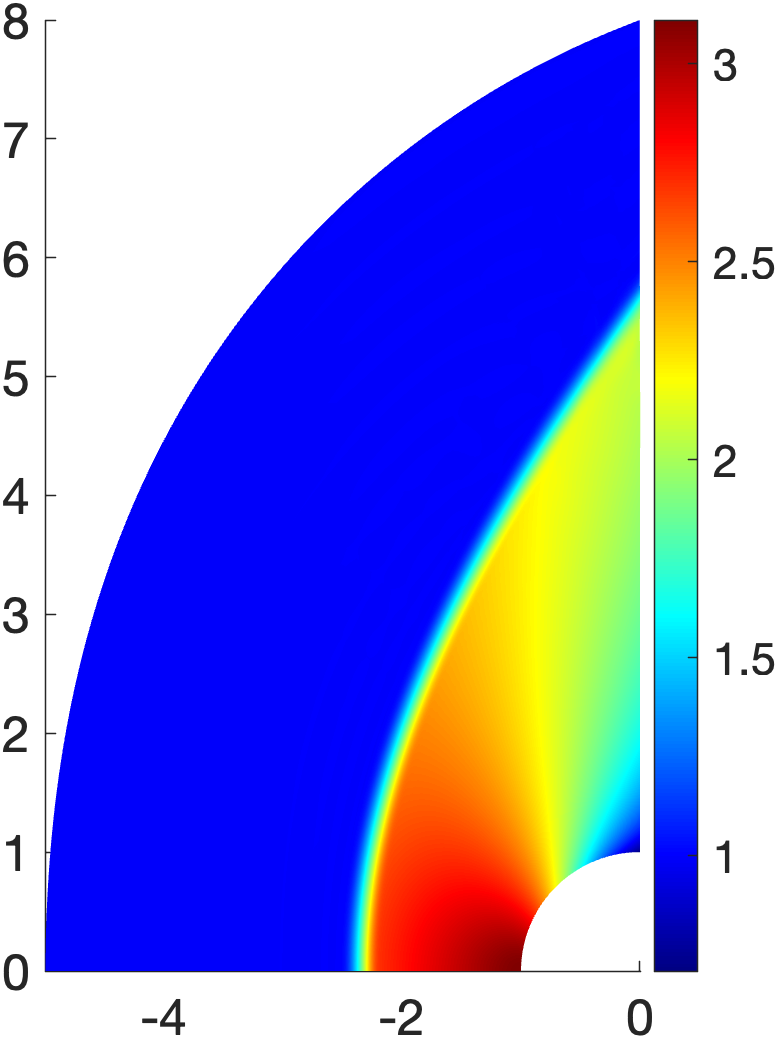}\label{fig:algsoladapt}} 
\caption{Visualization of steps of Algorithm \ref{alg:meshadapt}.}\label{fig:algdemo}
\end{center}
\end{figure}

To clarify each step of Algorithm \ref{alg:meshadapt}, we consider a sample flow field with a single bow shock in front of a cylindrical body. Representative inputs and outputs are shown in Figure \ref{fig:algdemo} for a free-stream Mach number of $\bm \mu = 2.0$.

\section{Reduced Order Modeling}
\label{sec:ROM}
\subsection{Model reduction on the reference domain}
\label{sec:ROMref}

In addition to allowing for mesh adaptation to be performed during the training of a reduced-order model, the optimal transport method also makes it convenient to map the numerical solution $u_h(\bm x; \bm \mu_i)$ from the $r$-adaptive mesh $\mathcal{T}_h^{\bm \mu}$ to the reference domain $\mathcal{T}_h^0$.  In $r$-adaptivity, mesh points are neither created nor destroyed, connectivity structure does not need to be modified.  We exploit this property of $r$-adaptivity to construct reduced order models. 

We perform Algorithm \ref{alg:meshadapt} for a given set of training parameters $\{\bm \mu_i\}_{i=1}^{n_{\text{train}}}$ to obtain snapshots of the numerical solution $\{\bm u_h(\bm x; \bm \mu_i)\}_{i=1}^{n_{\text{train}}}$ and the mapping, $\{\bm \phi_h(\bm x; \bm \mu_i)\}_{i=1}^{n_{\text{train}}}$. 
Note that  all of the mapping snapshots  $\{\bm \phi_h(\bm x; \bm \mu_i)\}_{i=1}^{n_{\text{train}}}$ are defined on the reference mesh $\mathcal{T}_h^0$, whereas  the  solution snapshots $\{\bm u_h(\bm x; \bm \mu_i)\}_{i=1}^{n_{\text{train}}}$ are defined on different  $\bm \mu$-dependent adaptive meshes $\mathcal{T}_h^{\bm \mu_i}$.  If $\bm x$ is mesh point in $\mathcal{T}_h^0$ then $\bm \phi_h(\bm x; \bm \mu_i)$ is a mesh point in $\mathcal{T}_h^{\bm \mu_i}$. We then introduce 
\begin{equation}
    \tilde{\bm{u}}_h(\bm x; \bm \mu_i) = \bm {u}_h(\bm \phi (\bm x; \bm \mu_i); \bm \mu_i), \quad \forall \bm x \in \mathcal{T}_h^0. 
\end{equation}
The resulting snapshots $\{\tilde{\bm{u}}_h(\bm x; \bm \mu_i)\}_{i=1}^{n_{\text{train}}}$ are defined on the reference mesh $\mathcal{T}_h^0$. On the reference mesh, we anticipate that sharp features will be smoothed out and brought into closer alignment, making model reduction more effective. Model reduction can then be performed on the snapshots of $\bm \phi_h$, which are anticipated to be smooth due to the elliptic nature of the Monge-Ampère equation, and the snapshots of $\tilde{\bm u}_h$, which are anticipated to be more regular than $\bm u_h$ and with less parametric variation. 

\begin{figure}[htb!]
\begin{center}
\subfigure[Solutions $\bm u_h(\bm x; \bm \mu_i)$ (left) and adapted meshes $\mathcal{T}_h^{\bm \mu_i}$  (right) for $\bm \mu_i \in \{ 2, 3, 4 \}$. Plotted field is physical density.]
{
\includegraphics[width=0.12\textwidth]{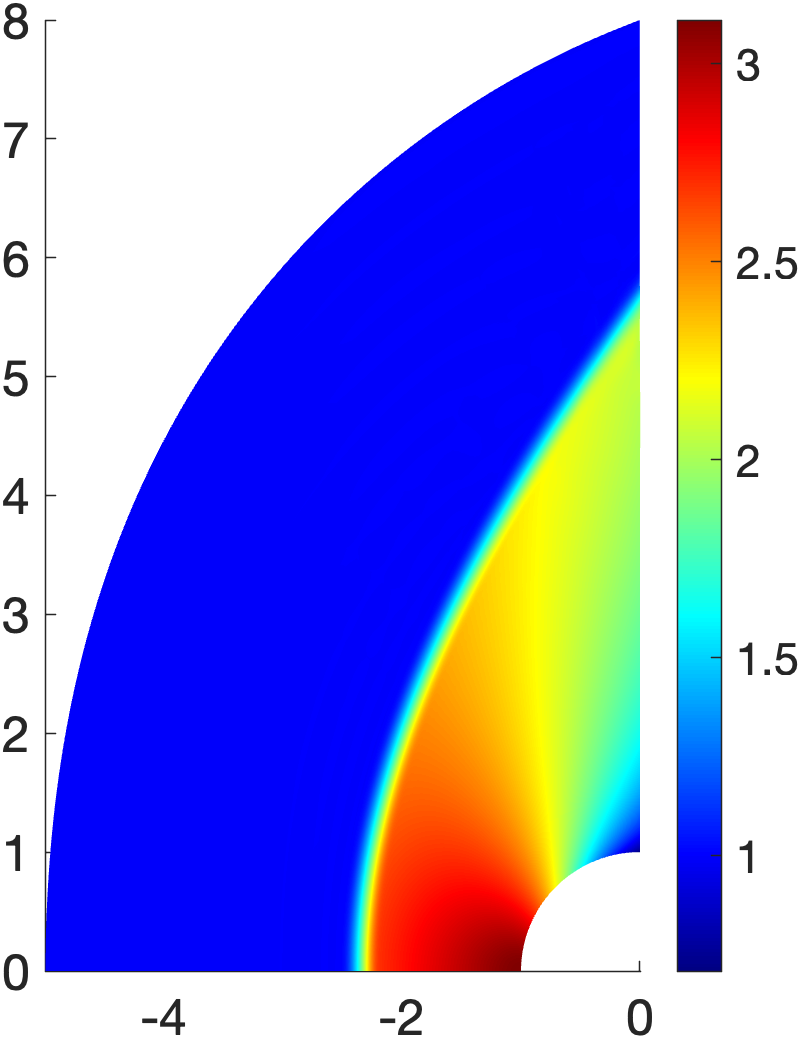}\label{fig:2sol}
\includegraphics[height=2.7cm,width=0.1\textwidth]{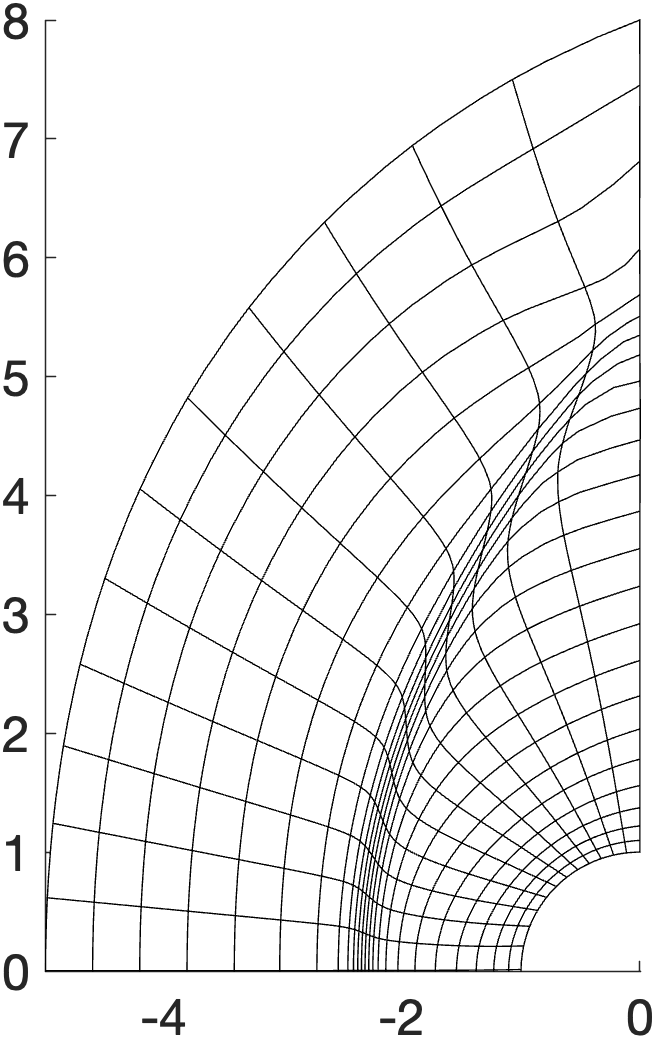}\label{fig:mach2mesh}
\qquad
\includegraphics[width=0.12\textwidth]{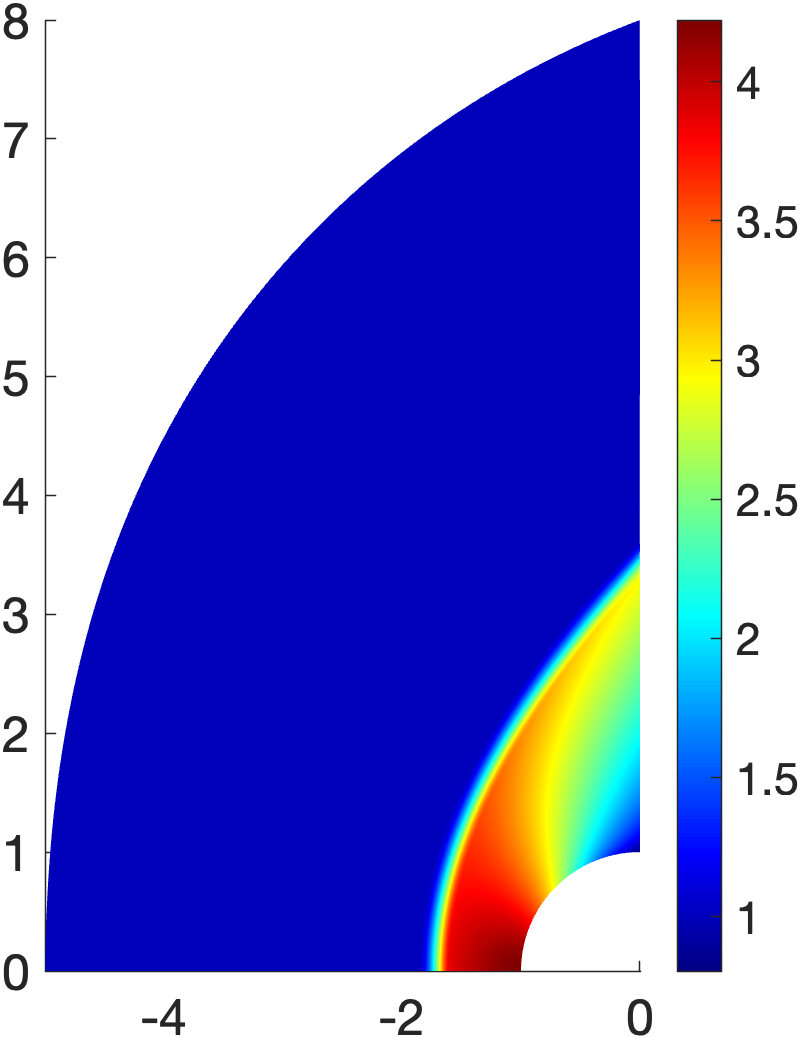}\label{fig:3sol}
\includegraphics[height=2.7cm,width=0.1\textwidth]{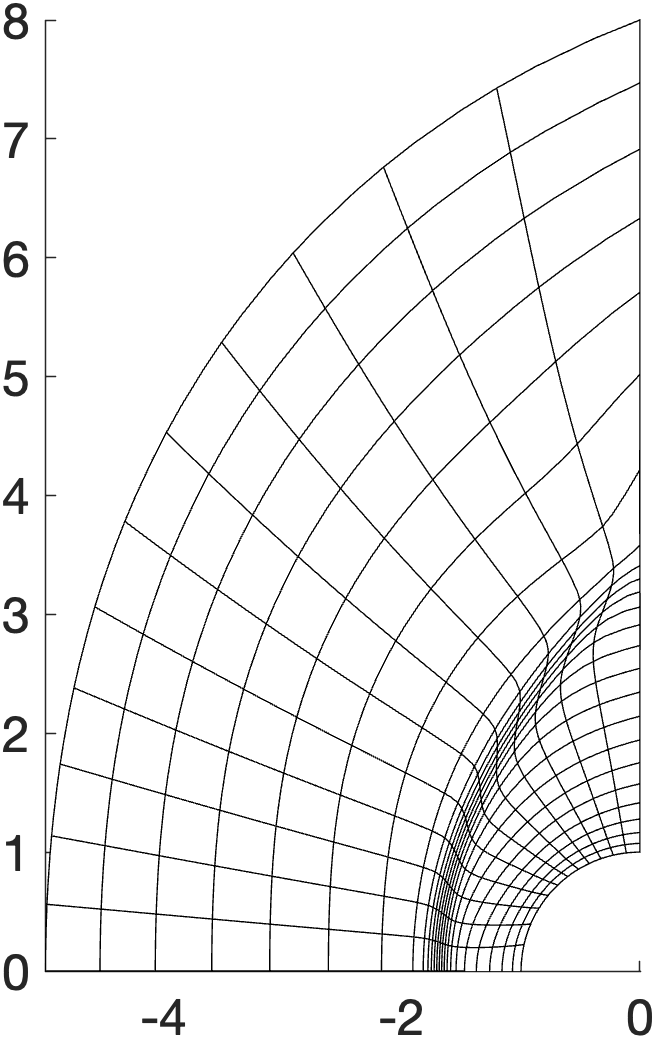}\label{fig:mach3mesh}
\qquad
\includegraphics[width=0.12\textwidth]{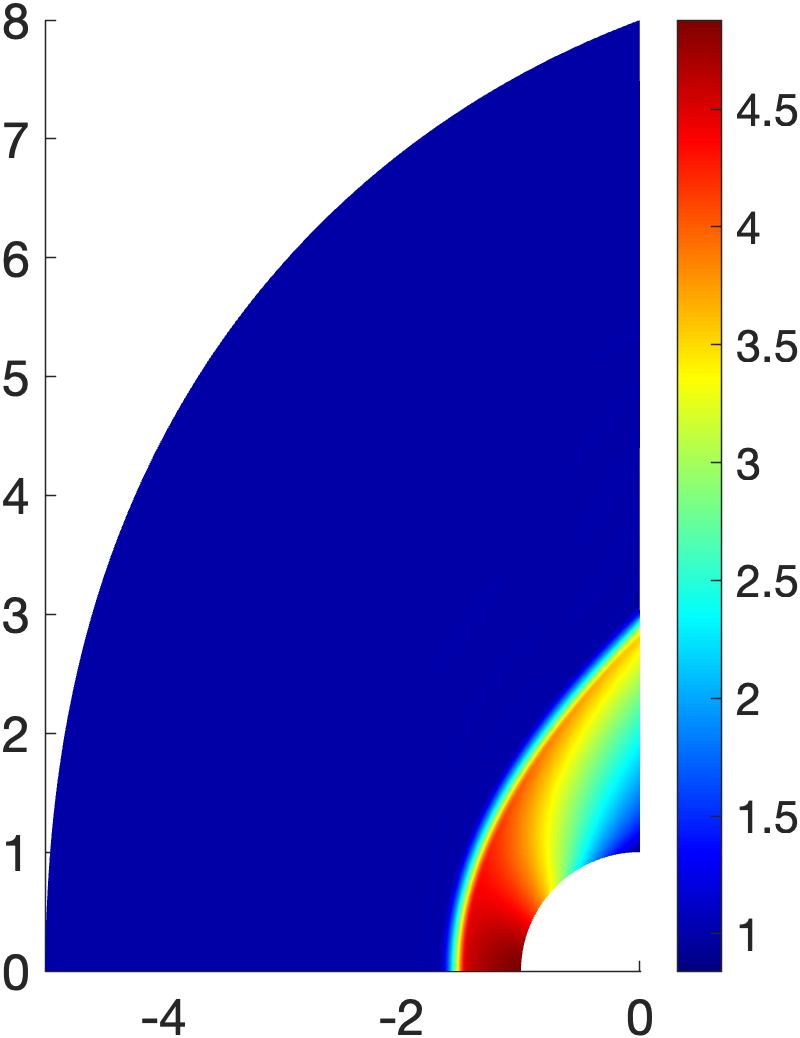}\label{fig:4sol}
\includegraphics[height=2.7cm,width=0.1\textwidth]{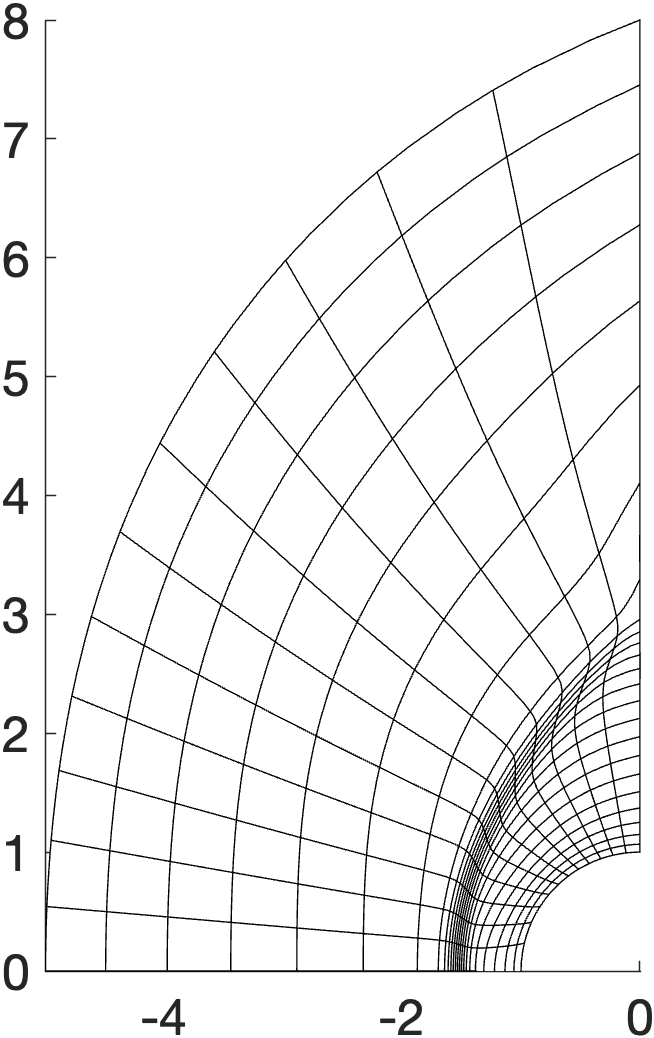}\label{fig:mach4mesh}
}\\

\subfigure[$x$-component (left) and $y$-component (right) of $\phi_h(\bm x; \bm \mu_i)$ for $\bm \mu_i \in \{2, 3, 4\}$]
{\includegraphics[width=0.12\textwidth]{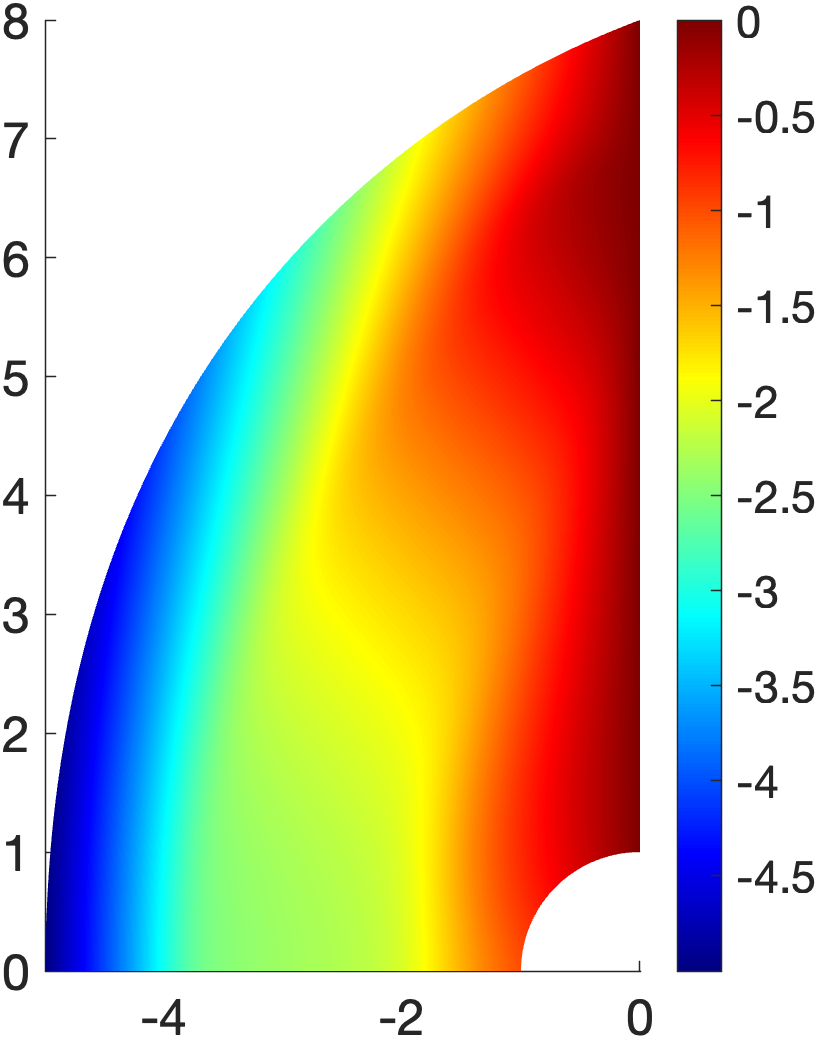}\label{fig:2qx} 
\includegraphics[width=0.12\textwidth]{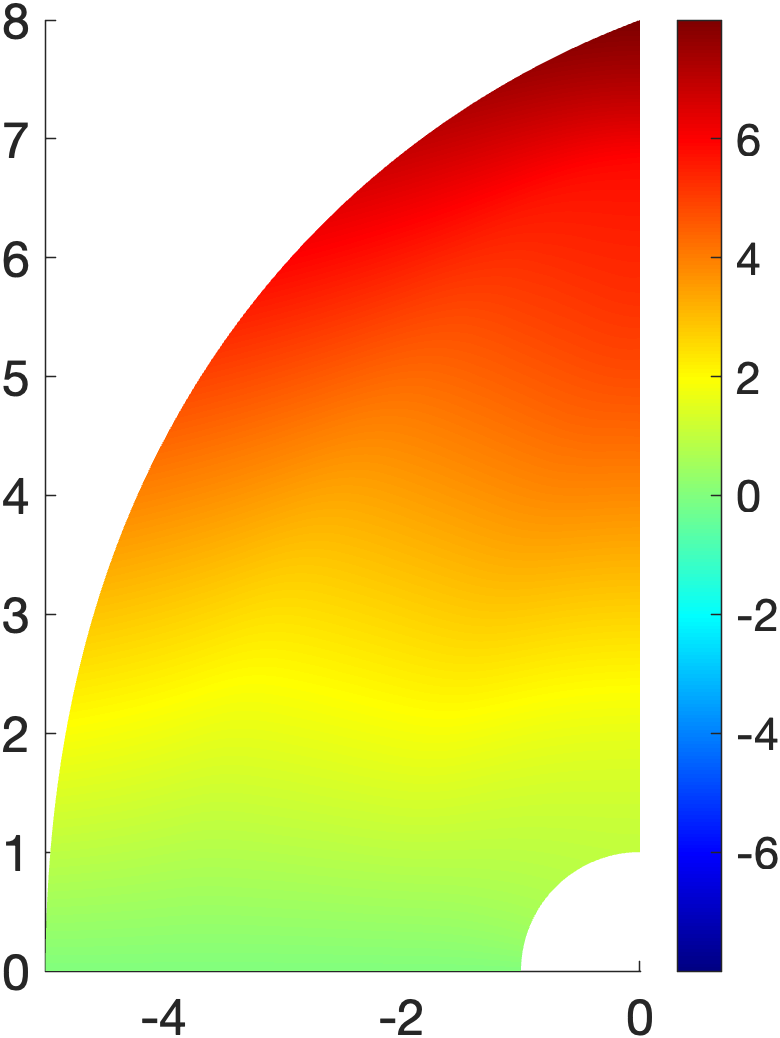}\label{fig:2qy} 
\qquad
\includegraphics[width=0.12\textwidth]{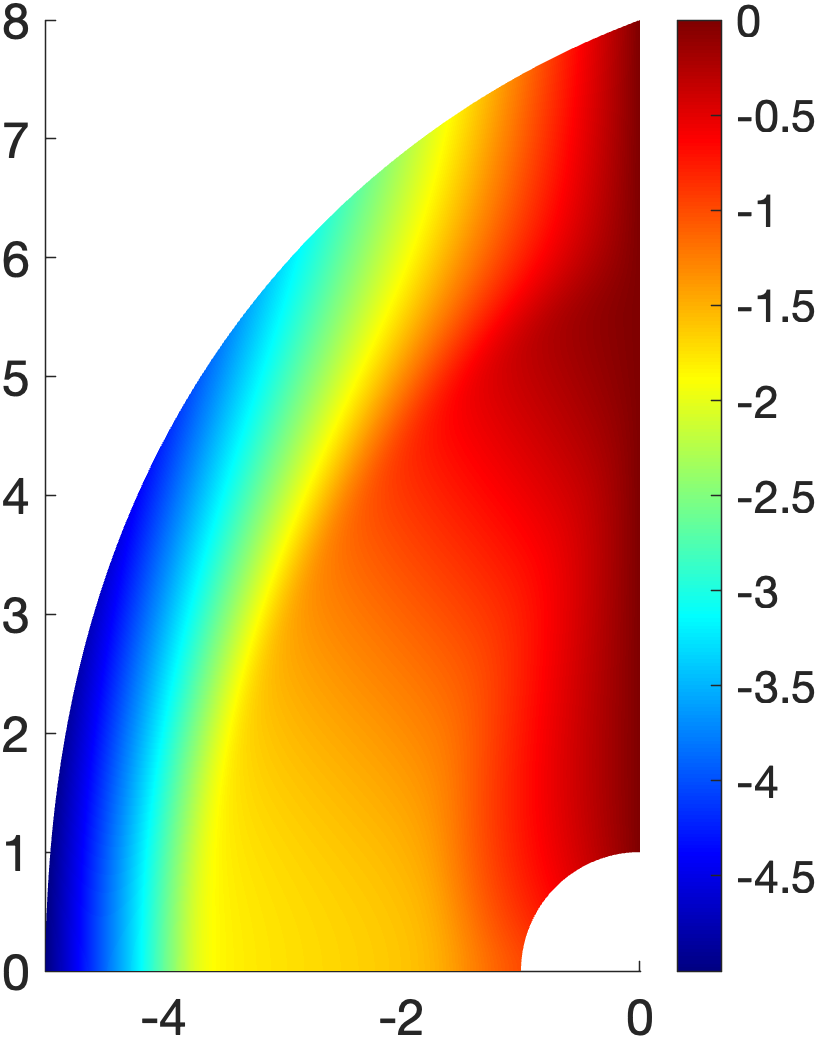}\label{fig:3qx} 
\includegraphics[width=0.12\textwidth]{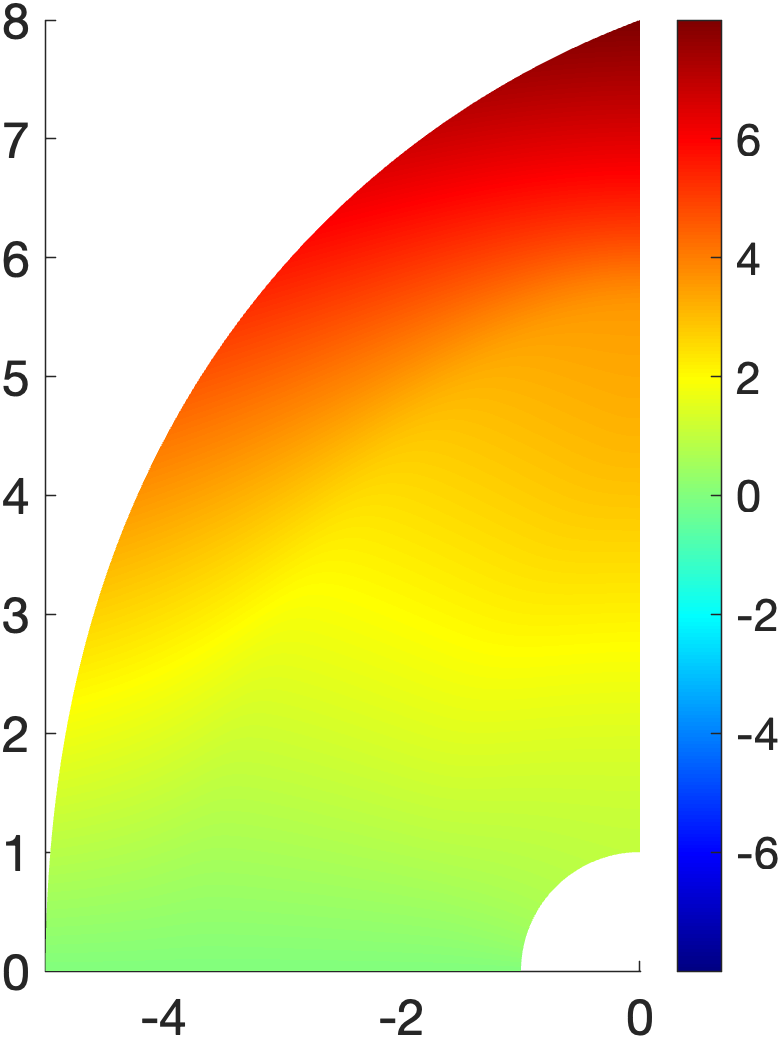}\label{fig:3qy}
\qquad
\includegraphics[width=0.12\textwidth]{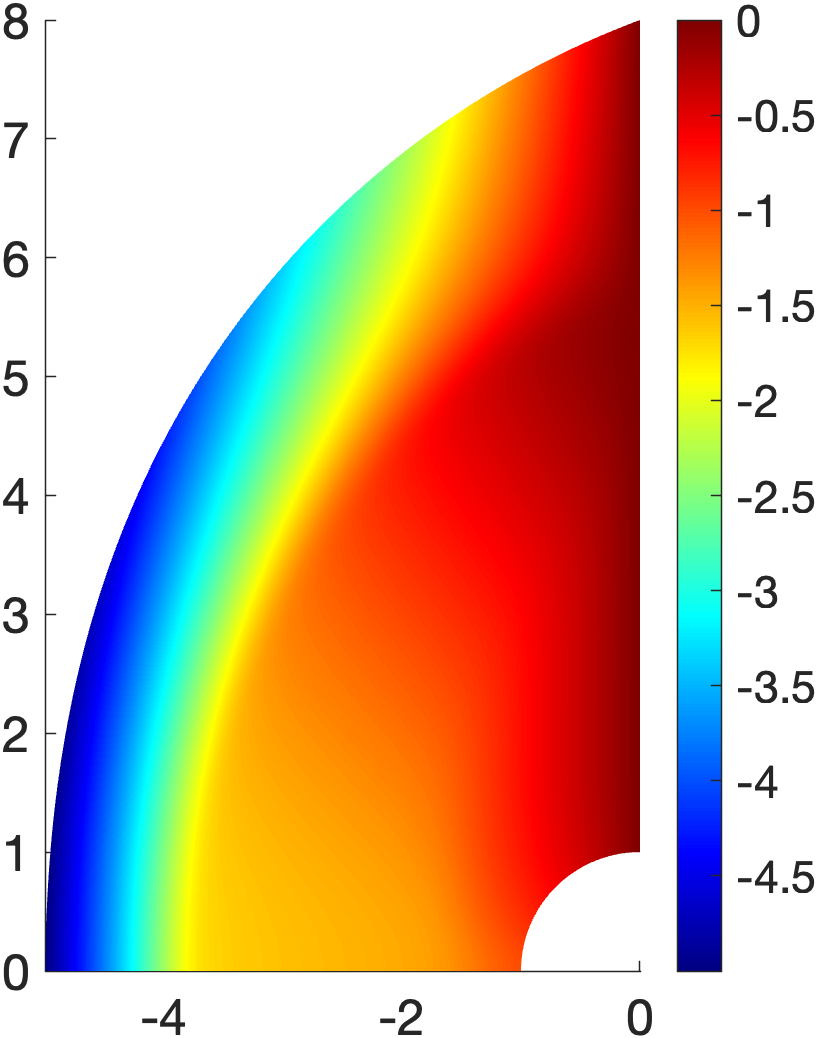}\label{fig:4qx} 
\includegraphics[width=0.12\textwidth]{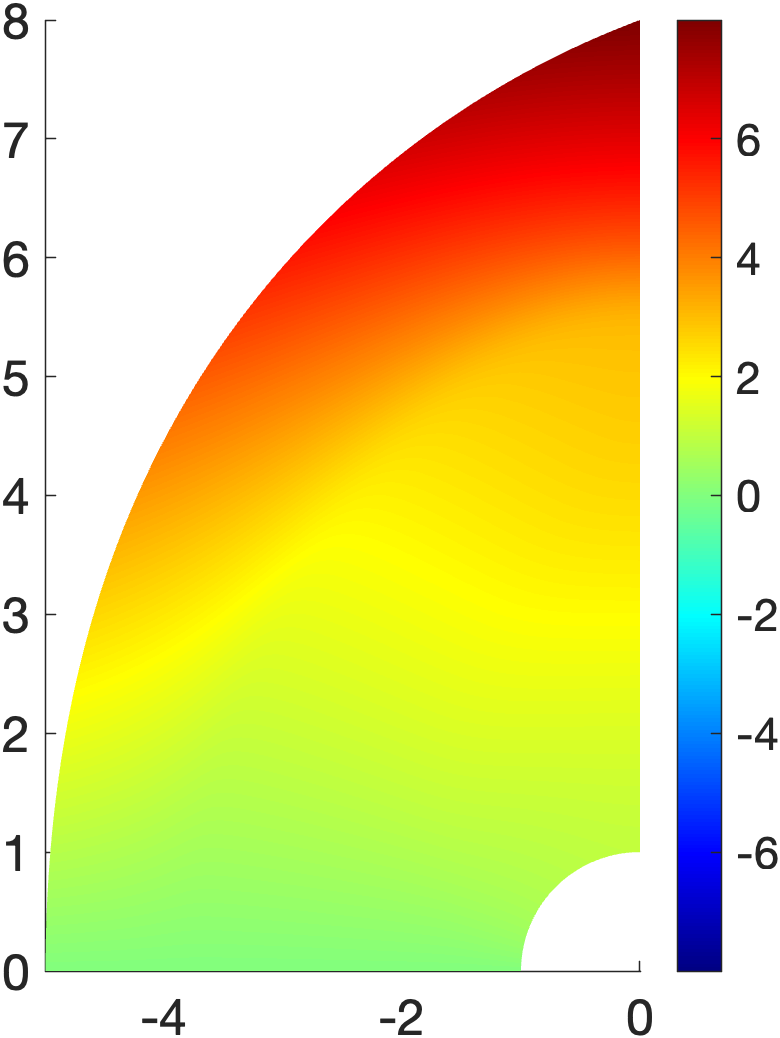}\label{fig:4qy}
}
\\
\subfigure[Mapped solutions $\tilde{\bm{u}}_h(\bm x; \bm \mu_i)$ for $\bm \mu_i \in \{2, 3, 4\}$. Plotted quantities are density (left) and pressure (right). ]
{
\includegraphics[width=0.12\textwidth]{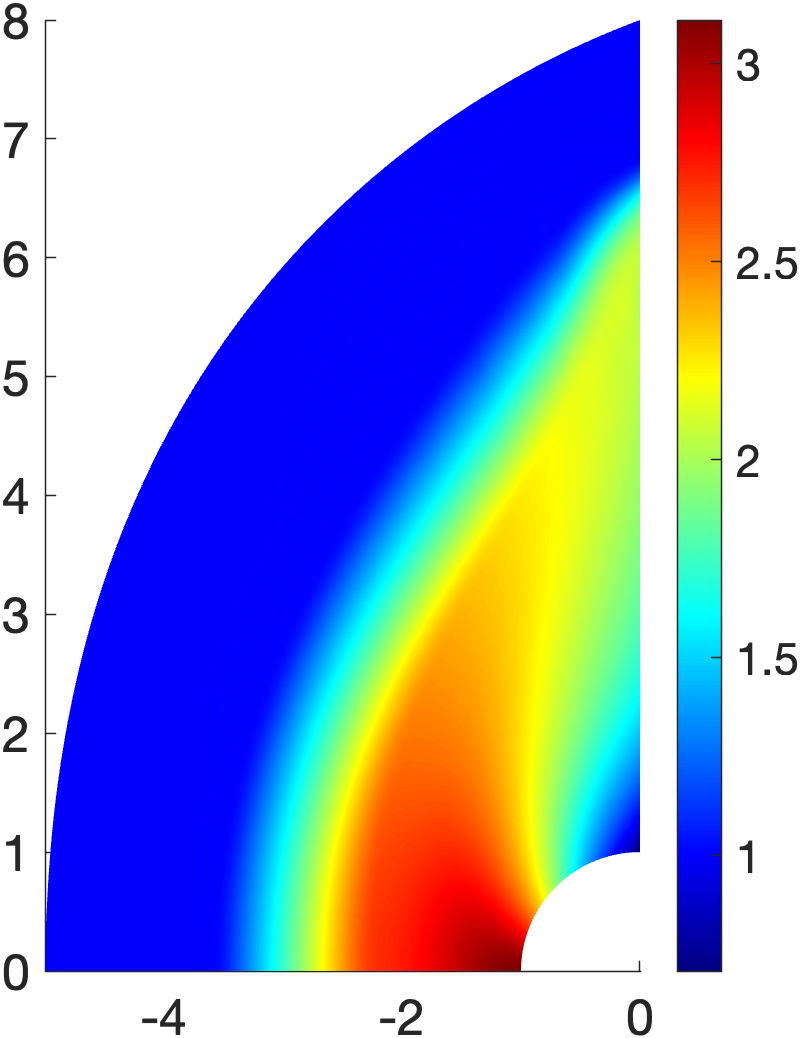}\label{fig:2solmap}
\includegraphics[width=0.12\textwidth]{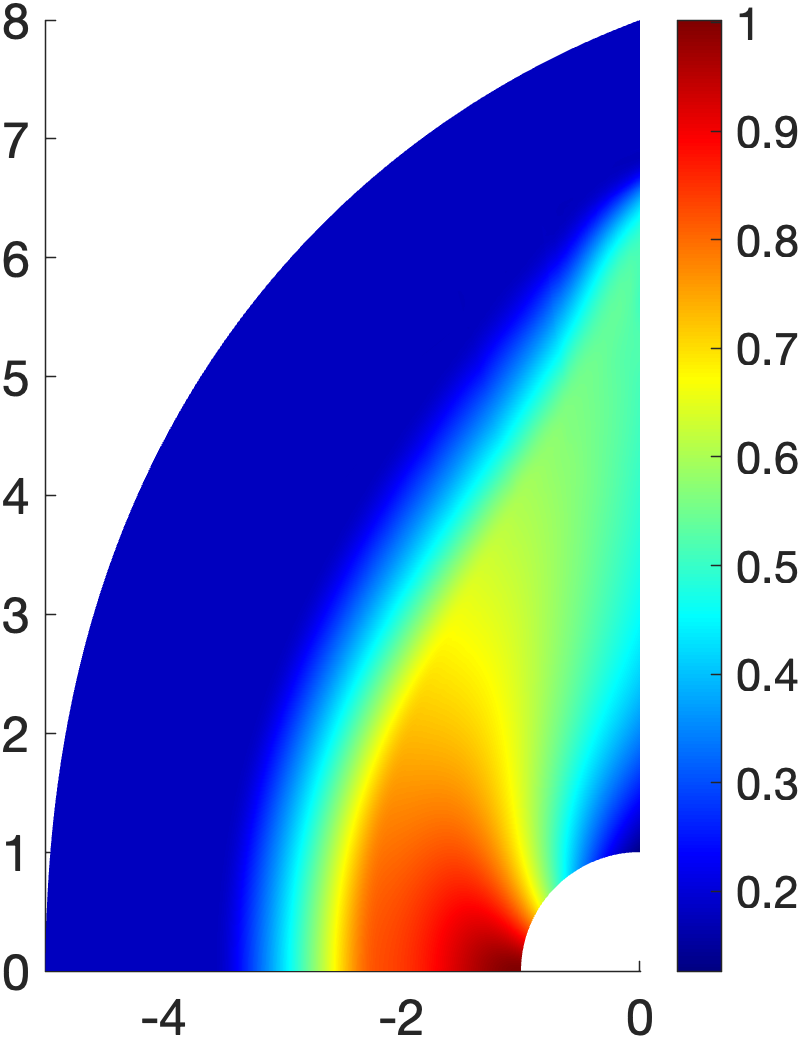}\label{fig:2pressmap}
\qquad
\includegraphics[width=0.12\textwidth]{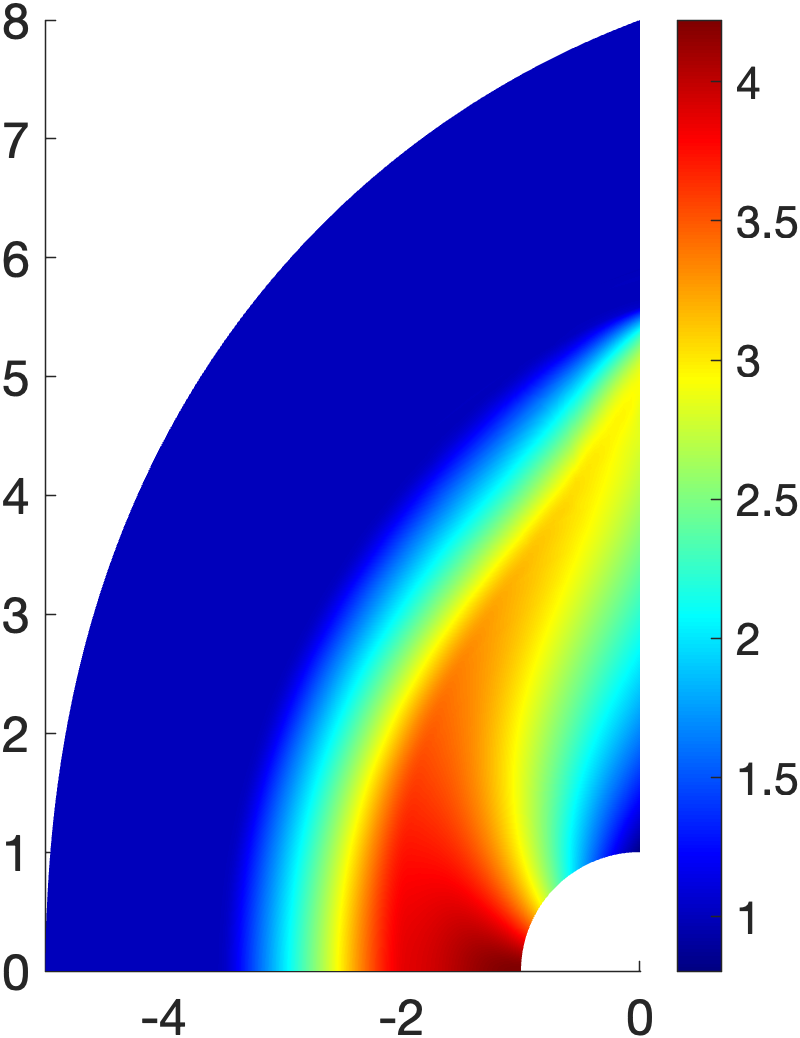}\label{fig:3solmap}
\includegraphics[width=0.12\textwidth]{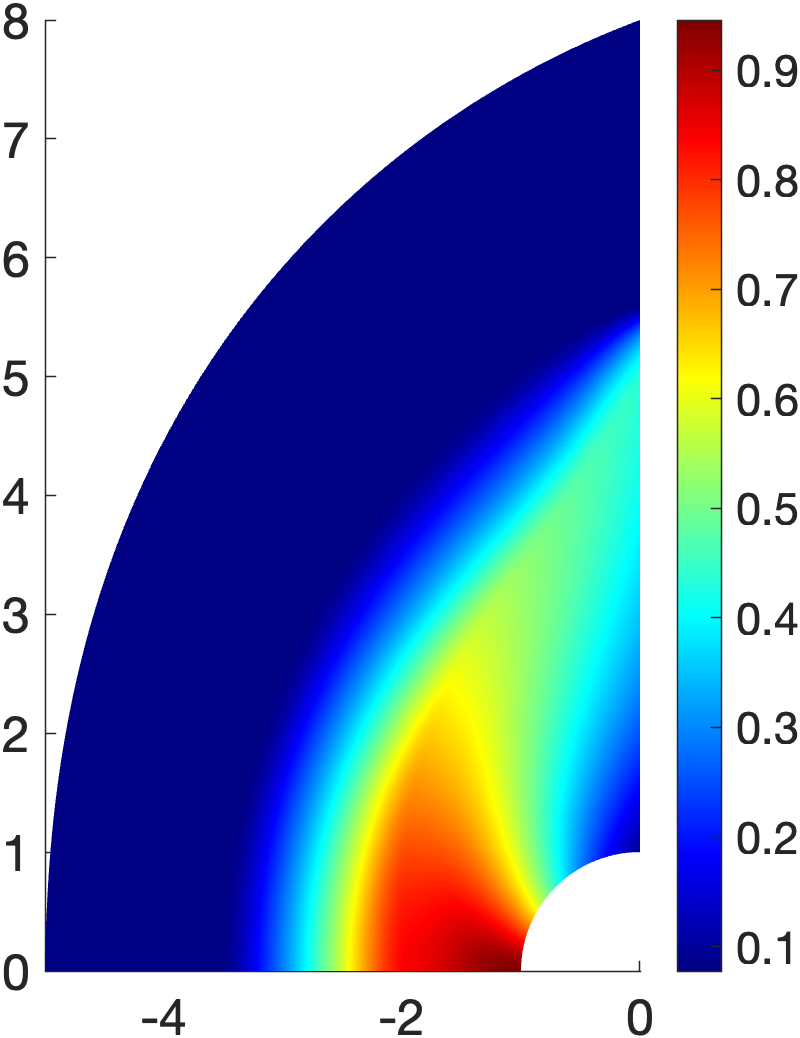}\label{fig:3pressmap}
\qquad
\includegraphics[width=0.12\textwidth]{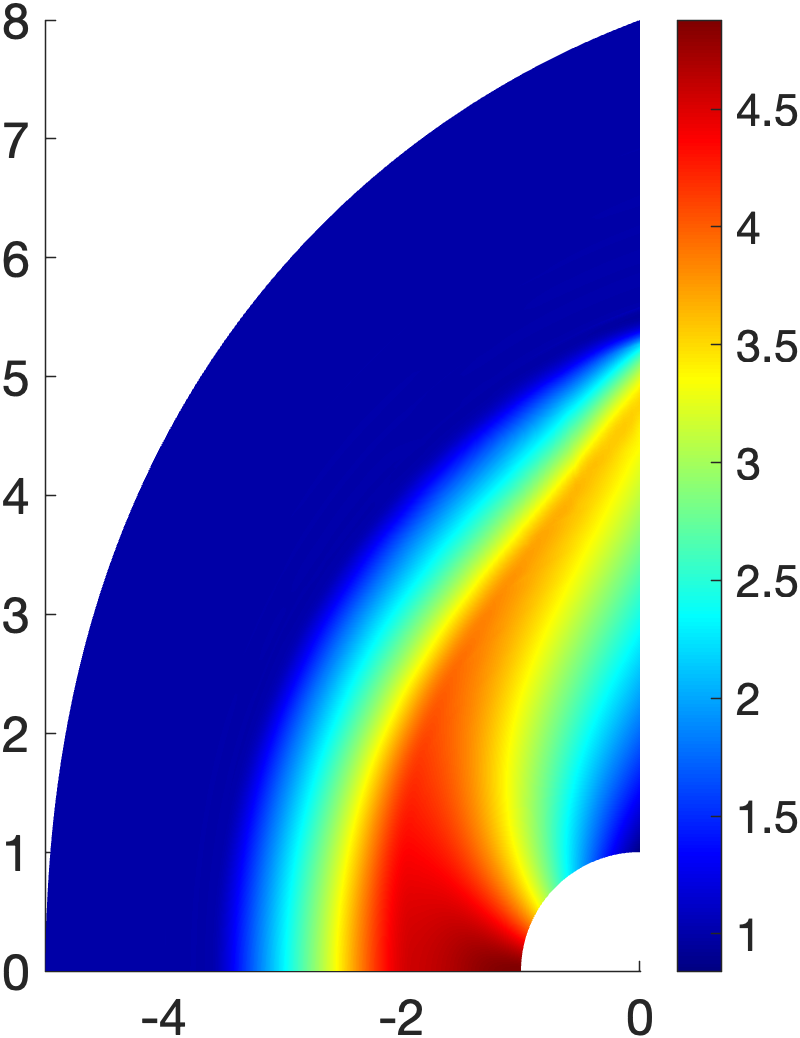}\label{fig:4solmap}
\includegraphics[width=0.12\textwidth]{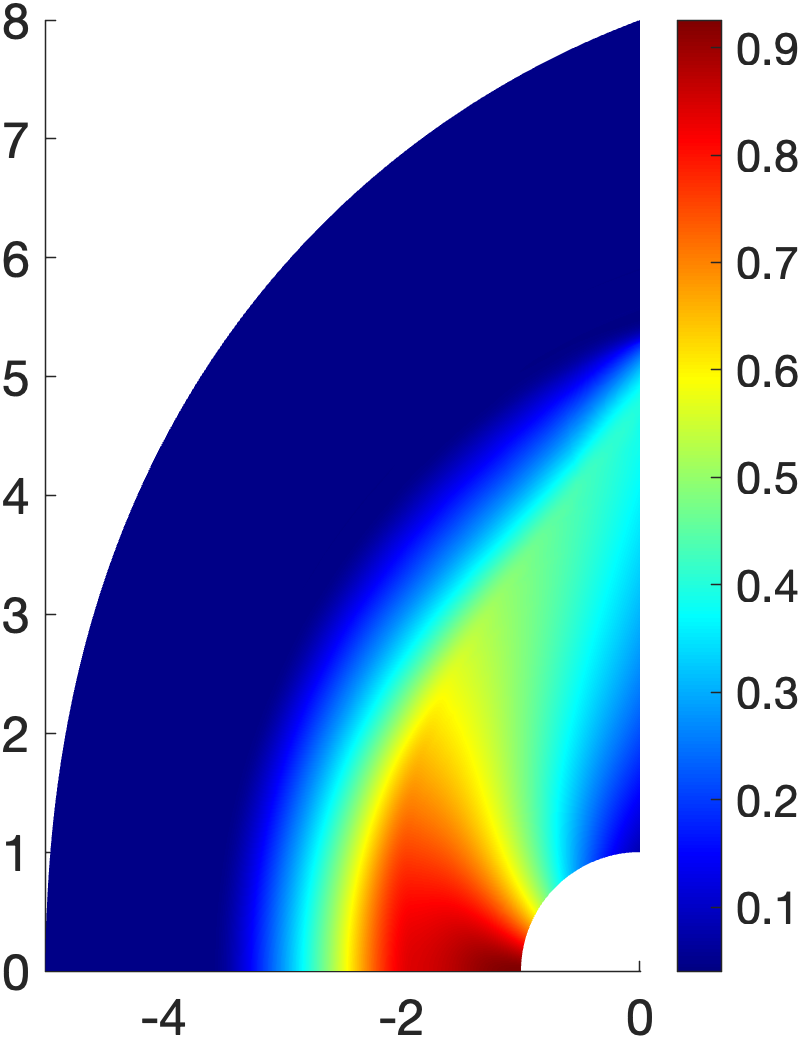}\label{fig:4pressmap}

}
\caption{Sample quantities used to build ROM for different parameter samples. Top row contains the outputs of Algorithm 1, solutions $\{\bm u_h(\bm x; \bm \mu_i)\}_{i=1}^{3}$ and adapted meshes $\{ \mathcal{T}_h^{\bm \mu_i} \}_{i=1}^{3}$. 
Middle row shows the corresponding mappings $\{\bm \phi_h(\bm x; \bm \mu_i)\}_{i=1}^3$ for each adapted mesh. Bottom row displays solutions mapped back to the reference domain, $\{\tilde{\bm{u}}_h(\bm x; \bm \mu_i)\}_{i=1}^{3}$. }\label{fig:trainingset}
\end{center}
\end{figure}

The quantities used to construct these ROMs are shown in Figure \ref{fig:trainingset} using the same high-speed cylinder flow introduced in Subsection \ref{sec:adapt} with a training set of $\mathcal{D} = \{2, 3, 4\}$.  

\subsection{Interpolation ROMs}
\label{sec:interp}
We construct nonintrusive interpolation-based ROMs for the mapping and the solution both defined on the reference mesh. 

For simplicity, we consider external supersonic and hypersonic configurations with one parameter allowed to vary, the free stream Mach number $\bm \mu = \mathrm{Ma}_{\infty}$, and with the training set fixed for each problem.

We compress the snapshots of the solution and mapping with the proper orthogonal decomposition approach. 
Snapshots of $\tilde{\bm u}_h$ and $\bm \phi_h$ are arranged into matrices $\bm X_{\bm \tilde{u}}$ and $\bm X_{\bm \phi}$ respectively. 
The SVD is taken of each matrix to extract the leading modes
\begin{equation}
    \bm U_{\tilde{\bm u}} \bm \Sigma_{ \tilde{\bm u}} \bm V^T_{\tilde{\bm u}} = \bm X_{\bm \tilde{u}}, \; \; \bm U_{\bm \phi} \bm \Sigma_{\bm \phi} \bm V^T_{\bm \phi} = \bm X_{\bm \phi}, 
\end{equation}
and the reduced bases for the solution and maps $\bm W_{\tilde{\bm u}}$ and $\bm W_{\bm \phi}$ are defined as the $N$ leading columns of $\bm U$, where $N \le n_{\rm train}$ is the reduced basis dimension. 

For any given $\bm \mu \in \mathcal{D}$, the best approximations of the solution and the mapping are given by
\begin{equation}
\tilde{\bm u}^{\ast}_N(\bm \mu) = \bm W_{ \tilde{\bm u}} \bm \alpha^{\ast}_{\tilde{\bm u}}(\bm \mu), \qquad 
    \bm \phi^{\ast}_N(\bm \mu) = \bm W_{\bm \phi} \bm \alpha^{\ast}_{\bm \phi}(\bm \mu),
\end{equation}
where the $\bm \alpha^*$ vectors are the coefficients of the best approximations  and are given by
\begin{equation}
\bm \alpha^{\ast}_{\tilde{\bm u}}(\bm \mu) = \bm W_{\tilde{\bm u}}^T \tilde{\bm  u}_h(\bm \mu), \qquad 
\bm \alpha^{\ast}_{\bm \phi}(\bm \mu) = \bm W_{\bm \phi}^T \bm \phi_h(\bm \mu) .
\end{equation}
Because the coefficients of the best approximations require the solution and the mapping, they can be only computed for $\bm \mu = \bm \mu_i, 1 \le i \le n_{\rm train},$ in the training set. As a result, the best approximations are for those parameter vectors in the training set.

For any $\bm \mu \in \mathcal{D}$ which does not reside in the training set, we compute the ROM approximations of the solution and mapping as follows
\begin{equation}
    \bm \tilde{u}_N(\bm \mu) = \bm W_{ \tilde{\bm 
 u}} \bm \alpha_{\tilde{\bm u}}(\bm \mu), \qquad 
    \bm \phi_N(\bm \mu) = \bm W_{\bm \phi} \bm \alpha_{\bm \phi}(\bm \mu). 
\end{equation}
The $\bm \alpha$ vectors are the coefficients of the solution expanded in the reduced basis $\bm W$ and are referred to as the generalized coordinates. 

We use the simple radial basis function (RBF) surrogate used as an initial guess for an intrusive ROM in \cite{ching2022gridtailor}. 
More sophisticated RBF ROMs are possible, usually involving some calibration with a mix of training and testing sets \cite{ dolci2016proper, xiao2015non, georgaka2020hybrid}.

The RBF interpolant is defined as 
\begin{equation}
    \alpha_n(\bm \mu) = \sum_{i=1}^{n_{\text{train}}} \beta_{ni} \Psi \left( \|\bm \mu - \bm \mu_i\|_2 \right), \quad 1 \le n \le N,
\end{equation}
where $\Psi$ are radial basis functions. For this work, we use a multiquadratic radial basis function with shape parameter set to 20. 
The coefficients $\beta_{ni}$ are determined by enforcing that $\bm \alpha(\bm \mu_i)$ are equal to $\bm \alpha^*(\bm \mu_i)$ for all $\bm \mu_i$ in the training set and solving a $n_{\text{train}} \times n_{\text{train}}$ linear system.
In \cite{ching2022gridtailor}, the resulting interpolant was found to be reasonably competitive with an intrusive ROM when interpolating in parameter space, considering the relative ease of training and online evaluations.

\section{Numerical examples}
\label{sec:results}

For each case we compare a \textit{mapped ROM} to a \textit{fixed mesh ROM}. The mapped ROM uses the quantities described in Subsection \ref{sec:ROMref}, the mappings $\{\bm \phi_h(\bm x; \bm \mu_i)\}_{i=1}^{n_{\text{train}}}$ and mapped snapshots $\{ \tilde{\bm u}_h(\bm x; \bm \mu_i)\}_{i=1}^{n_{\text{train}}}$. The fixed mesh ROM builds a ROM using snapshots of $\{ \bm u^0_h(\bm x; \bm \mu_i)\}_{i=1}^{n_{\text{train}}}$, from step 1 of Algorithm \ref{alg:meshadapt}. 
We use a reasonable starting mesh for each case so that the snapshots of $\bm u_h^0$ are of acceptable quality.

We report the relative $\ell^2$ errors of a vector field $\bm v_h(\bm x; \bm \mu)$. For the mixed mesh ROM, $\bm v_h$ is $\bm u^0(\bm x; \bm \mu)$. For the mapped mesh ROMs there are two fields to consider, $\bm \phi_h(\bm x; \bm \mu)$ and $\tilde{\bm u}_h(\bm x; \bm \mu)$. The $\ell^2$ error over a test set $\mathcal{S}_{\text{test}} \in \mathcal{D}$ on the reference domain $\Omega$ is given by 
\begin{equation}
    E_v(\bm \mu) = \sqrt{\frac{\int_{\Omega} \| \bm v_h(\bm x; \bm \mu) - \bm W_{\bm v} \bm \alpha_{\bm v} (\bm \mu)\|_2^2}{\int_{\Omega} \| \bm v(\bm x; \bm \mu) \|_2^2 }}, \; \; \forall \bm \mu \in \mathcal{S}_{\text{test}}.
    \label{eq:error}
\end{equation}

Errors are computed on the reference mesh due to the fact that there may be non-negligible error in the mappings, so a reference solution at any test parameter will have a slightly different grid than the ROM. 
Furthermore, the errors on the reference mesh will more clearly show whether the solution $\tilde{\bm u}_h$ is indeed more amenable to model reduction on the reference domain for practical mesh mappings. 
For examples with a wide range of Mach numbers, we evaluate at test sets that are equally spaced in between the training sets in order to account for the full variation of flow conditions. 
In all cases, the mean and max errors over the test set are reported. 

\subsection{Flow around a cylinder}
We consider inviscid flow around a cylinder geometry with parametrically varying Mach number. 
These configurations are a common benchmark for shock capturing methods for high-speed flow \cite{fernandez2018physics, ching2019shock, lal2023aeroheating, bai2022continuous}.
The flow exhibits a strong curved bow shock that varies parametrically with the Mach number. The free-stream Mach number is varied between supersonic and hypersonic ranges, with $\mathcal{D} = [2, 10]$. 

\begin{figure}[htb!]
    \centering
    \includegraphics[width=0.7\textwidth]{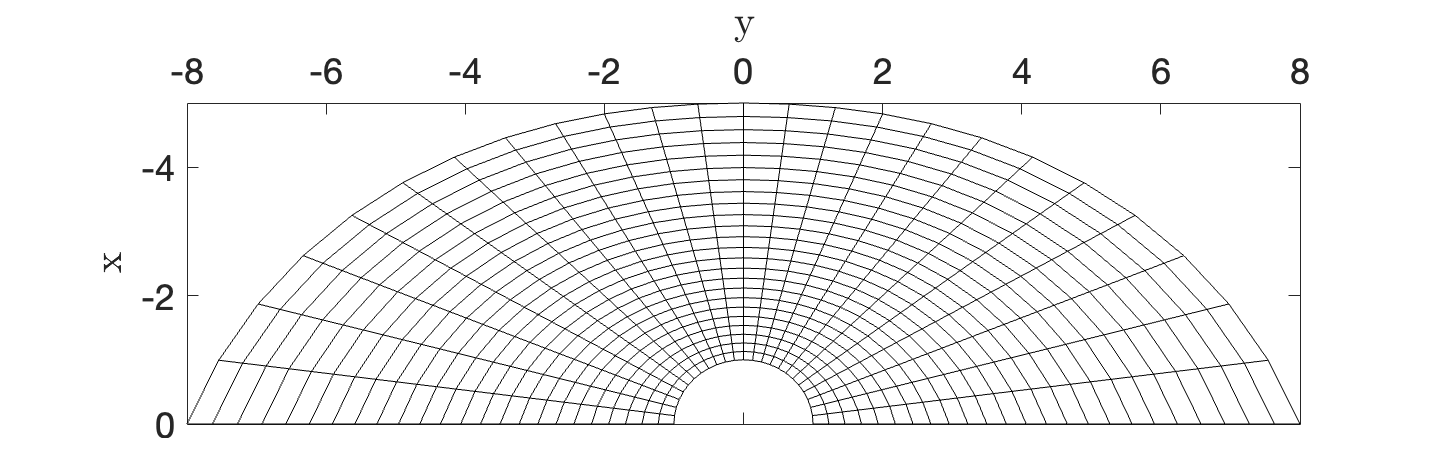}
    \caption{Starting structured grid of $\Omega$ for the high-speed cylinder.}
    \label{fig:cylmesh}
\end{figure}

For each FOM solution, we begin with a reference mesh with 25 elements in the radial and axial directions, with polynomial order 3. The reference mesh is shown in Figure \ref{fig:cylmesh}.  Each FOM is solved from a uniform solution initialized to the free-stream values with no polynomial order continuation.  

\begin{figure}[htb!]
\begin{center}
\subfigure[$\mathrm{Ma}_{\infty} = 2.0$]
{
\includegraphics[width=0.4\textwidth]{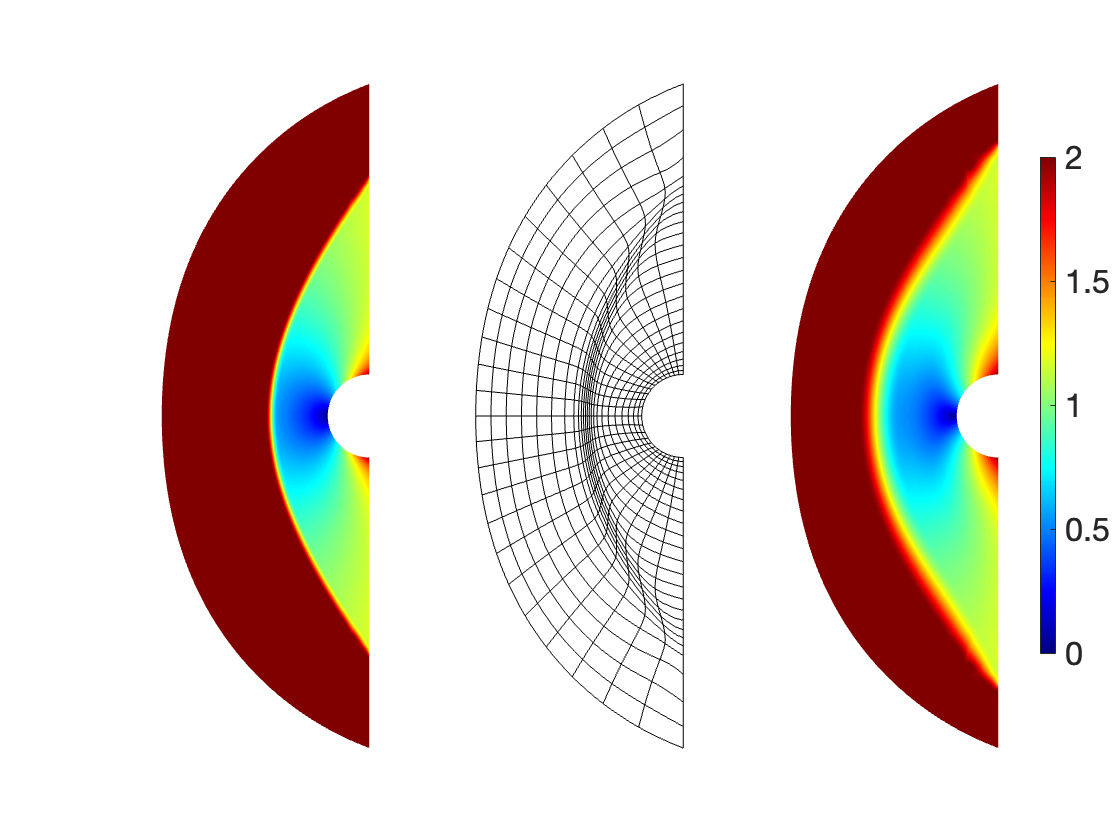}\label{fig:cyl2}}
\subfigure[$\mathrm{Ma}_{\infty} = 4.0$]
{
\includegraphics[width=0.4\textwidth]{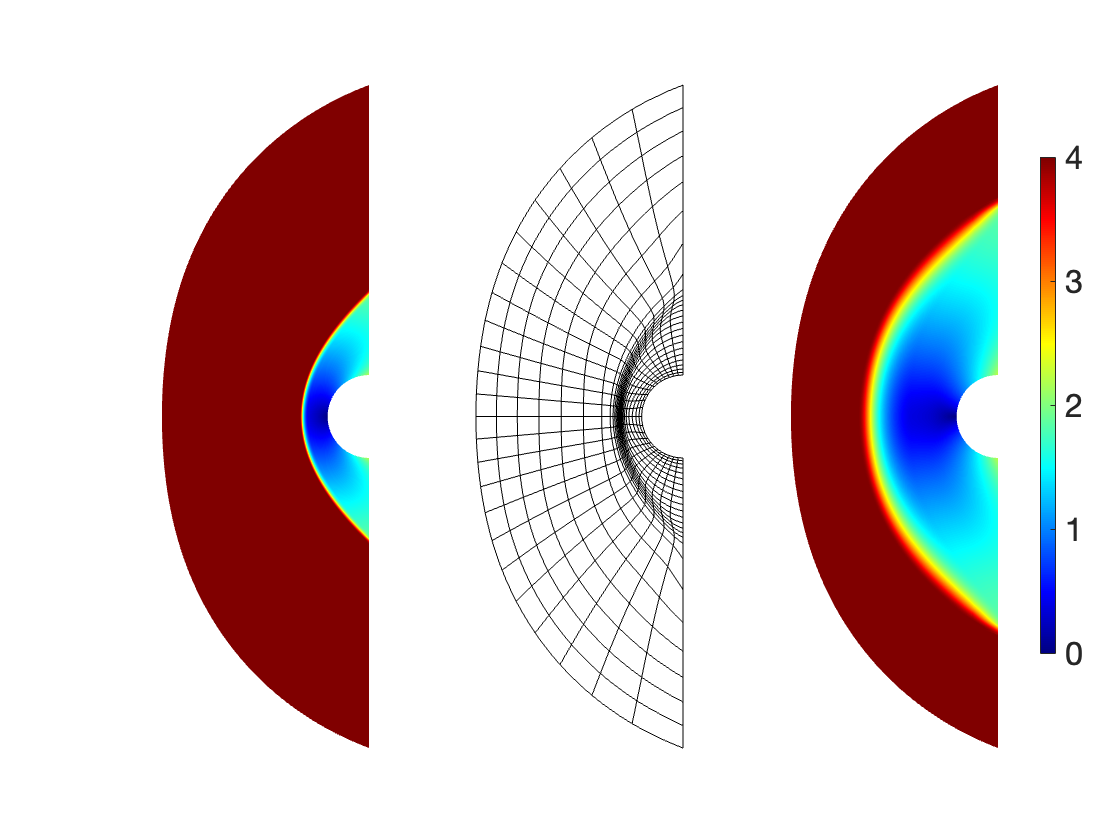}\label{fig:cyl4}}
\subfigure[$\mathrm{Ma}_{\infty} = 10.0$]
{\includegraphics[width=0.4\textwidth]{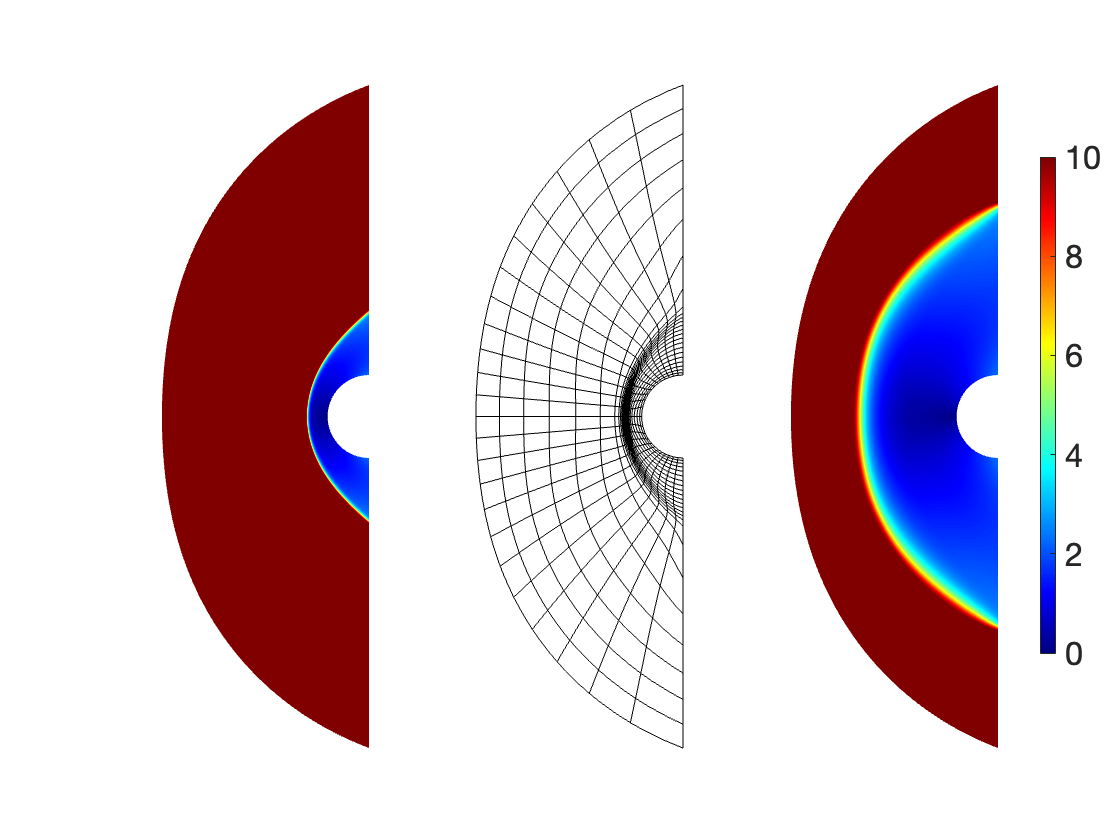}\label{fig:cyl10}} 
\caption{For different Mach numbers, a visualization of the solution $\bm u_h(\bm x)$, the corresponding adapted mesh $\mathcal{T}_h^{\bm \mu}$ and the solution mapped back to the reference domain $\tilde{\bm u}_h$. Plotted quantity is the Mach number.}\label{fig:cyl_sweep}
\end{center}
\end{figure}

The solutions $\bm u_h(\bm x; \bm \mu)$,  adapted meshes $\mathcal{T}_h^{\bm \mu}$, and mapped solutions $\tilde{\bm u}_h(\bm x; \bm \mu)$ for representative values choices of $\bm \mu$ are shown in Figure \ref{fig:cyl_sweep}. 

Although the snaphots of $\tilde{\bm u_h}$ are not perfectly aligned, the parametric shock motion is substantially reduced.

\begin{figure}[!htb]
\begin{center}
\subfigure
{
\includegraphics[width=0.45\textwidth]{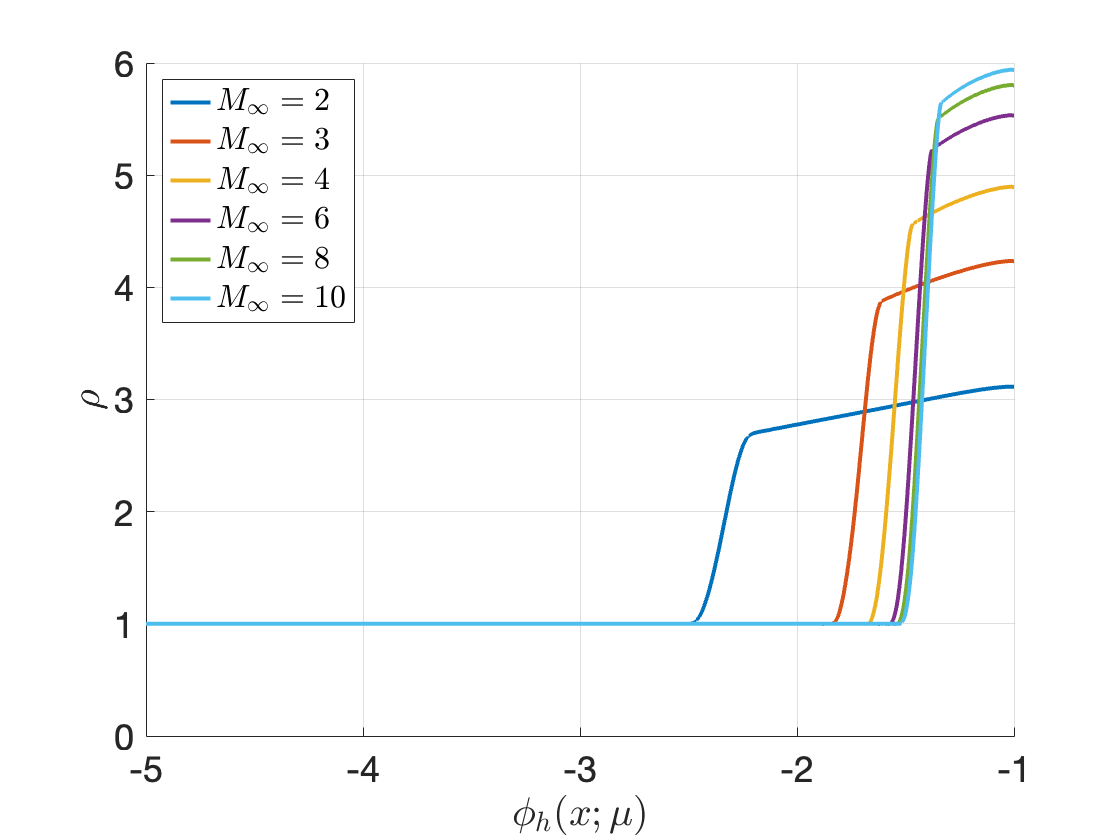}\label{fig:stagfixed}}
\subfigure
{
\includegraphics[width=0.45\textwidth]{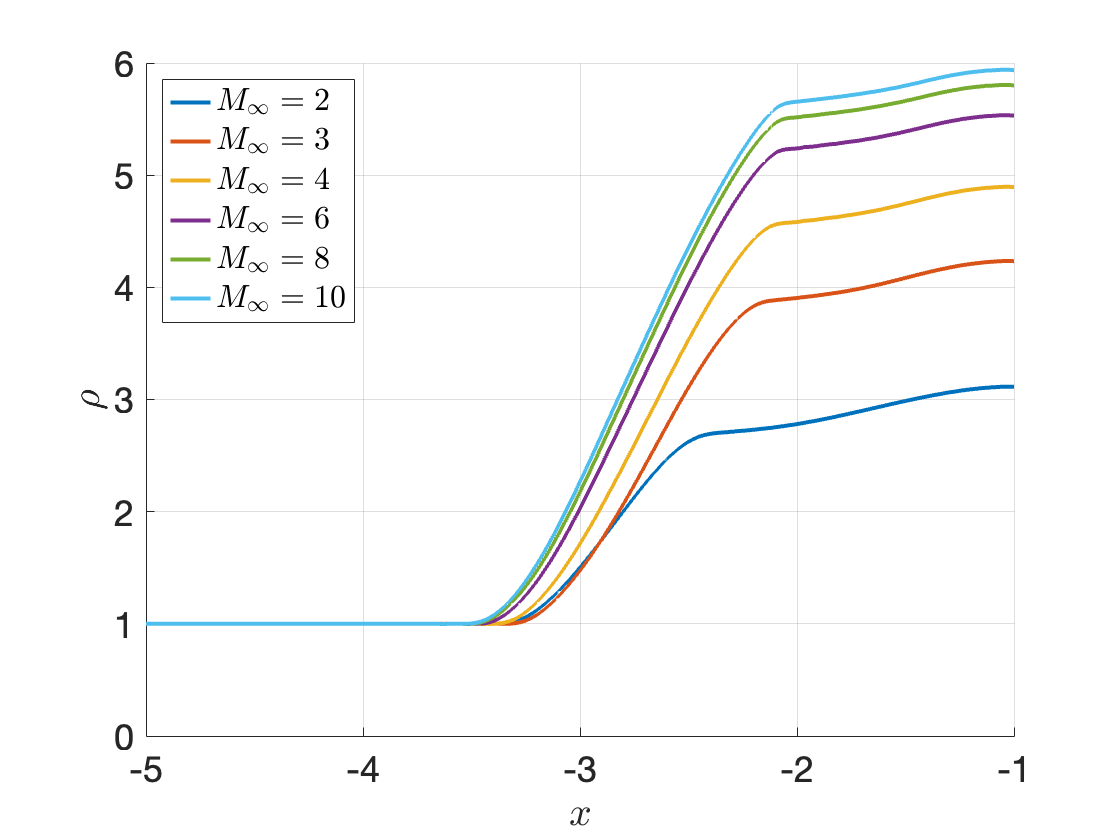}\label{fig:stagmapped}}
\caption{Density along the stagnation line on the adapted meshes $\mathcal{T}_h^{\bm \mu}$ and pulled back to the reference mesh $\mathcal{T}_h^0$ . The jump in the solution is well-aligned and is spread out over the domain.}\label{fig:staglines}
\end{center}
\end{figure}

Most of the shock movement occurs towards the outflow boundaries; along the stagnation line, the shock is brought into near alignment. See Figure \ref{fig:staglines}. 
For shocks that are stabilized with some artificial viscosity, reversing the mapping also has the effect of smearing out sharp features.
The combination of alignment and smearing should abet classical model reduction procedures. 

\begin{figure}
    \centering
    \includegraphics[width=0.5\textwidth]{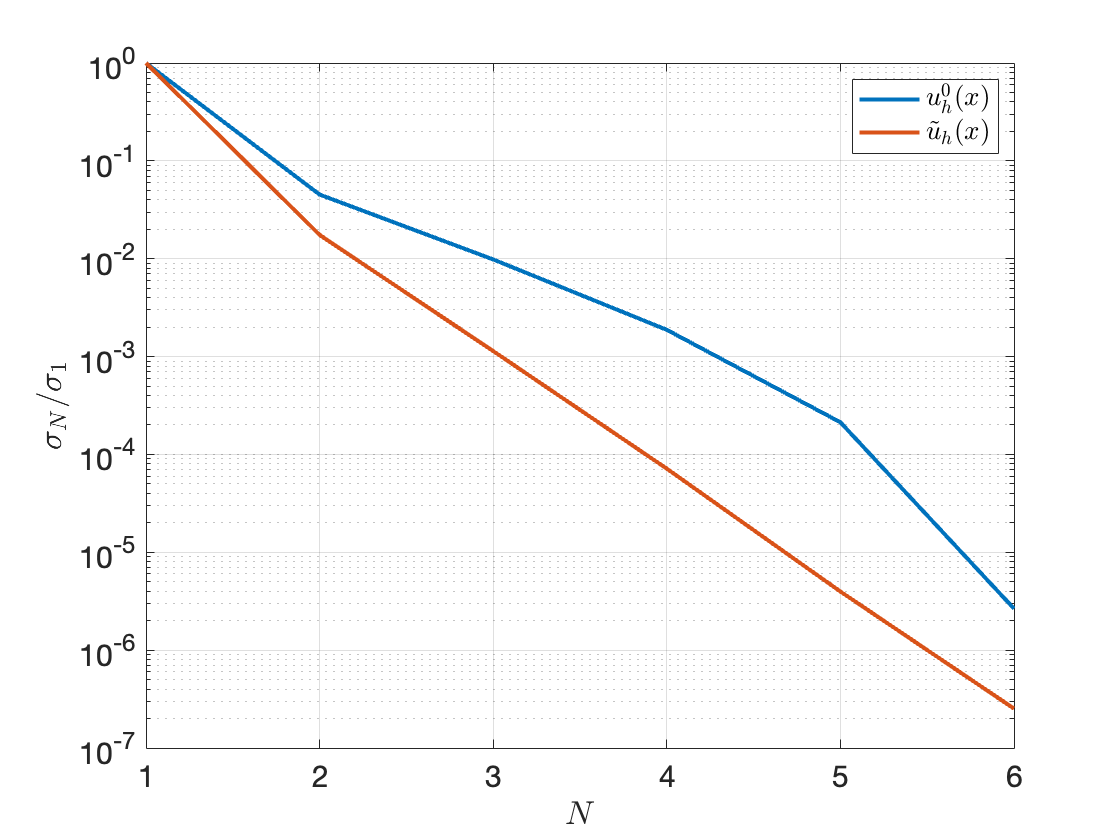}
    \caption{Singular value decay for the cylinder problem comparing fixed mesh solutions to mapped mesh solutions.}
    \label{fig:svdecaycyl}
\end{figure}

For model reduction, the training set is $\mathcal{S}_{\text{train}} = [2, 3, 4, 6, 8, 10]$. We note that the training set is slightly finer for lower Mach numbers. Changing the Mach number moves the shock more in the supersonic regime than in the hypersonic regime, so the parametric problem is actually more challenging on the lower end of the parameter range. Still, this is a relatively coarse training set for the Mach numbers considered. 

We see some evidence that the mapped snapshots are more amenable to model reduction in Figure \ref{fig:svdecaycyl} by plotting the decay of the singular values for the snapshots computed on a fixed mesh, $\bm u_h$, and those computed on an adapted mesh and pulled back to the reference mesh, $\tilde{\bm u}_h$. 
The normalized singular values of the mapped solutions are consistently lower for the solutions mapped back to the reference domain. The difference between the singular value decay is less pronounced than it would be for some other optimal transport-oriented ROMs, as alignment between snapshots is not specifically enforced. Still, the improvement is non-negligible and attained without any extra offline calibration. 

\begin{figure}[!htb]
\begin{center}
\subfigure
{
\includegraphics[width=0.45\textwidth]{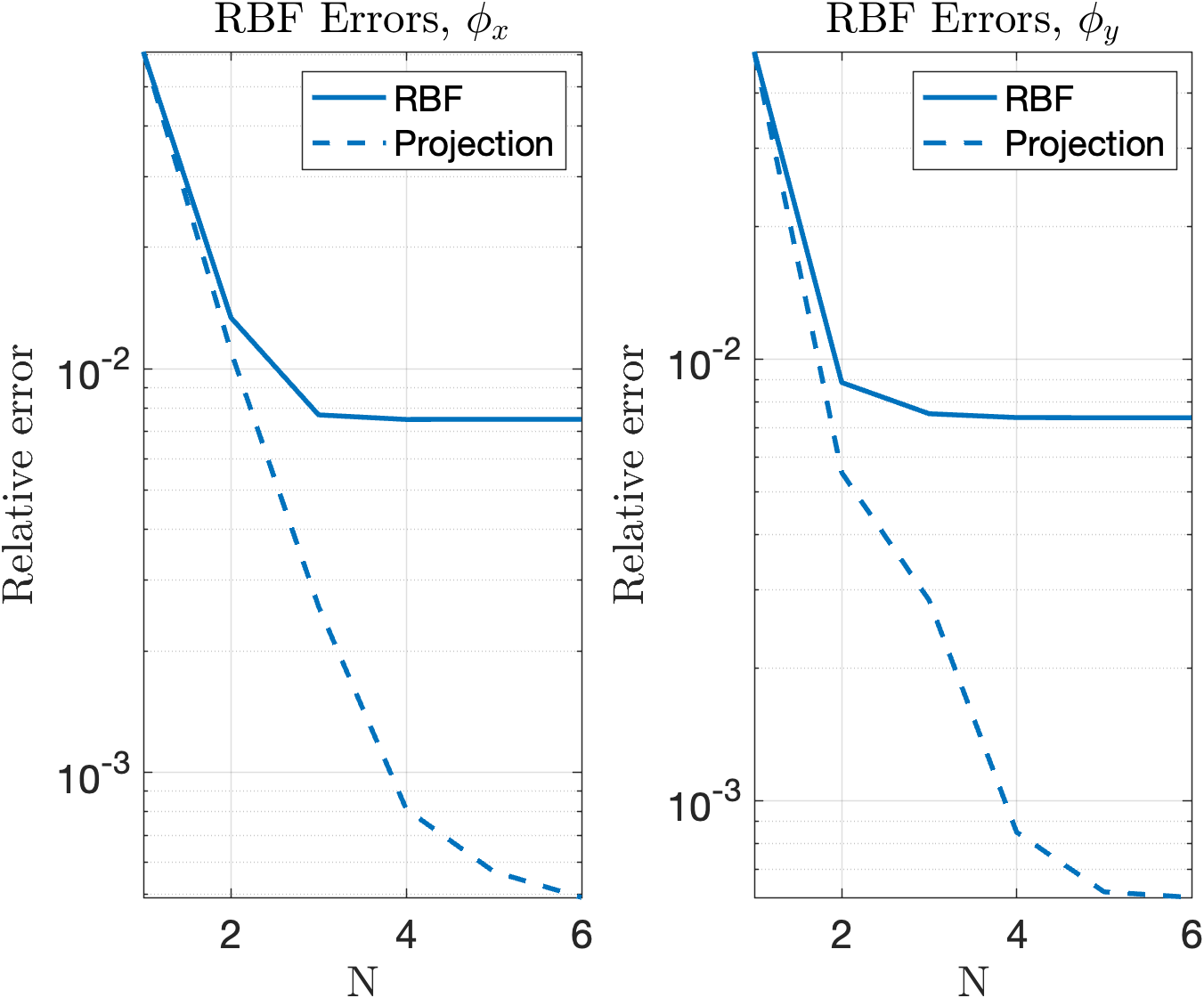}\label{fig:map_errs_cyl}
}
\subfigure
{
\includegraphics[width=0.45\textwidth]{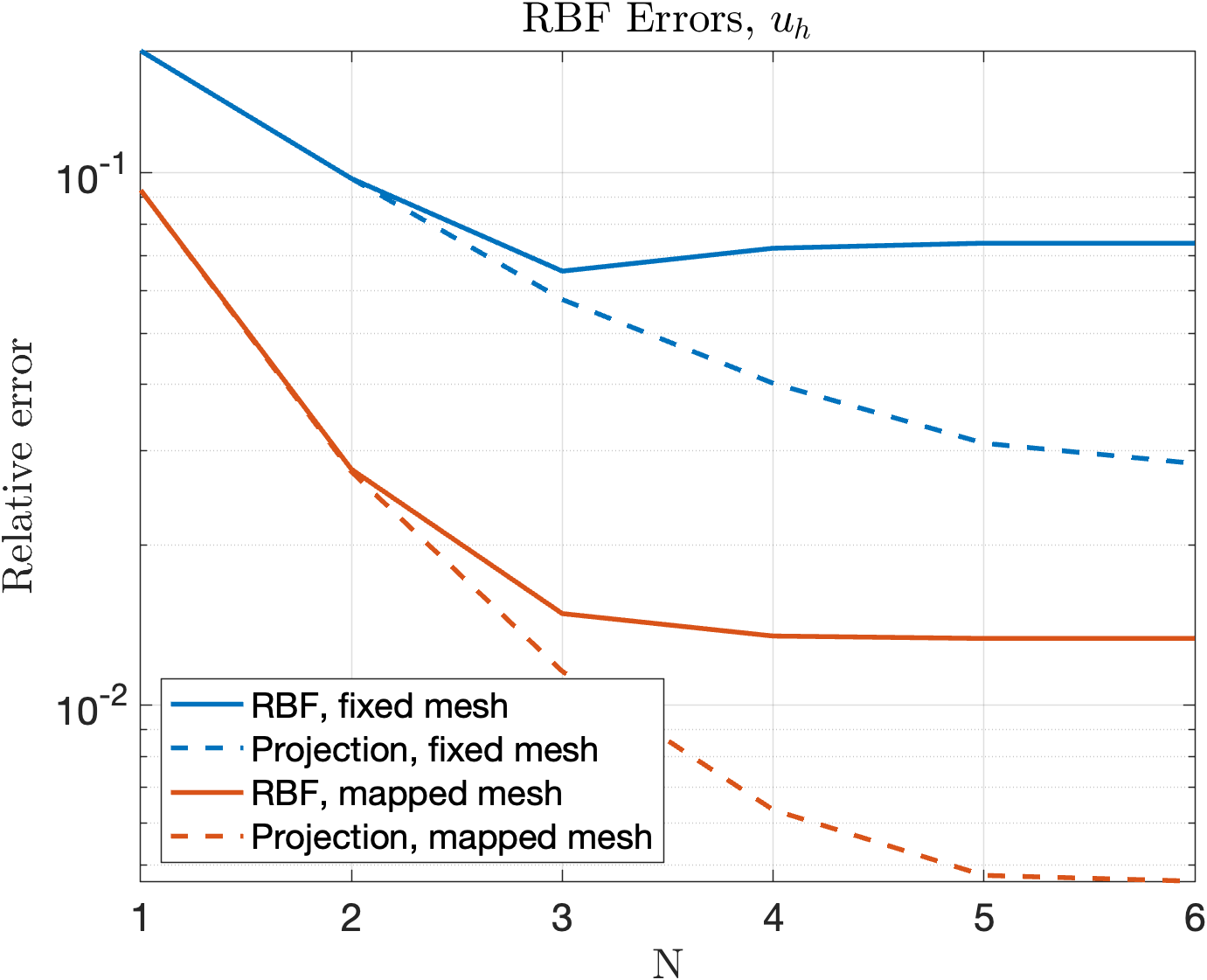}\label{fig:sol_errs_cyl}
} 
\caption{Average RBF ROM relative errors for mesh mapping $\bm \phi_h$ (left) and solution $\tilde{\bm u}_h$ (right) for high speed cylinder case.} \label{fig:cyl_errs}
\end{center}
\end{figure}

The ROMs are tested with a test set $\mathcal{S}_{\text{cyl}} \in \mathcal{D}$ consisting of three equally spaced test parameters between each training parameter, for a total of 15 test points.
The accuracy of the RBF ROM for the mappings and solution is shown in Figures \ref{fig:map_errs_cyl} and \ref{fig:sol_errs_cyl} for different choices of reduced basis size and summarized in Table \ref{table:tablecirc} for $N = 6$. The RBF ROM for the mapping is able to reach below 1\% relative error for both $x$ and $y$ components.  The RBF ROM shows improvements over the fixed mesh surrogate in terms of average errors on the fixed mesh, with neither surrogate reaching projection error.

\begin{table}[h!]
\begin{center}
    \begin{tabular}{|c||c|c|}
        \hline
        Quantity &  Average & Maximum\\
        \hline
        $E_{\bm u_h^0}$ & 0.078 & 0.149 \\
        \hline
        \hline
        $E_{\bm \tilde{\bm u}_h}$ & 0.011& 0.032 \\
        $E_{\bm \phi_h}$ & 0.008 & 0.018 \\
        \hline
    \end{tabular}
    \caption{Mean and max relative errors over $\bm \mu \in \mathcal{S}_{\text{cyl}}$ of fixed mesh ROM for $\bm u^0_h(\bm x)$ (top row) and mapped mesh ROM for $\tilde{\bm u}_h(\bm x)$ and $\bm \phi_h(\bm x)$ (bottom two rows) for flow over a cylinder, with no truncation ($N=6$).
    }
    \label{table:tablecirc}
\end{center}
\end{table}

The solution profiles themselves, taking into account the mapping, show a dramatic improvement. While the RBF surrogates on the fixed mesh solution are able to achieve reasonable L2 errors for a nonintrusive surrogate, they result in entirely unphysical solutions. 
On the other hand, the mapped ROM is able to produce solutions that qualitatively match the true solution, with only some visible differences in the outflow boundaries towards the extremes of the parameter range.

\begin{figure}[!htb]
\begin{center}
\subfigure[$\mathrm{Ma}_{\infty} = 2.25$]
{
\includegraphics[width=0.45\textwidth]{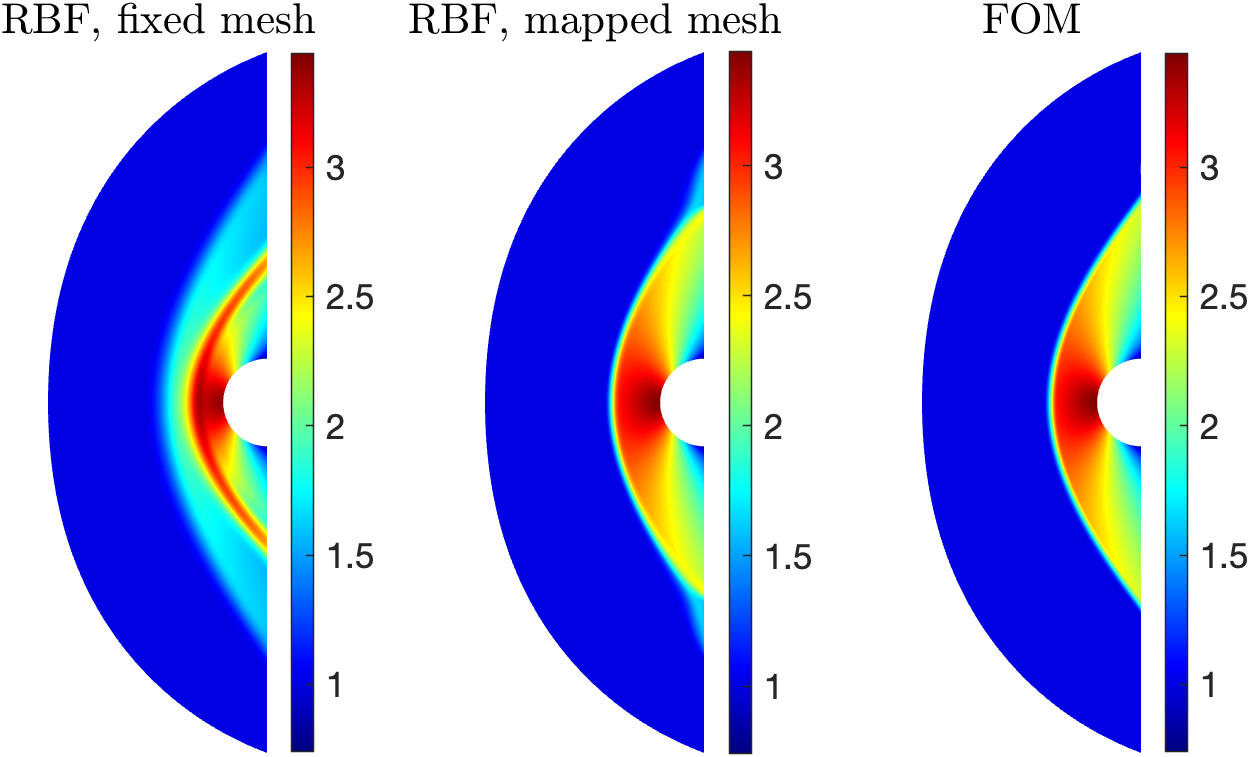}\label{fig:cyl_compare_2.25}}
\subfigure[$\mathrm{Ma}_{\infty} = 9.5$]
{
\includegraphics[width=0.45\textwidth]{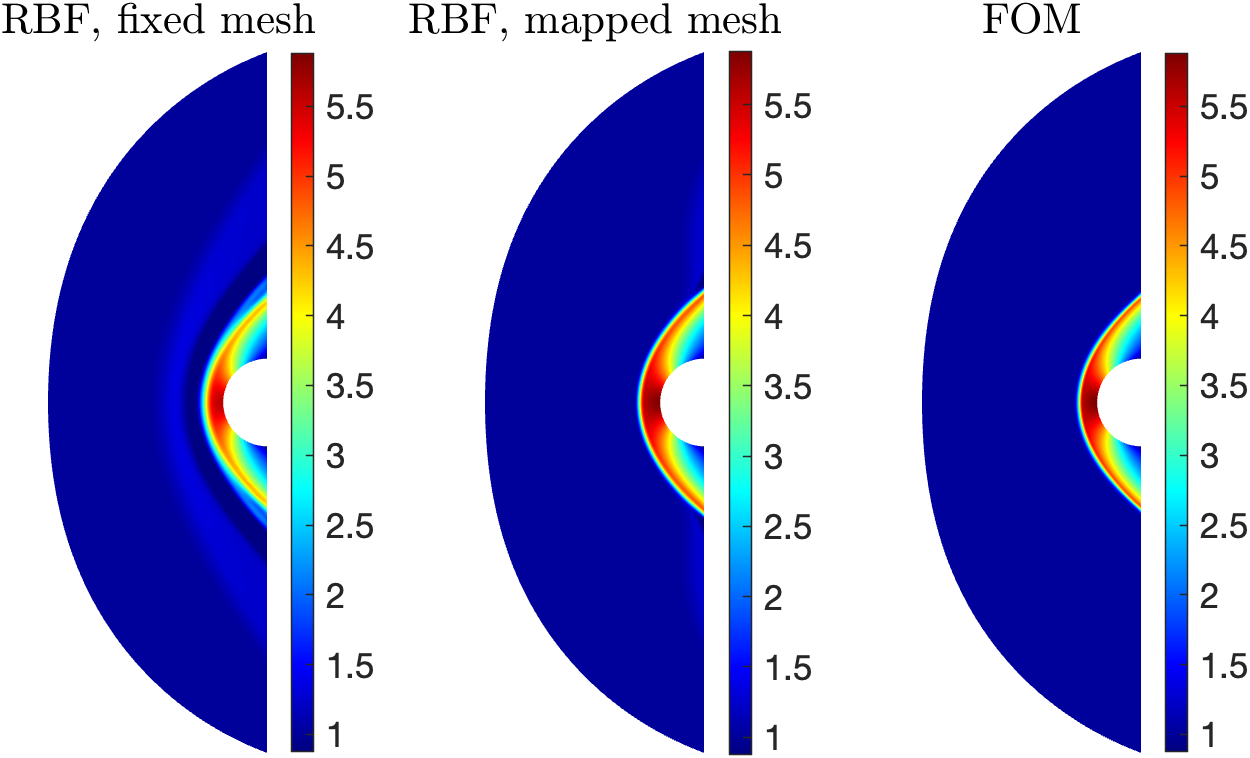}\label{fig:cyl_compare_9.5}}\\
\caption{Comparison of fixed mesh and mapped ROMs for parameters with the largest errors in the test set for flow around a cylinder. }\label{fig:cyl_rbf_compare}
\end{center}
\end{figure}
 
\begin{figure}
    \centering
    \includegraphics[width=0.35\textwidth]{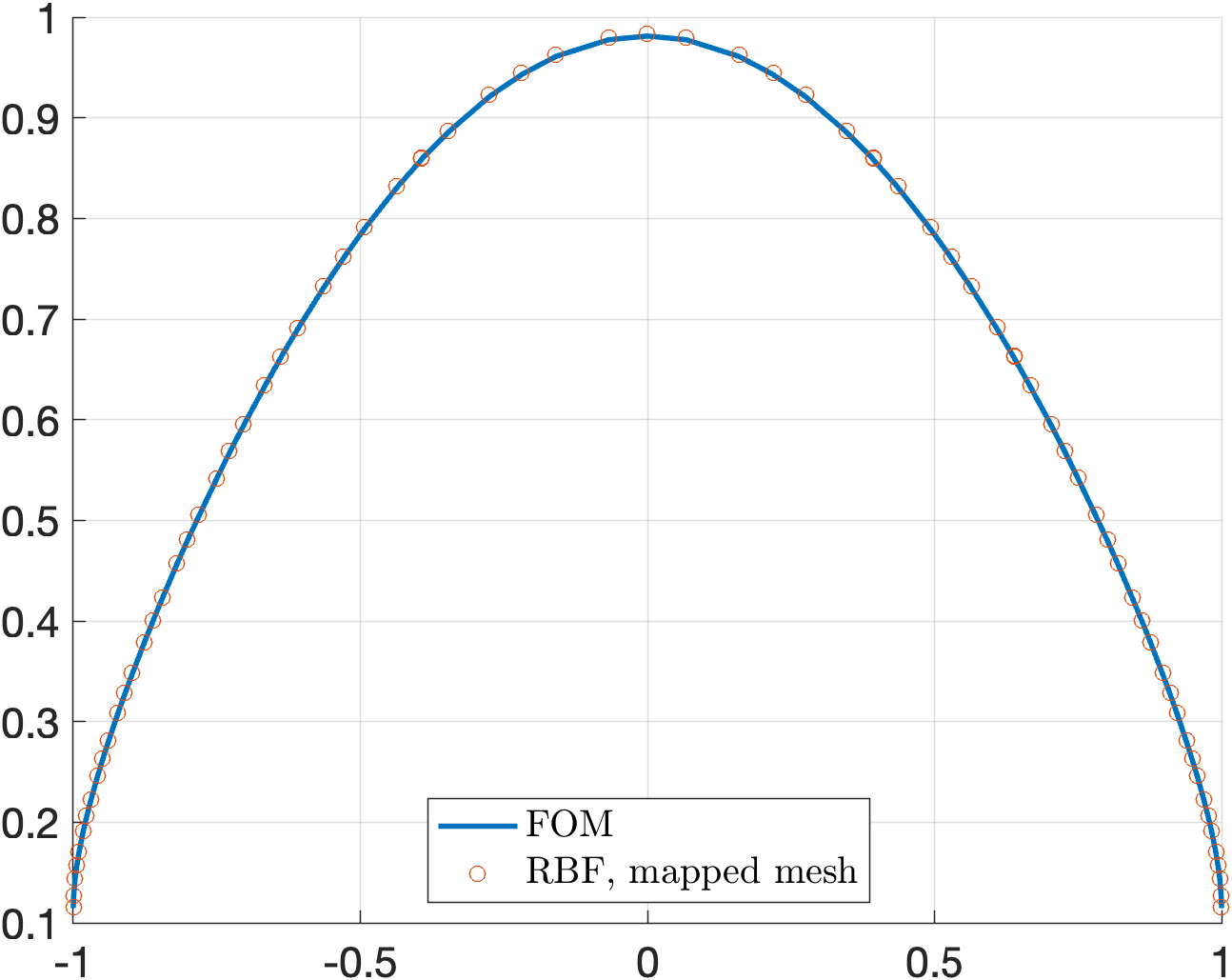}
    \caption{Nondimensionalized wall pressure with FOM and RBF ROM at $\mathrm{Ma}_{\infty} = 2.25$}
    \label{fig:wall}
\end{figure}

These irregularities at the outflow are not expected to impact wall quantities; in Figure \ref{fig:wall} we can see that wall pressure is indeed captured well by the ROM for the parameter with the largest errors, $\bm \mu = 2.25$.

\subsection{Flow over a smooth bump}
\begin{figure}[!htb]
\begin{center}
\subfigure[Mesh]
{
\includegraphics[width=0.46\textwidth]{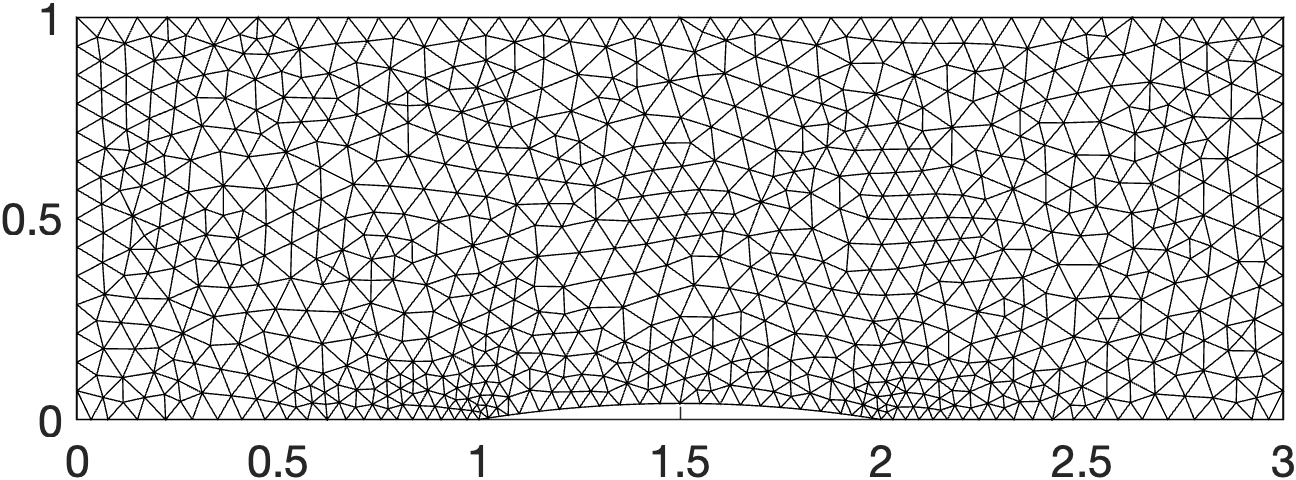}\label{fig:bump_mesh1}
}
\subfigure[Adapted mesh]
{
\includegraphics[width=0.46\textwidth]{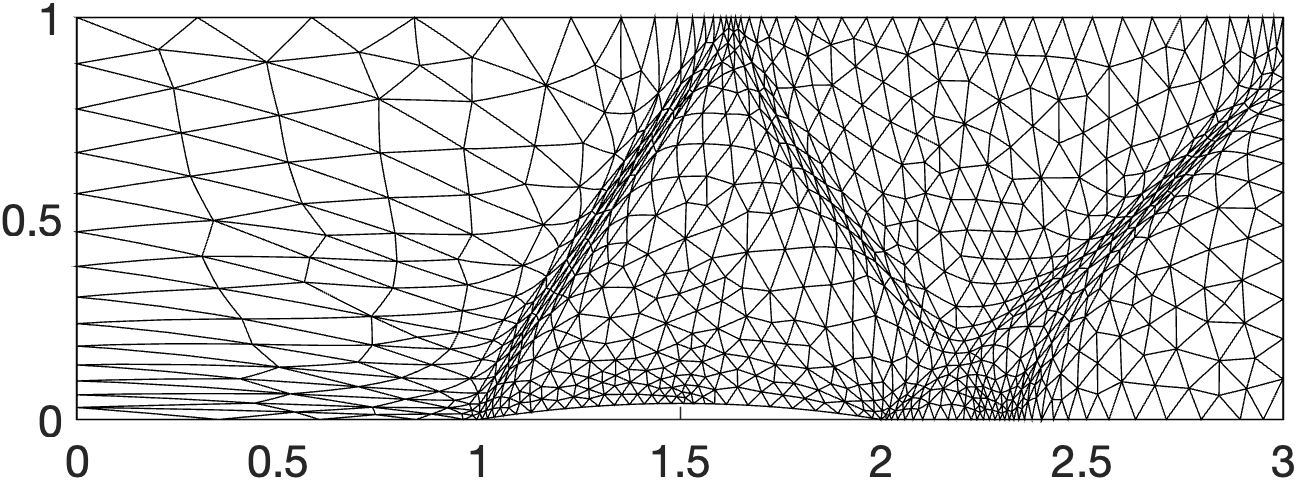}\label{fig:bump_meshadapt}
}\\
\subfigure[Physical density on starting mesh]
{
\includegraphics[width=0.47\textwidth]{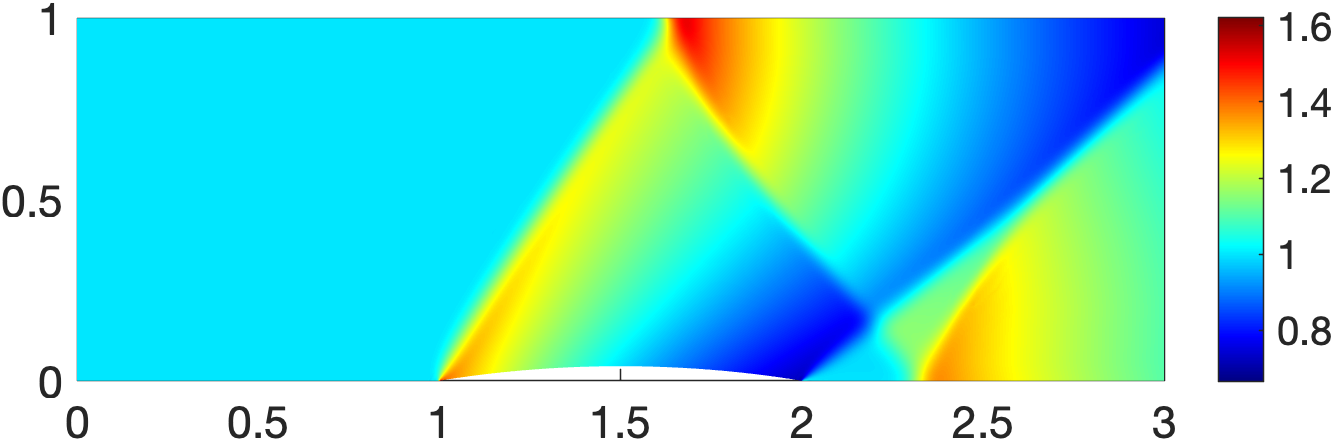}\label{fig:bump_noadapt}
}
\subfigure[Physical density on adaptive mesh]
{
\includegraphics[width=0.47\textwidth]{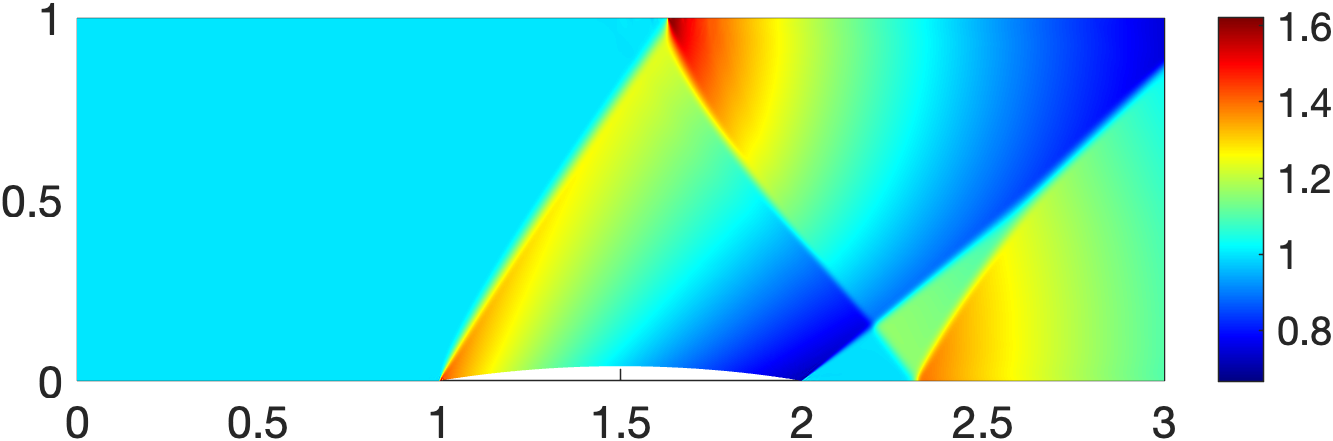}\label{fig:bump_adapt}
}\\
\subfigure[Physical density along $y = 0.7$]
{
\includegraphics[width=0.44\textwidth]{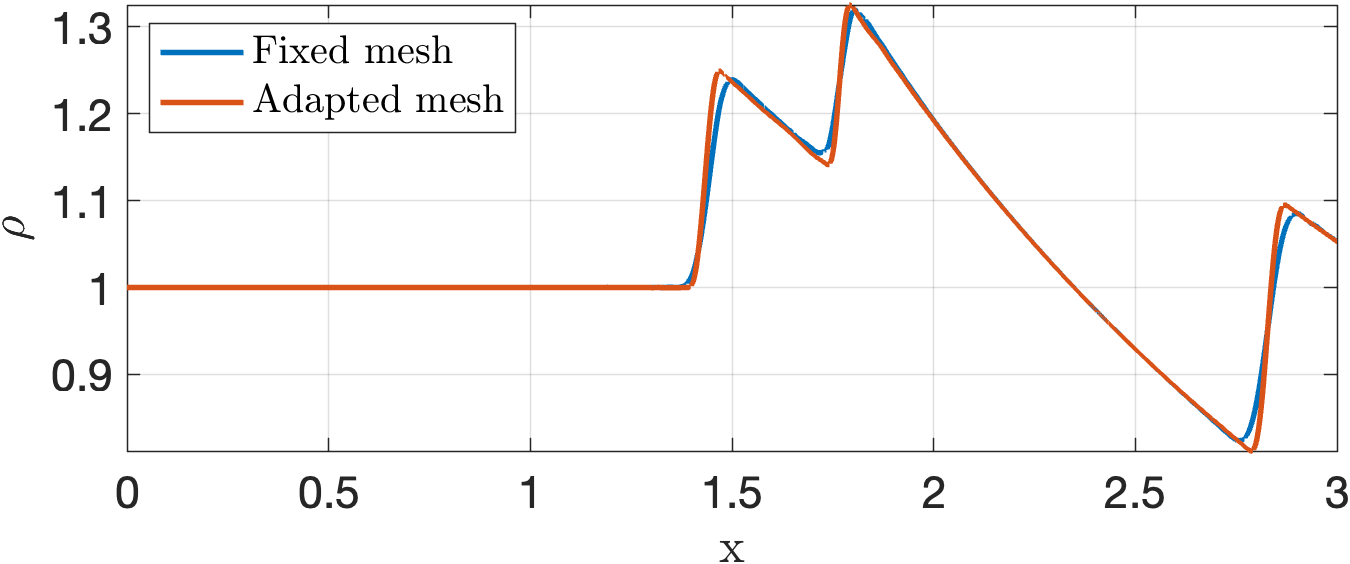}\label{fig:bump_slice_high}
}
\subfigure[Physical density along $y = 0.25$]
{
\includegraphics[width=0.44\textwidth]{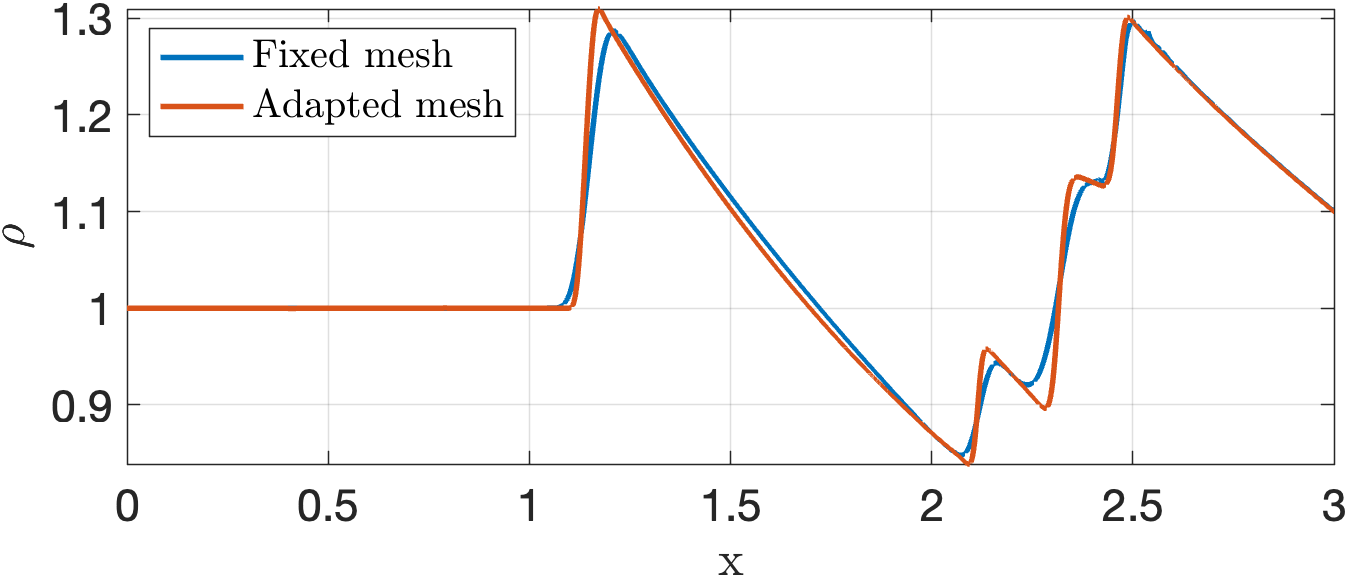}\label{fig:bump_slice_low}
}
\caption{Fixed mesh and adaptive mesh solutions of $\mathrm{Ma}_{\infty} = 1.4$ flow over the smooth bump.}\label{fig:bump_summary}
\end{center}
\end{figure}

We repeat the procedure on a flow with more complicated features. 
This flow over a smooth bump is similar to the problem considered in \cite{ferrero2022registration} and is a common benchmark for shock capturing schemes \cite{nguyen2011adaptive, burbeau2001problem, lynn1995multigrid}. 
We consider the case from \cite{nguyen2011adaptive} of supersonic flow over a smooth bump with  height 4\% in a rectangular domain with length 3 and height 1. Inlet and outlet conditions are set at the left and right boundaries while the top and bottom are inviscid walls. The geometry and initial mesh with 1048 elements is shown in Figure \ref{fig:bump_mesh1}. We use polynomial order 4 and viscosity continuation for all cases. 
The Monge-Ampère equation is solved on a rectangular domain and the resulting adapted grid is mapped onto the physical geometry with a simple 
transformation.
 
For $\mathrm{Ma}_{\infty} = 1.4$, the solution with and without the adapted mesh is shown in Figures \ref{fig:bump_noadapt} and \ref{fig:bump_adapt}. 
The mesh is able to refine towards the initial shocks, the reflections off of the top wall, and the resulting interactions, leading to sharper shocks throughout the domain.

For the parametric study, we let the Mach number vary from 1.4 to 1.5, which is a smaller parametric variation than in other studies.

\begin{figure}
    \centering
    \includegraphics[width=0.5\textwidth]{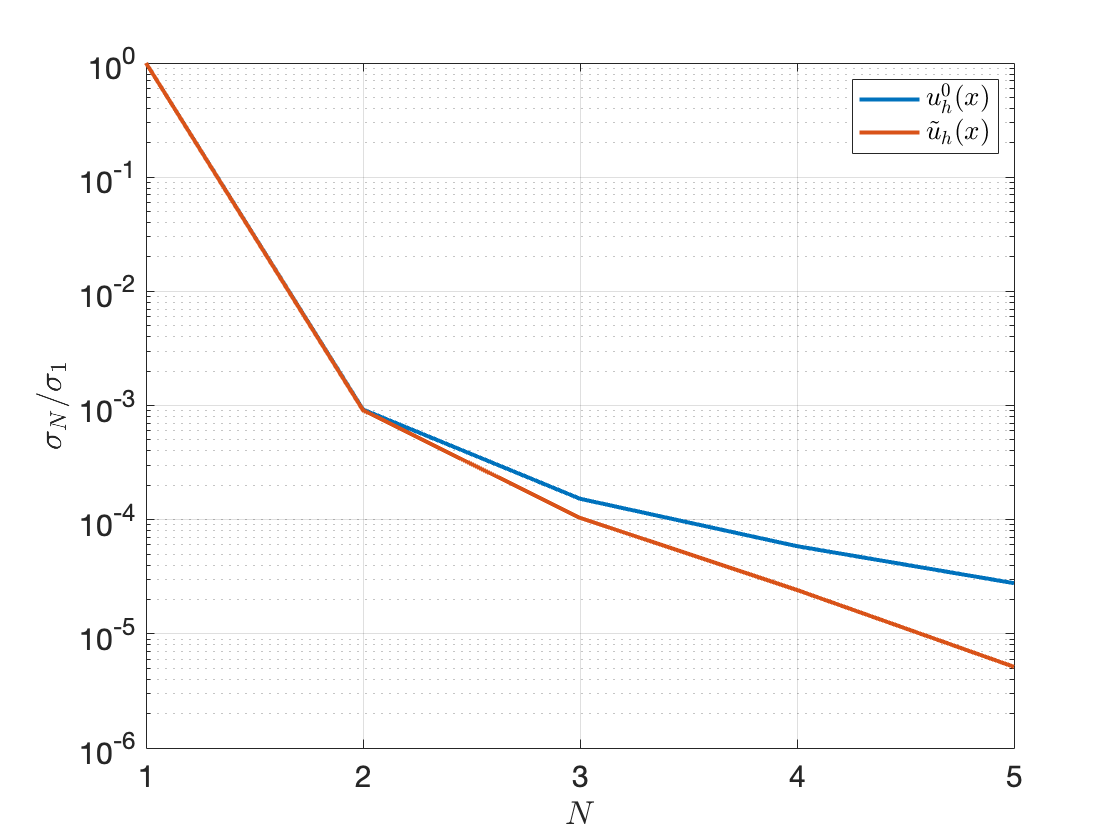}
    \caption{Decay of singular values for bump case with and without mesh mapping for training set $\mathcal{S}_{\text{train}}^2$}. 
    \label{fig:svalbump}
\end{figure}

\begin{figure}[!htb]
\begin{center}
\subfigure[$\rho$ for $\mathrm{Ma}_{\infty} = 1.4127$]
{
\includegraphics[width=0.48\textwidth]{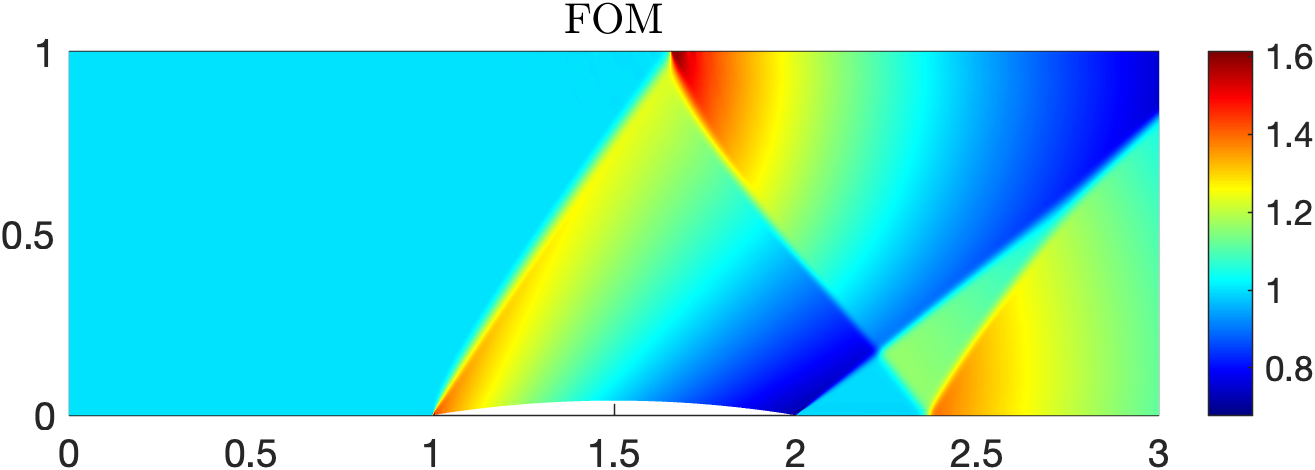}\label{fig:bump_compare_fom}
}\\
\subfigure[RBF ROMs trained on $\mathcal{S}_{\text{train}}^1$]
{
\includegraphics[width=0.46\textwidth]{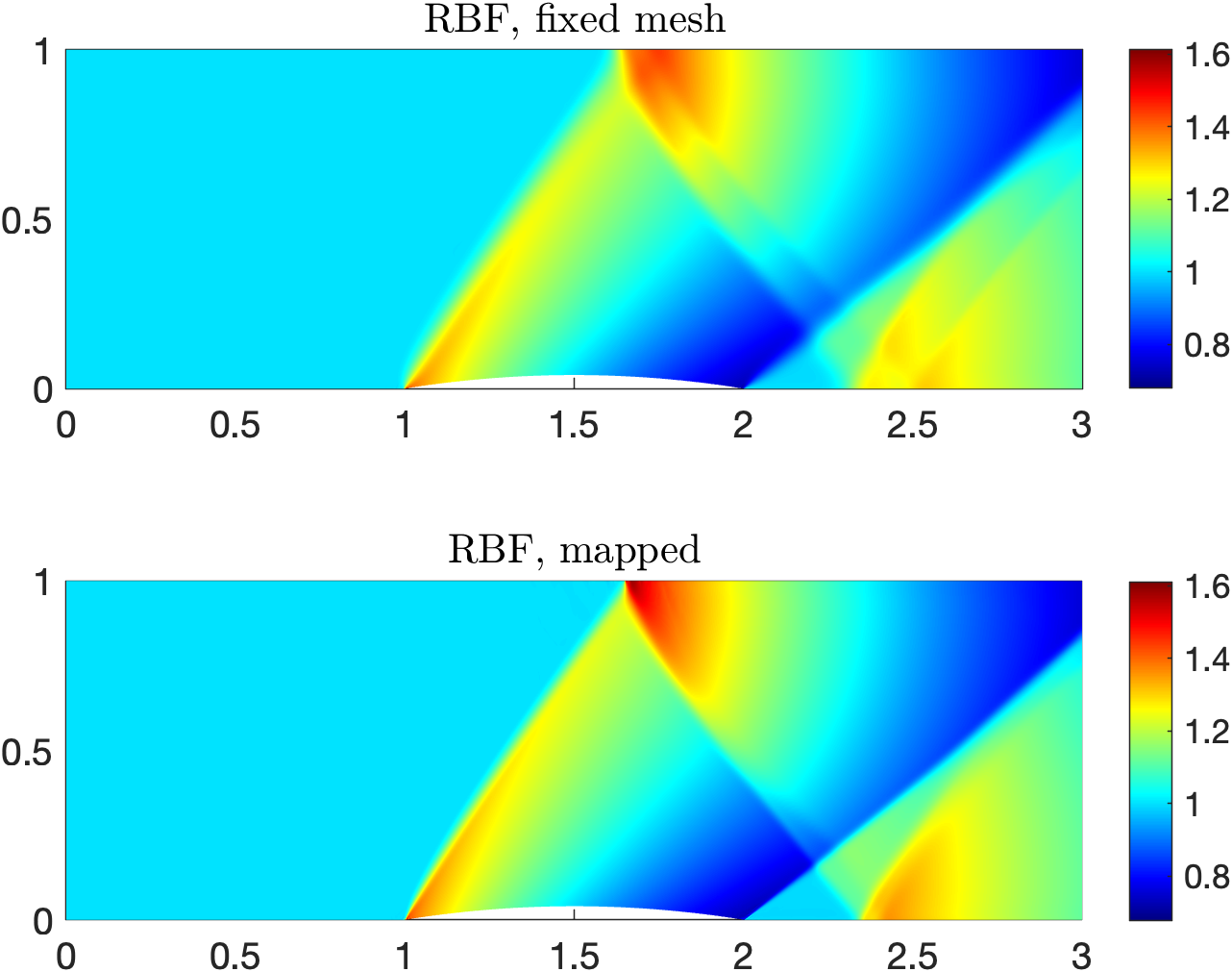}\label{fig:sparse_compare}
}
\subfigure[RBF ROMs trained on $\mathcal{S}_{\text{train}}^1$]
{
\includegraphics[width=0.46\textwidth]{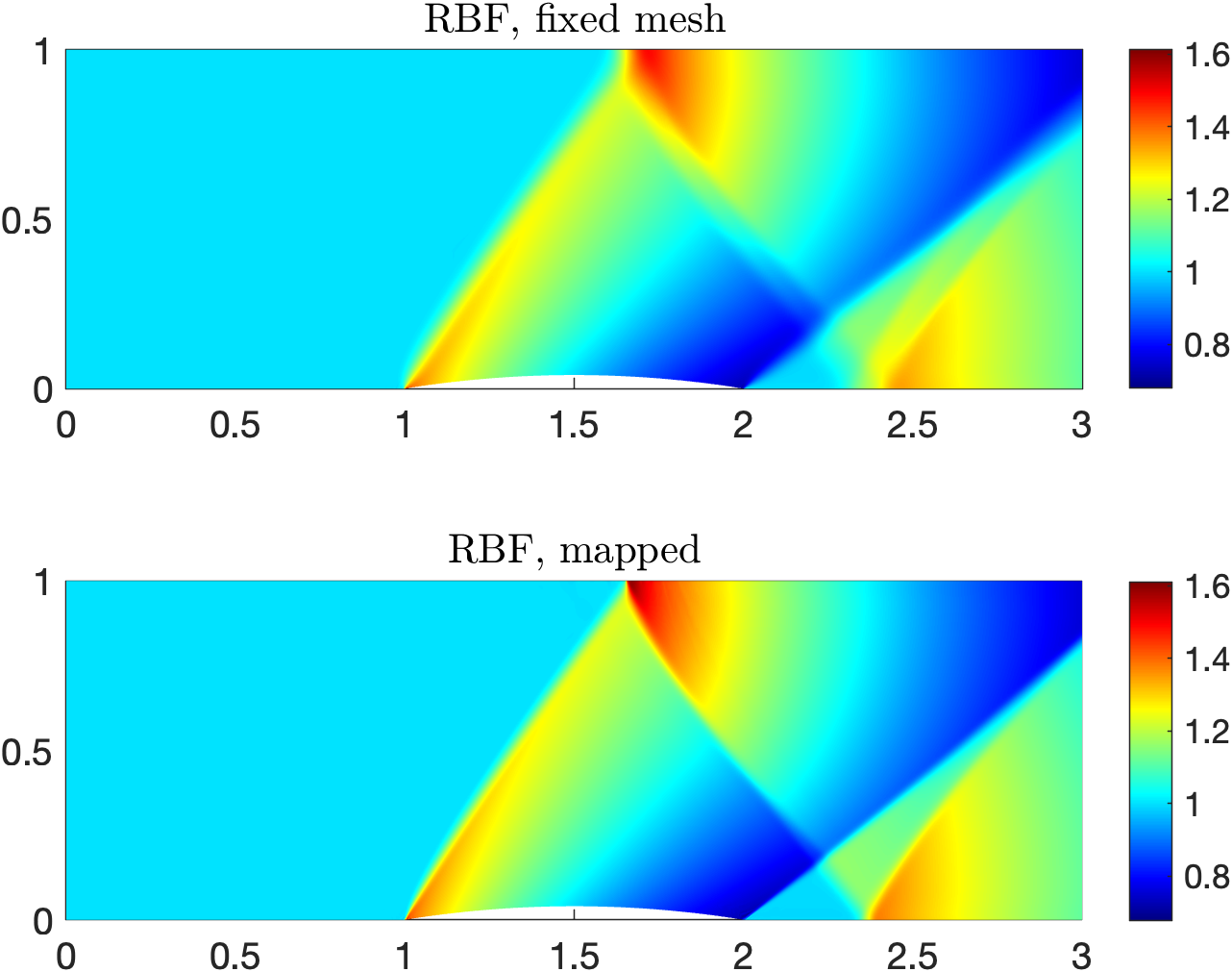}\label{fig:dense_compare}
}\\
\subfigure[$\rho$ at $y = 0.7$ for ROMs trained on $\mathcal{S}_{\text{train}}^2$]
{
\includegraphics[width=0.46\textwidth]{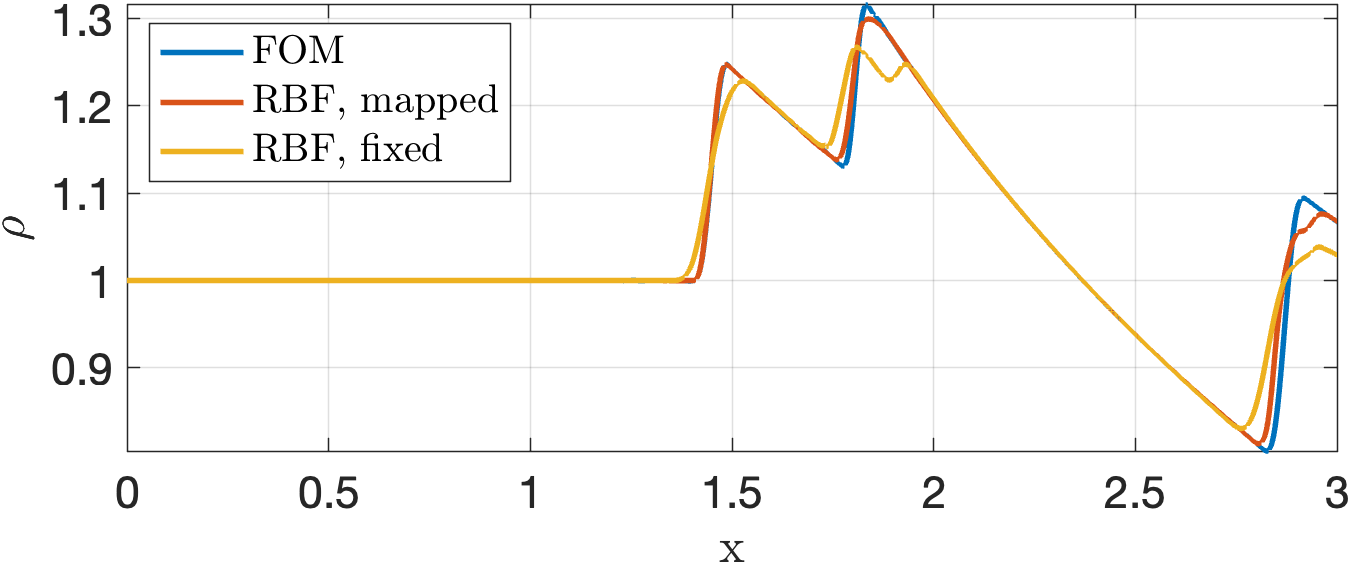}\label{fig:bump_clice2_mu2}
}
\subfigure[$\rho$ at $y = 0.7$ for ROMs trained on $\mathcal{S}_{\text{train}}^2$]{
\includegraphics[width=0.46\textwidth]{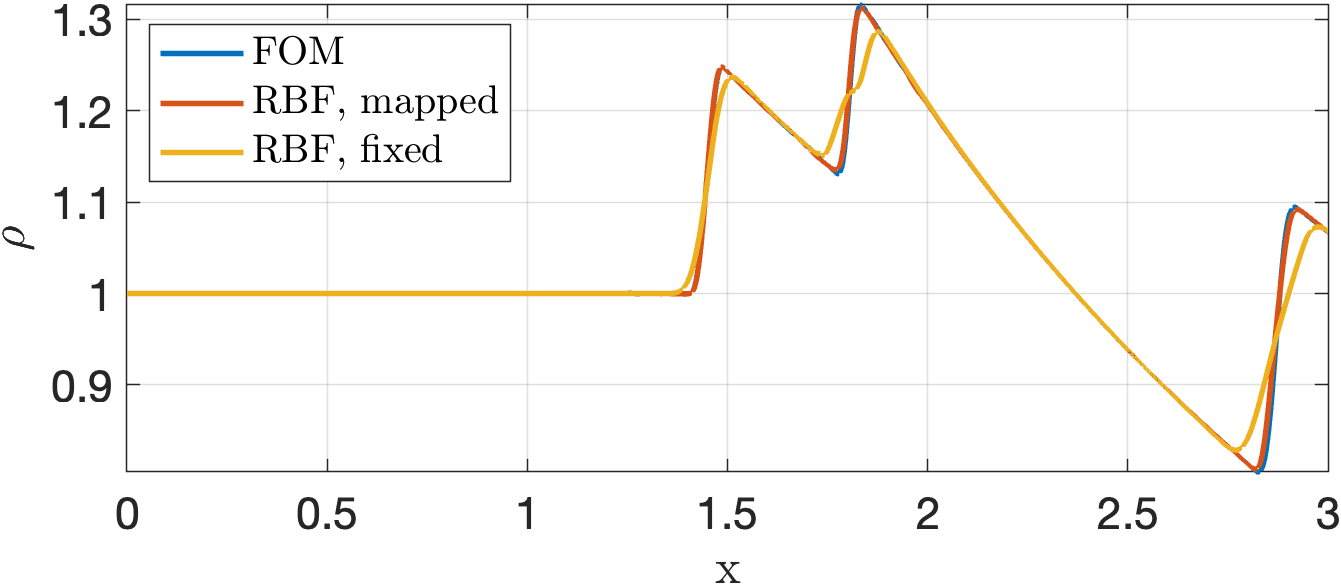}\label{fig:bump_slice1_mu2}
}\\
\subfigure[$\rho$ at $y = 0.25$ for ROMs trained on $\mathcal{S}_{\text{train}}^1$.]
{
\includegraphics[width=0.46\textwidth]{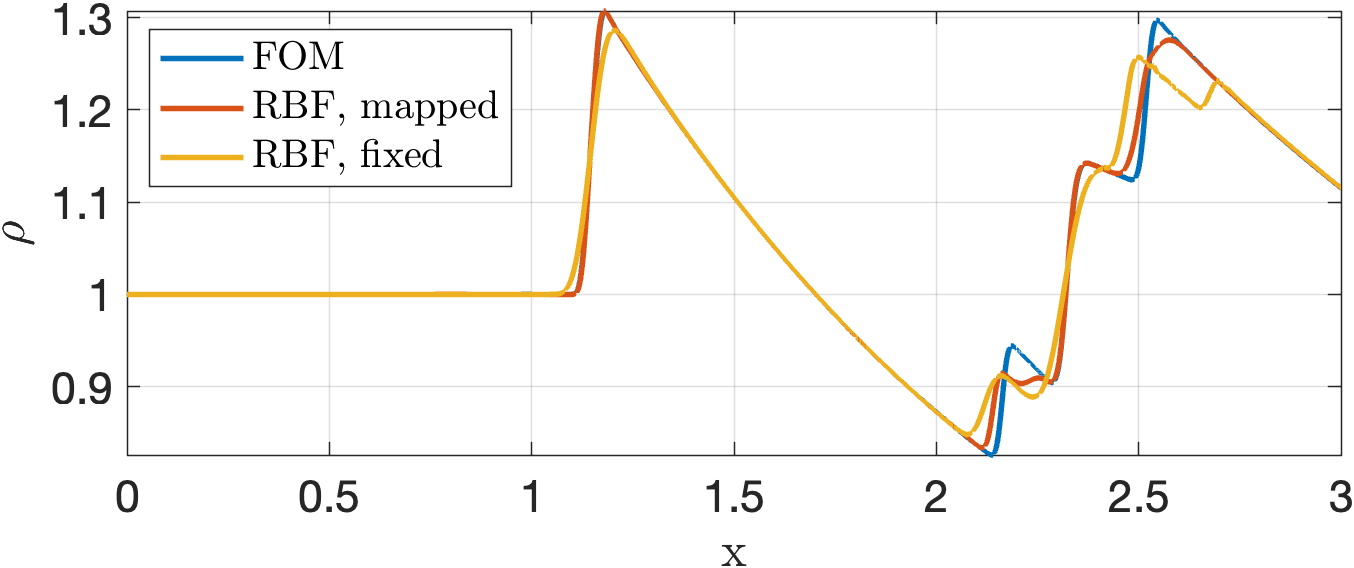}\label{fig:nump_slice2_mu1}
}
\subfigure[$\rho$ at $y = 0.25$ for ROMs trained on $\mathcal{S}_{\text{train}}^2$.]{
\includegraphics[width=0.46\textwidth]{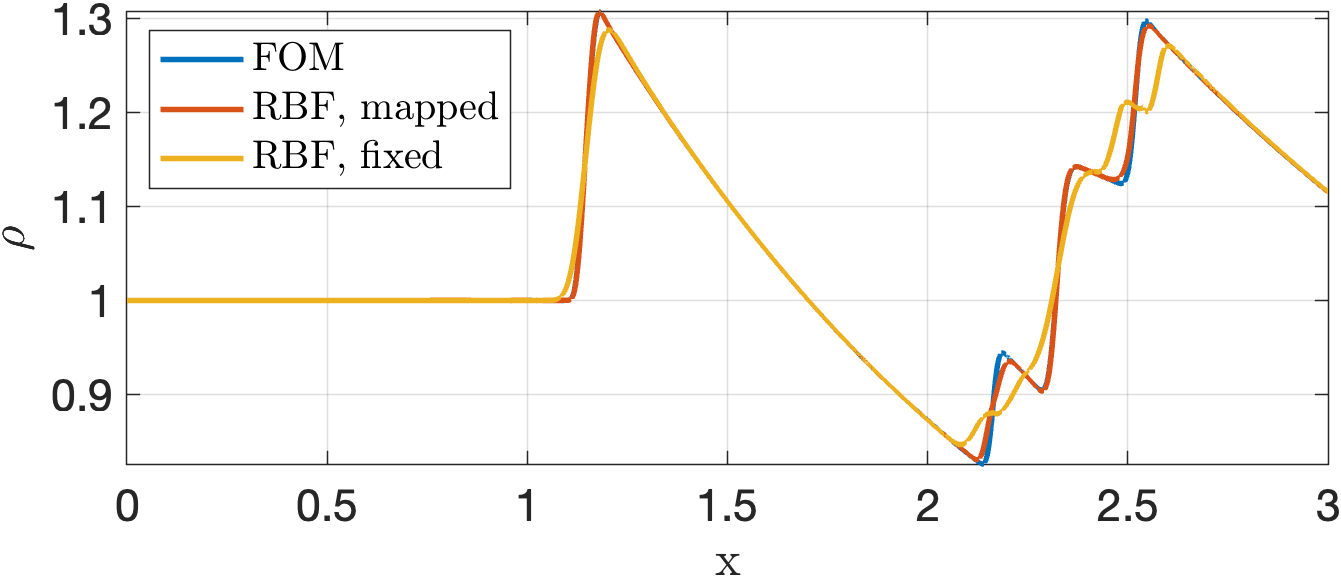}\label{fig:bump_slice2_mu2}
}
\caption{ROM solution profiles for $\bm \mu$ with the largest error, trained with different parameter sets.}\label{fig:bump_rbf_errs}
\end{center}
\end{figure}

We consider two sets of training data of $\mathcal{S}_{\text{train}}^1 = [1.4, 1.45, 1.5]$ and $\mathcal{S}_{\text{train}}^2 = [1.4, 1.425, 1.45, 1.475, 1.5]$.  
The performance of the ROM is tested at 10 random points in this Mach number range. The average and max errors are given in Table \ref{table:tablebump}. Although there is a consistent reduction in the error, the improvements are much more modest than in the cylinder case. 
As can be seen in the singular value decay in Figure \ref{fig:svalbump} for $\mathcal{S}_{\text{train}}^2$, the representability of the solution in a linear basis does not change as much as for the cylinder case.  
Despite this, in Figure \ref{fig:bump_rbf_errs}, we again see that the qualitative nature of the solution is improved with the mapped ROMs. 
There are fewer noticeable ``staircase''-like patterns that are commonly observed in linear model reduction of localized features, although some features are not accurately captured with the sparser training data set.  

\begin{table}[h!]
\begin{center}
    \begin{tabular}{|c||c|c|}
        \hline
        Quantity &  Average & Maximum\\
        \hline
        $E_{\bm u_h^0}$ & 0.008 & 0.009 \\
        \hline
        \hline
        $E_{\bm \tilde{\bm u}_h}$ & 0.002& 0.003 \\
        $E_{\bm \phi_h}$ & 0.004 & 0.005 \\
        \hline
    \end{tabular}
    \caption{Mean and max relative errors over $\bm \mu \in \mathcal{S}_{\text{cyl}}$ of fixed mesh ROM for $\bm u^0_h(\bm x)$ (top row) and mapped mesh ROM for $\tilde{\bm u}_h(\bm x)$ and $\bm \phi_h(\bm x)$ (bottom two rows) for supersonic flow over a bump. Trained on $\mathcal{S}_{\text{train}}^2$ with no truncation ($N=5$).
}
\label{table:tablebump}
\end{center}
\end{table}

\subsection{Double wedge}

\begin{figure}[!htb]
    \centering
    \includegraphics[width=0.3\textwidth]{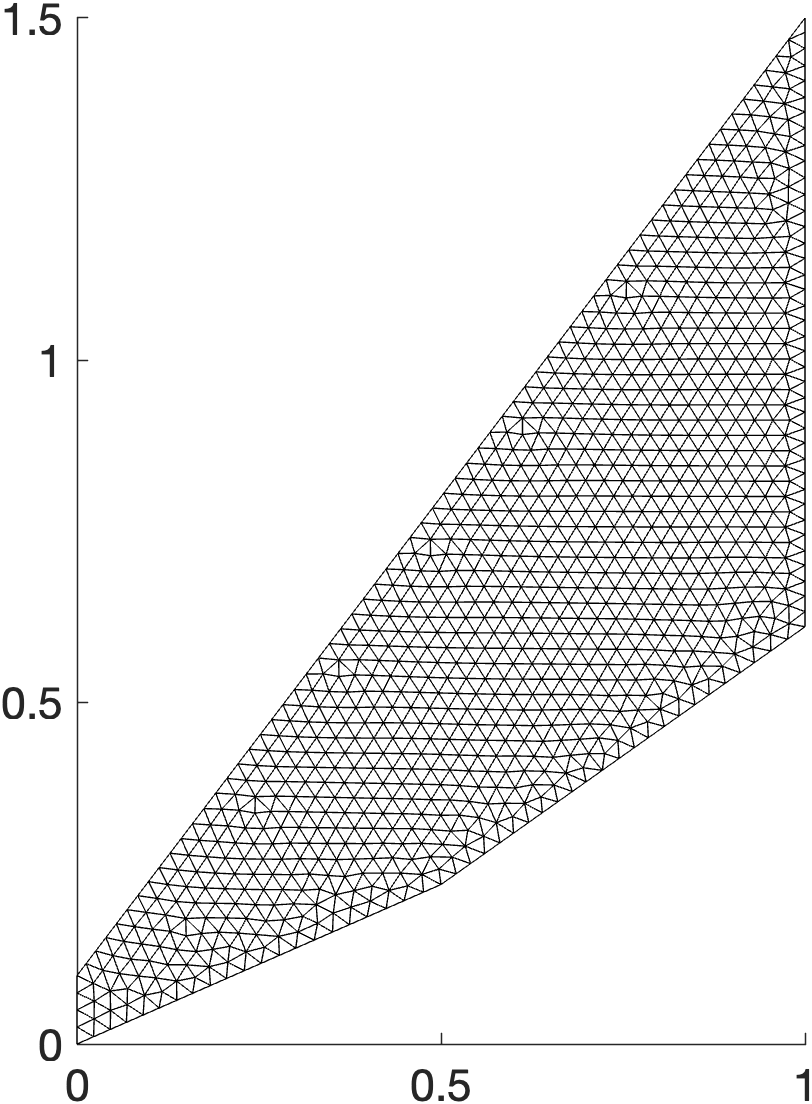}
    \caption{Starting unstructured mesh for double wedge problem}
    \label{fig:wedge_mesh}
\end{figure}

The final flow setup is taken from \cite{carnes2019code} and is meant to be an initial step towards double ramp and cone flows. The geometry is a 2D double wedge with angles $\theta_1 = 25^{\circ}$ and $\theta_2 = 37^{\circ}$. See Figure \ref{fig:wedge_mesh} for the starting mesh. 
The bottom wall is assigned an inviscid wall boundary condition, the rightmost wall is assigned an outlet boundary condition, and the other boundaries are treated as supersonic inlets. 
We use an initial unstructured mesh with 1842 tetrahedral elements with polynomial order 3. The FOM is solved with the viscosity continuation starting from uniform free-stream values. 

\begin{figure}[!htb]
\begin{center}
\subfigure[$\mathrm{Ma}_{\infty} = 4$]
{
\includegraphics[width=0.6\textwidth]{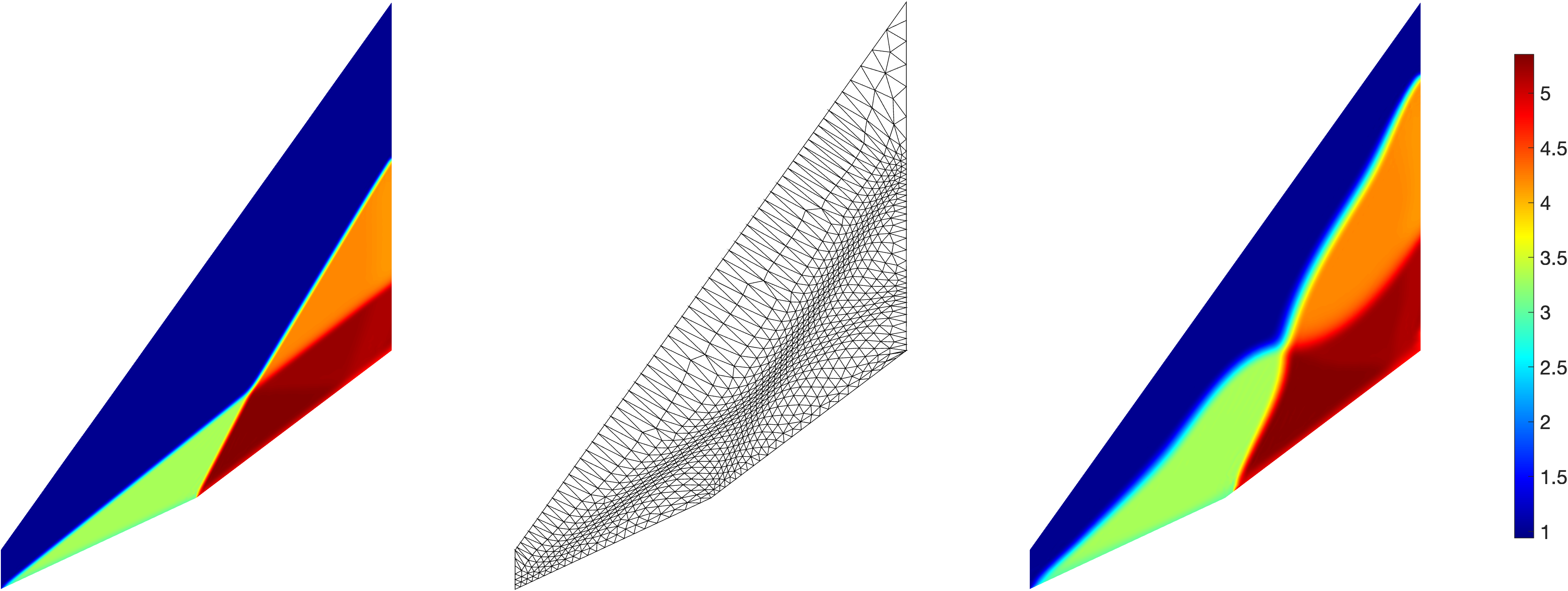}\label{fig:wedge4}}
\subfigure[$\mathrm{Ma}_{\infty} = 5$]
{
\includegraphics[width=0.6\textwidth]{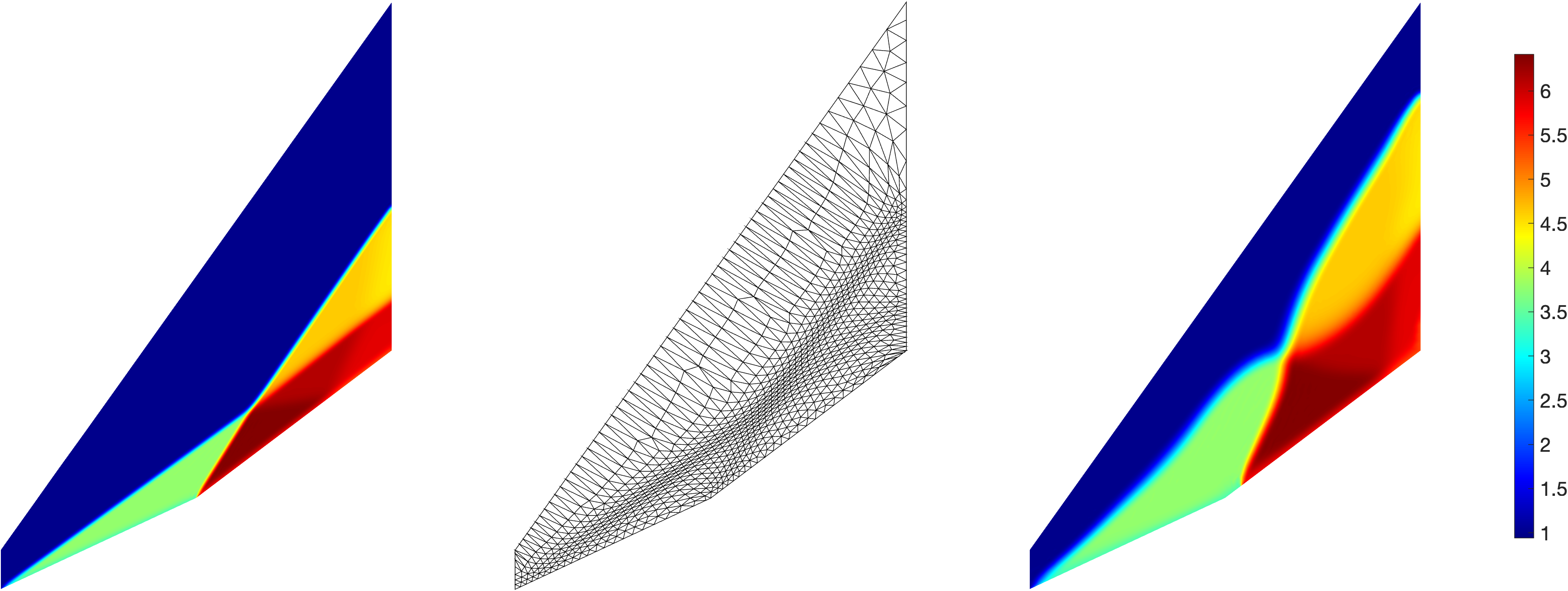}\label{fig:wedge55}}
\subfigure[$\mathrm{Ma}_{\infty} = 7$]
{\includegraphics[width=0.6\textwidth]{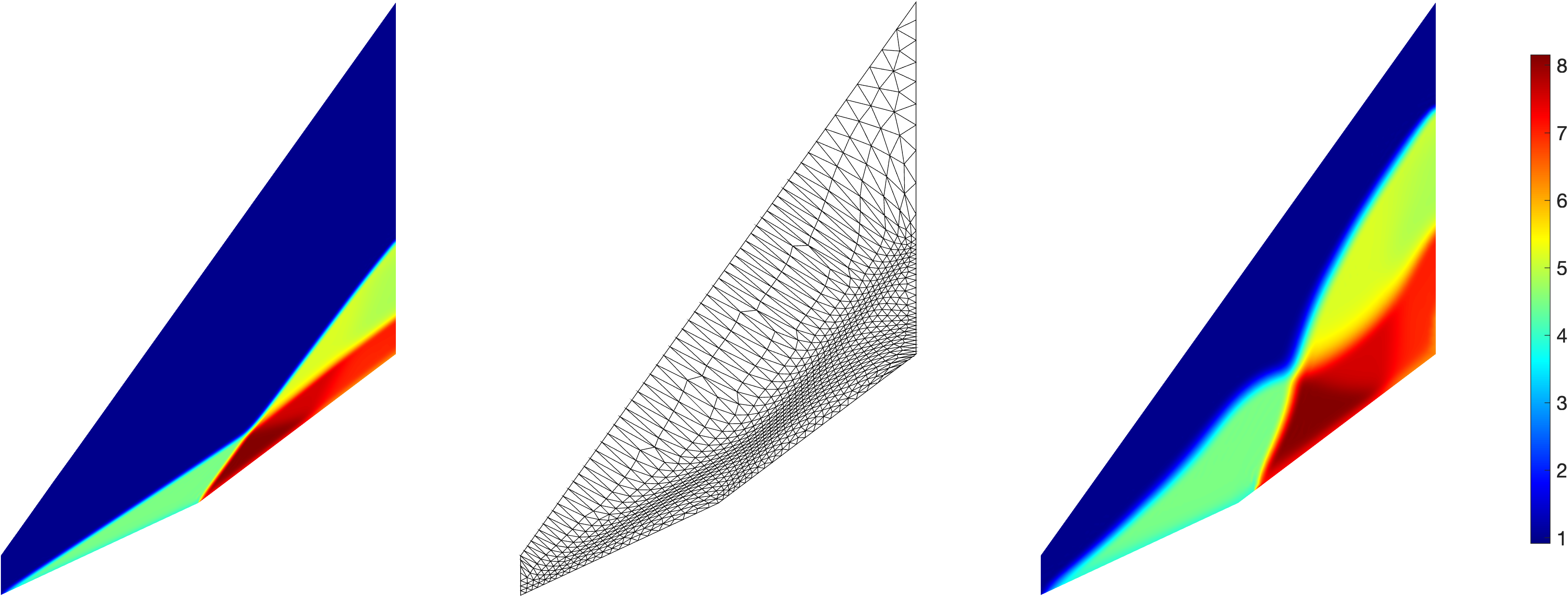}\label{fig:wedge7}} 
\caption{For different Mach numbers, a visualization of the solution $\bm u_h(\bm x)$, the corresponding adapted mesh $\mathcal{T}_h^{\bm \mu}$ and the solution mapped back to the reference mesh $\tilde{\bm u}_h$. Plotted quantity is density.}\label{fig:wedge_sweep}
\end{center}
\end{figure}

The ROM is trained on $\mathcal{S}_{\text{train}} = [4, 5, 6, 7]$ and tested at 12 equally spaced Mach numbers between the training set. Visualizations of the FOM and the mesh mapping can be seen in Figure \ref{fig:wedge_sweep}. 
As before, the flow features are brought into closer alignment on $\Omega$ and the features are aggressively smoothed out.

\begin{figure}
    \centering
    \includegraphics[width=0.5\textwidth]{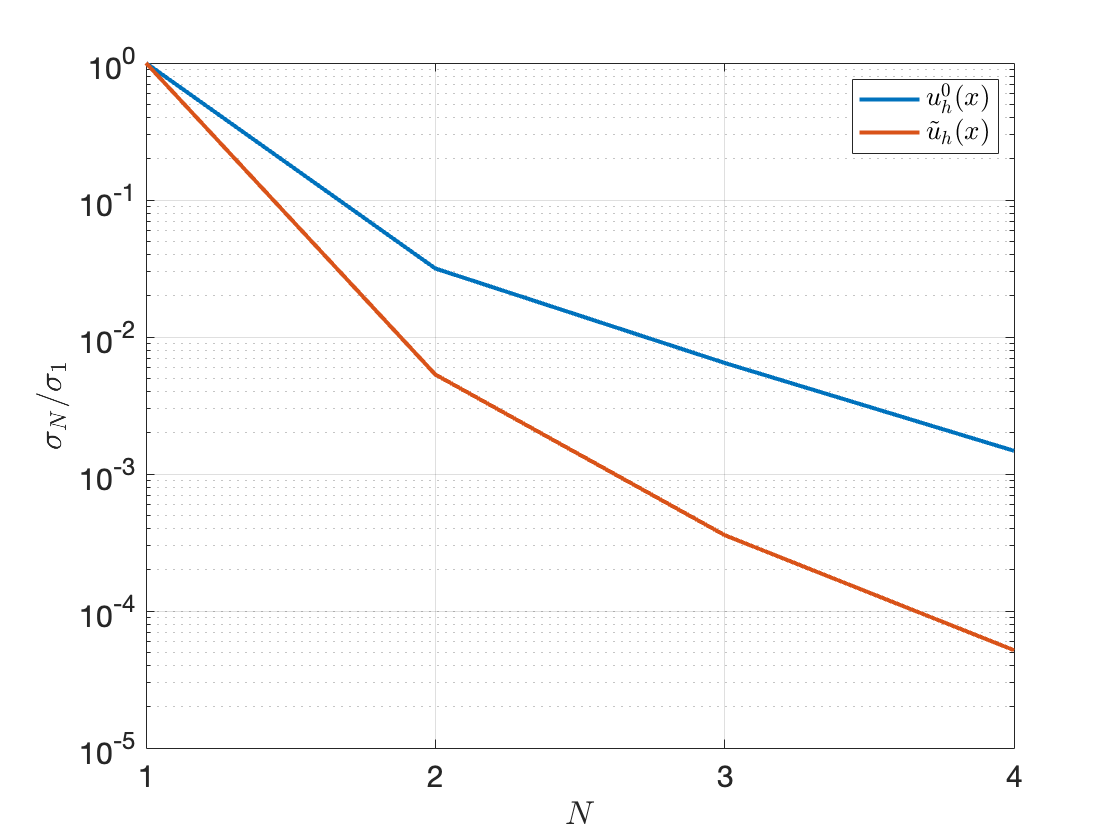}
    \caption{Decay of singular values for double wedge case with and without mesh mapping} 
    \label{fig:sval_wedge}
\end{figure}

We examine the decay of the singular values for the original snapshots and the mapped snapshots to see if this makes the ensemble more amenable to model reduction. Figure \ref{fig:sval_wedge} shows that the singular values of the mapped solutions are consistently about an order of magnitude smaller than those of the fixed mesh solutions. 

\begin{figure}[!htb]
\begin{center}
\subfigure[Mapping errors]
{
\includegraphics[width=0.45\textwidth]{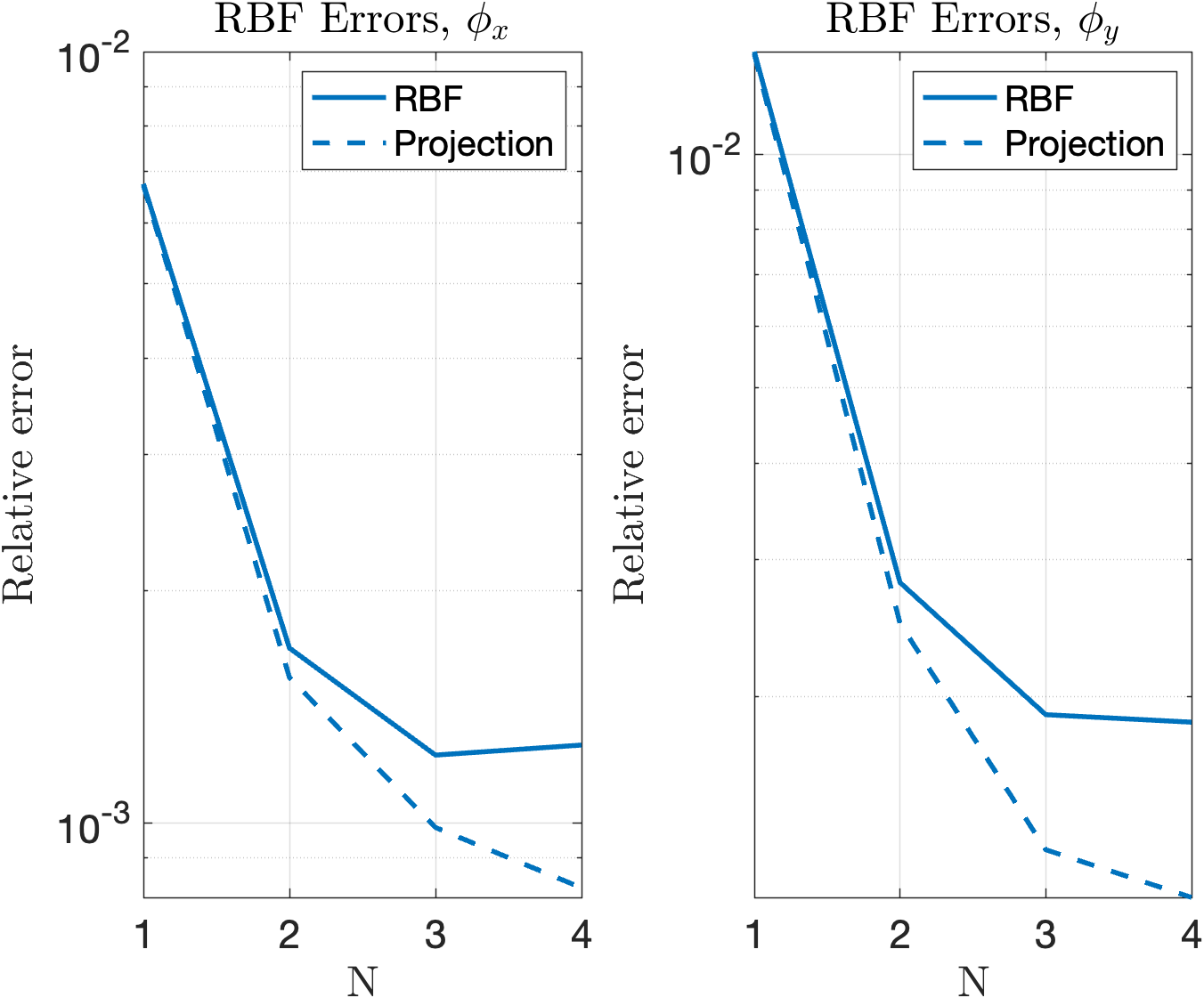}\label{fig:map_errs_wedge}
}
\subfigure[Solution errors]
{
\includegraphics[width=0.45\textwidth]{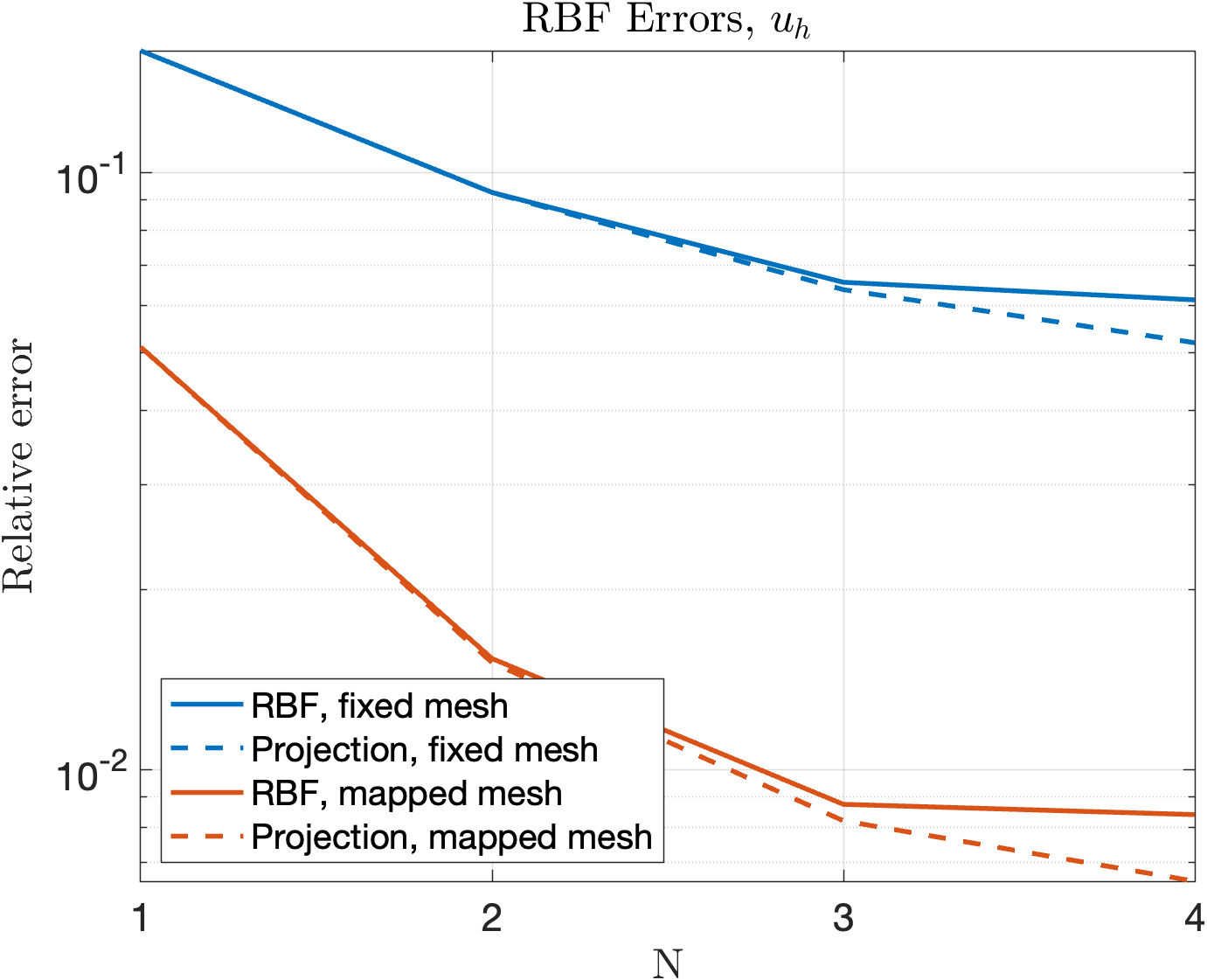}\label{fig:sol_errs_wedge}
} 
\caption{RBF ROM relative errors for mesh mapping $\bm \phi_h$ and solution $\tilde{\bm u}_h$}\label{fig:errs_wedge}
\end{center}
\end{figure}

\begin{figure}
    \centering
    \includegraphics[width=1.0\textwidth]{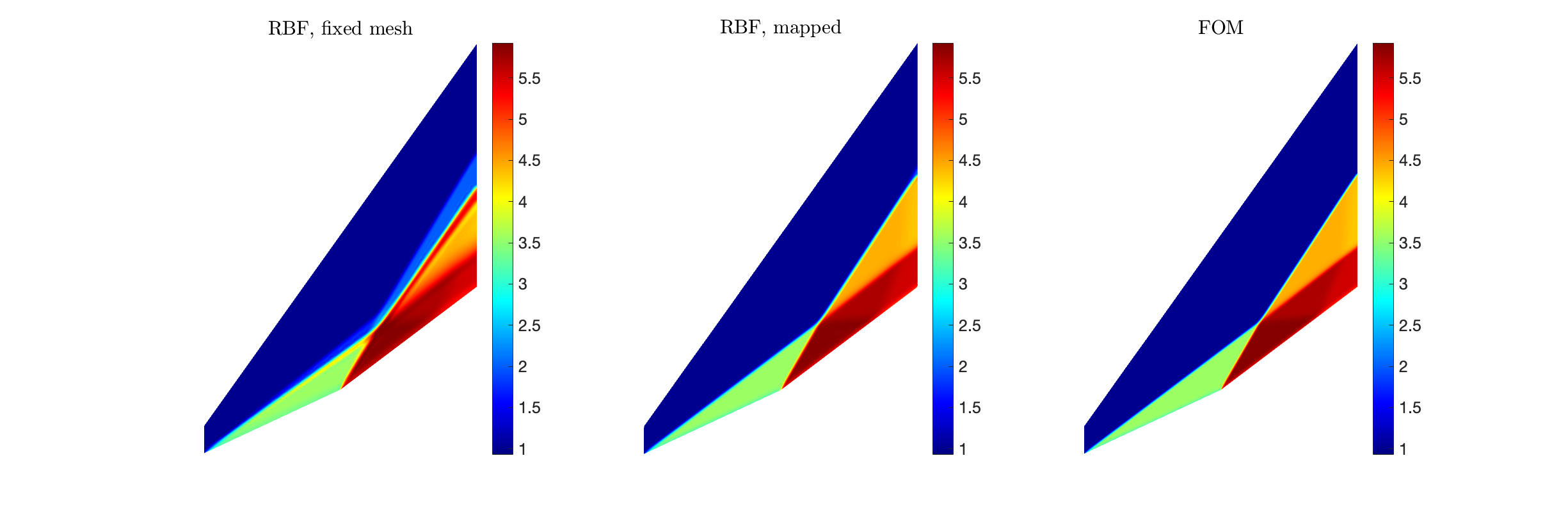}
    \caption{Comparison of the fixed mesh and mapped mesh ROMs for $\mathrm{Ma}_{\infty} = 4.5$}
    \label{fig:wedge_rbf}
\end{figure}

The RBF ROM trained on the mapped values is able to produce much more accurate results, as shown in Table \ref{table:ramp}. We see in Figure \ref{fig:errs_wedge} that the average solution errors are substantially smaller for all values of $N$ and that the relative mapping errors are also consistently below 1\%. Figure \ref{fig:wedge_rbf} shows again that the RBF ROM with the mesh mapping tends to do a better job producing physically meaningful solutions.
\begin{table}[h!]
\begin{center}
    \begin{tabular}{|c||c|c|}
        \hline
        Quantity &  Average & Maximum\\
        \hline
        $E_{\bm u_h^0}$ & 0.061 & 0.116 \\
        \hline
        \hline
        $E_{\bm \tilde{\bm u}_h}$ & 0.008& 0.012 \\
        $E_{\bm \phi_h}$ & 0.002 & 0.003 \\
        \hline
    \end{tabular}
    \caption{Mean and max relative errors over $\bm \mu \in \mathcal{S}_{\text{cyl}}$ of fixed mesh ROM for $\bm u^0_h(\bm x)$ (top row) and mapped mesh ROM for $\tilde{\bm u}_h(\bm x)$ and $\bm \phi_h(\bm x)$ (bottom two rows) for hypersonic flow over a double wedge with no truncation ($N=4$).
}
\label{table:ramp}
\end{center}
\end{table}

\section{Conclusions}
\label{sec:conclusions}
We have applied an $r$-adaptive method using high-order solutions of the Monge-Ampère equation to the high-fidelity solution and model reduction of parametrized systems of PDEs with shocks. 
Using a high-order method for the FOM allows us to iteratively solve and refine on anisotropically adapted, curved grids that track relevant solution features. 
The usefulness of this approach was then demonstrated on parametrized supersonic and hypersonic flow problems with moving shocks. 
The approach does not depend on the discretization of the FOM and should work with continuous and discontinuous schemes, as long as they can solve on curved meshes. 

Future work includes the application of this method to problems with higher dimensional parameter spaces and more complicated phenomena, such as viscous effects or chemical reactions. 
The sensor in this work is not designed specifically to track shocks and should be able to refine towards other features. 
One challenge is that for viscous problems useful targets for refinement of the FOM, such as a boundary layer, may not be the features that need to be tracked in a model reduction process. 

In order to apply the Monge-Ampère equation to more challenging problems, its efficiency needs to be addressed. 
The Newton method from \cite{nguyen2023monge} should be more efficient than the fixed point method used here, but it requires the linearization of the solution-dependent sensor. 
Another target for efficiency gains is the evaluation of the sensor itself during the solution process, which is now handled with high-order interpolation. 
A more efficient approach would be to evaluate the sensor on a fine structured grid allowing for more rapid interpolation, as is done with the sensors in \cite{iollo2022mapping}. 

More study is needed for the solution of the Monge-Ampère equation on non-convex domains and domains with corners. 
In these cases, it might be preferable to use different PDEs for $r$-adaptivity such as a variational approach that directly controls node locations and element quality \cite{fortunato2016high}. 
While the Monge-Ampère equation is an effective method for mesh adaptivity with good guarantees on mesh quality, the described model reduction approach will work for any PDE-based $r$-adaptivity method. 

Finally, this approach can be coupled with intrusive ROMs and automated sampling procedures for greater accuracy and more robust performance when extrapolating. 

\section*{Acknowledgements} \label{}

We gratefully acknowledge the United States  Department of Energy under contract DE-NA0003965. 

This article has been authored by an employee of National Technology \& Engineering Solutions of Sandia,LLC under Contract No. DE-NA0003525 with the U.S. Department of Energy (DOE). The employee owns all right,title and interest in and to the article and is solely responsible for its contents. The United States Government retains and the publisher, by accepting the article for publication, acknowledges that the United States Government retains a non-exclusive, paid-up, irrevocable, world-wide license to publish or reproduce the published form of this article or allow others to do so, for United States Government purposes. The DOE will provide public access to these results of federally sponsored research in accordance with the DOE Public Access Plan https://www.energy.gov/downloads/doe-public-access-plan. SAND2023-10277O.

\bibliography{library}{}
\bibliographystyle{plain}
\end{document}